\documentclass[10pt,a4paper,reqno]{amsart}
\usepackage{amssymb}
\usepackage{amscd}
\usepackage{hyperref}
\numberwithin{equation}{section}

\theoremstyle{plain}

\newtheorem{theorem}[subsection]{Theorem}
\newtheorem{proposition}[subsection]{Proposition}
\newtheorem{lemma}[subsection]{Lemma}
\newtheorem{corollary}[subsection]{Corollary}
\newtheorem{remark}[subsection]{Remark}

\theoremstyle{definition}

\newtheorem{remarks}[subsection]{Remarks}
\newtheorem{example}[subsection]{Example}
\newtheorem{examples}[subsection]{Examples}

\renewcommand{\leq}{\leqslant}
\renewcommand{\geq}{\geqslant}

\newsavebox{\proofbox}
\savebox{\proofbox}{\begin{picture}(7,7)%
  \put(0,0){\framebox(7,7){}}\end{picture}}

\newcommand\Z{\mathbb{Z}}
\newcommand\R{\mathbb{R}}
\newcommand\T{\mathbb{T}}
\newcommand\C{\mathbb{C}}
\newcommand\N{\mathbb{N}}

\def\proof{\noindent\textit{Proof. }}

\def\remark{\noindent\textit{Remark. }}

\def\endproof{\hfill{\usebox{\proofbox}}}

\begin{document}

\title[GWP of Cubic NLS]{Global well-posedness of the cubic nonlinear Schr\"odinger equation on closed manifolds}
\author{Zaher Hani}

\address{Courant Institute for Mathematical Sciences, 251 Mercer Street, New York, NY 10012.}
\email{hani@cims.nyu.edu}

\maketitle

\begin{abstract}
We consider the defocusing cubic non-linear Schr\"odinger equation on general closed (compact without boundary) Riemannian surfaces. The problem was shown to be locally well-posed in $H^s(M)$ for $s>\frac{1}{2}$ in \cite{BGT}. Global well-posedness for $s\geq 1$ follows easily from conservation of energy and standard arguments. In this work, we extend the range of global well-posedness to $s>\frac{2}{3}$. This generalizes, without any loss in regularity, the results in \cite{B4}\cite{dSPST}, where the same result is proved for the torus $\T^2$. The proof relies on the I-method of Colliander, Keel, Staffilani, Takaoka, and Tao, a semi-classical bilinear Strichartz estimate proved by the author in \cite{H}, and spectral localization estimates for products of eigenfunctions, which is essential to develop multilinear spectral analysis on general compact manifolds. 
\end{abstract}
\tableofcontents

\section{Introduction}

We consider the defocusing cubic non-linear Schr\"odinger equation on a 2D-Riemannian manifold $(M,g_{\alpha \beta})$ given by:

\begin{eqnarray}\label{cubic NLS}
i\partial_tu +\Delta_g u &=|u|^2u\\
u(0,x)&=u_0(x)
\end{eqnarray}

where $t \in \R$, $x \in M$, and $u(t,x)$ is a complex-valued function on $\R\times M$. $\Delta_g$ is the Laplace-Beltrami operator on $M$ (negative operator) with respect to the metric $g_{\alpha \beta}$. This equation is an important model in several areas of physics, including laser optics, plasma physics, and Bose-Einstein condensates (see \cite{SS}). We will be most interested in the case when $M$ is closed, i.e. compact without boundary.

\subsection{Local theory:} The state of the art of the local existence question for $\eqref{cubic NLS}$ depends on the manifold $M$ under consideration. When $M=\R^2$, Strichartz estimates allow to prove local well-posedness in $H^s(\R^2)$ for $s\geq 0$ (see \cite{CW} or \cite{C}\cite{T} for a survey). The well-posedness for $s>0$ is described as \emph{subcritical}, i.e. the time of existence guaranteed by the local theory depends only on the size of the initial data in $H^s(\R^2)$. In contrast, the local well-posedness for $s=0$ is \emph{critical} in the sense that the time of existence depends on the profile of the data. This difference stems from the fact that $s=0$ is the critical regularity for which the $\dot H^s$ norm is left invariant by the scaling symmetry enjoyed by the equation:

\begin{equation}\label{scaling symmetry} 
u(t,x)\to \frac{1}{\lambda}u(\frac{t}{\lambda^2},\frac{x}{\lambda}).
\end{equation}

The range $s>0$ is the \emph{scaling subcritical} range and it determines the range of \emph{subcritical local well-posedness} as well. Moreover, the equation is known to be ill-posed in a certain sense for $s<0$ (the flow is not uniformly continuous in any neighborhood of 0) \cite{CCT2}. 

The first answer to the local existence question on a compact manifold was provided by Bourgain in the case of the torus \cite{B1}. In this paper, he derives periodic Strichartz estimates for the linear Schr\"odinger equation and uses them to prove subcritical local well-posedness for scaling-subcritical equations, including \eqref{cubic NLS} on $\T^2$ for $s>0$. The loss of $\epsilon$ derivatives (for arbitrary small $\epsilon$) in the critical Strichartz estimates forbids one from proving local well-posedness in $L^2$ by arguing as in the $\R^2$ case. This problem is still open. As on $\R^2$, the equation is known to be ill-posed for $s<0$ \cite{CCT1}\cite{CCT2}. As a result, one concludes that on $\T^2$, as on $\R^2$, the scaling subcritical range $s>0$ determines the range of subcritical local-wellposedness.

The situation changes if one considers more complex geometries for the compact manifold $M^2$. In fact, one expects that, in contrast to the case of the wave equation, the infinite speed of propagation of the linear Schr\"odinger equation would cause the geometry to play a role in the local theory. This is the case when $M$ is the 2-sphere $S^2$. Here, Burq, Gerard, and Tzvetkov proved that \eqref{cubic NLS} is subcritically locally well-posed for $s>\frac{1}{4}$\cite{BEE} and $C^3-$ill-posed (i.e. the flow is not $C^3$ at the origin) for $s<\frac{1}{4}$ cf. \cite{ BEE,BGT4}. The proof of well-posedness was based on a sharp bilinear Strichartz estimate which in turn is based on bilinear eigenfunction estimates and the sharp localization of the eigenvalues of the Laplace-Beltrami operator on $S^2$. 

For general compact manifolds without boundary, Burq, Gerard, and Tzvetkov established in \cite{BGT} (see also \cite{ST}) Strichartz estimates for the Schr\"odinger equation on closed manifolds with a loss of derivative. This allowed them to prove local well-posedness of $\eqref{cubic NLS}$ for $s>\frac{1}{2}$. It is expected, as in the case of $\R^2$, $\T^2$, and $S^2$, that for any compact manifold $M$ there exists a critical regularity $s_0(M)\geq 0$ ($s=0$ being the scaling critical regularity) for which one has subcritical well-posedness for $s>s_0(M)$ and ill-posedness for $s<s_0(M)$. The exact relation of this regularity $s_0(M)$ to the geometry of the manifold is a very interesting thing still to understand. 

We end our discussion of the local theory by mentioning that in the case when $M$ is a compact manifold with boundary, local existence has been proved for all $s>\frac{2}{3}$ in \cite{J} building on the work of Blair, Smith and Sogge \cite{BSS}. We should also mention the work of Ivanovici \cite{Ivanovici} on Strichartz estimates for compact manifolds with strictly concave boundary and exterior domains.

\subsection{Global theory:}
$\eqref{cubic NLS}$ enjoys the following two important conservation laws\footnote{There is also a conservation law for the momentum, but we will not be using it.}

\begin{itemize}
\item Conservation of mass
\begin{equation}\label{mass conservation}
M[u](t):= \int_{M} |u(t,x)|^2 dx
\end{equation}

\item Conservation of the Hamiltonian
\begin{equation} \label{energy conservation}
E[u](t):=\int_{M}\frac{1}{2} |\nabla_g{u}(t,x)|^2+\frac{1}{4}|u(t,x)|^4 dx=E[u](0)
\end{equation}

\end{itemize}

These conservation laws and the subcritical nature of the well-posedness for $s=1$ imply global well-posedness in the energy space. Standard arguments imply global wellposedness for all $s\geq 1$
(see \cite{T},\cite{B3}). The first result to prove local well-posedness of $\eqref{cubic NLS}$ below the energy norm was due to Bourgain in the case when $M=\R^2$(\cite{B3}\cite{B2}). Roughly speaking, the idea was to split the solution into a low frequency part supported on frequencies $\leq N$ ($N$ is a parameter chosen at the end depending on an arbitrarily chosen time interval $[0,T]$) and a high frequency one supported on frequencies larger than $N$. One then evolves the low frequency part by the nonlinear flow of $\eqref{cubic NLS}$(which is global in time since low frequencies have finite energy) and the high frequency part by the linear flow $e^{it\Delta}$(which conserves the $H^s$ norm). While the sum of those two flows is far from being a solution to $\eqref{cubic NLS}$, the difference between the real solution $u$ and this sum can be shown to be smoother, in fact in $H^1$, which means that the solution cannot blowup in $H^s$\footnote{The situation is actually a bit more complicated. One has to do this analysis on small intervals for technical reasons, which requires iterating the above-mentioned procedure to conclude that the interval of existence $[0,T]$ can be chosen arbitrarily.}

This ``high-low" Fourier Truncation method of Bourgain was the impetus to a more powerful method of Colliander, Keel, Staffilani, Takaoka, and Tao based on the almost conservation of a modified energy functional(see \cite{CKSTT1} for instance). Here one considers a Fourier multiplier $I$ which is the \textbf{I}dentity for low frequencies and an \textbf{I}ntegration by an order of $(1-s)$(enough to make an $H^s$ function in $H^1$) for high frequencies. While the energy of the modified solution $Iu$ is not conserved in general, it is almost conserved and this is enough to close an iteration argument and prove that $E[Iu]$ does not blow up, and hence $||u(t)||_{H^s}$, being controlled by the latter, does not blow up either. The success of this strategy, commonly referred to as the ``I-method" for obvious reasons, is based on multilinear analysis in Fourier space and known local theory.

In the context of compact manifolds, global well-posedness below energy norm was first considered by Bourgain in \cite{B4} where he used the language of normal forms. The problem was also studied in \cite{dSPST} where a language similar to that on $\R^2$ was used. The result in these papers is global well-posedness of $\eqref{cubic NLS}$ in $H^s(\T^2)$ for all $s>\frac{2}{3}$. The proof relies heavily on a Fourier space analysis, which is one of the many advantages of working on an abelian group like $\T^2$ or $\R^2$. We should also mention the recent paper \cite{Akahori} which proves global well-posedness on sphere-like manifolds (Zoll manifolds) in $H^s$ for all $s>\frac{15}{16}$ as well as some recent developments on the local and global theory for energy critical problems on compact, exterior, or product domains (\cite{IvanoviciPlanchon, HTT, HTT1, Herr, IP, IP2}).

In this paper, we consider the problem of global well-posedness below the energy norm for $\eqref{cubic NLS}$ posed on a general closed (compact without boundary) Riemannian 2-manifold $M$. The first thing one misses when moving to the setting of compact manifolds is the Fourier transform. This loss is only partially compensated by the spectral resolution of the Laplace-Beltrami operator. Since the $I$ operator is now given as a spectral multiplier rather than a Fourier multiplier as in the case of the torus $\T^2$ (or on $\R^2$), there are considerable difficulties in running the multilinear analysis in the former case. Firstly, in order to be able to perform a multilinear analysis with Littlewood-Paley pieces at fine scales, one has to know the spectral localization (which stands for the frequency support) of products of eigenfunctions. The resolution of this difficulty will be discussed below and leads to implications that we think are interesting in their own right. Secondly, after estimating individual products of eigenfunctions, one is faced with the problem of summing over these products without incurring additional losses in regularity\footnote{Any use of Weyl's law leads to a loss that cannot be tolerated.}. This particular difficulty stems from the absence of a good analogue to Fourier inversion which allows the treatment of multilinear Fourier multipliers on $\R^n$ or $\T^n$ as convolution operators in the Fourier space. To resolve this latter difficulty, we use a Fourier series decomposition trick (cf. \cite{Chr} \cite{Hw}). The basic idea is to convert certain multilinear multipliers satisfying Coifman-Meyer type symbol estimates into tensored multipliers for which the problem can be easily solved because the sums decouple. This process is done using Fourier series expansions (see Section \ref{multilinear spectral multiplier lemma}). 

Another main difficulty encountered in the case of compact manifolds $M$ in comparison to the torus, is the additional loss of derivatives in Strichartz estimates. The loss of $\frac{1}{4}$ of a derivative in linear ($\frac{1}{2}$in bilinear) Strichartz estimates not only leads to limiting local well-posedness to the range $s>\frac{1}{2}$ but also poses itself as a serious problem in establishing global well-posedness here as well. This should be compared to case of the torus where only $\epsilon=0+$ derivatives are lost in the linear and bilinear estimates (see \cite{B1},\cite{B4},\cite{dSPST}) and to the slightly more suggestive case of the sphere $S^2$ where there is a necessary loss of $\frac{1}{4}$ of a derivative in the bilinear Strichartz estimate (\cite{BEE}). Overcoming this difficulty is done by using the semi-classical bilinear Strichartz estimates proved by the author in \cite{H}. In this latter paper, it is shown that at the semi-classical level, one can obtain bilinear improvement to linear Strichartz estimates similar (actually the same) as those refinements proved by Bourgain \cite{B2} on $\R^2$. While these short range refinements do not improve on the loss of half a derivative in the bilinear estimates on $[0,1]\times M$ (at least without a smarter way of summing over the short range pieces of the interval $[0,1]$), they \emph{do} offer crucial improvements for bilinear estimates on the rescaled (inflated) manifold $\lambda M$ for $\lambda>1$. By moving the analysis to $\lambda M$ rather than $M$ itself (for $\lambda$ growing large with the asymptotic parameter $N$), one is able to capture this gain and compensate (almost completely) the loss of derivatives on $M$ itself in comparison to the case of the torus $\T^2$.
 
After dealing with these new difficulties, we are able to run the I-method strategy and generalize without any loss in regularity the results proved in \cite{B4}\cite{dSPST} for $\T^2$:

\begin{theorem}\label{main theorem}
Let $M$ be any compact Riemannian 2-manifold without boundary. $\eqref{cubic NLS}$ is globally well-posed in $H^s$ for all $s>\frac{2}{3}$. Moreover, the $H^s$ norm of the solution satisfies the following polynomial bound:

\begin{equation}\label{polynomial bound}
||u(t)||_{H^s(M)}\lesssim_{||u_0||_{H^s}} t^{\frac{2s(1-s)}{3s-2}+}.
\end{equation}
\end{theorem}

\subsection{Bilinear Strichartz estimates:}\label{bilinear Strichartz estimates}
We recall the bilinear Strichartz estimate on $\R^2$ that was proved in \cite{B2} (and was first applied there in proving the above-mentioned ``high-low" Fourier truncation method): Suppose that $u_0,v_0 \in L^2(\R^2)$ are frequency localized at dyadic scales $N_1$ and $N_2$ respectively, i.e. $\operatorname{supp} \widehat{u_0}\subset \{\xi:|\xi|\sim N_1\}$ and $\operatorname{supp} \widehat{v_0}\subset \{\xi:|\xi|\sim N_2\}$ with $N_2 \leq N_1$. Then

\begin{equation}\label{bilinear Strichartz on R^2}
||e^{it\Delta}u_0 e^{it\Delta}v_0||_{L^2_{t,x}(\R\times \R^2)}\lesssim \left(\frac{N_2}{N_1}\right)^{1/2}||u_0||_{L^2(\R^2)}||v_0||_{L^2(\R^2)}.
\end{equation}

By a simple scaling and limiting argument, it is easy to see that $\eqref{bilinear Strichartz on R^2}$ is equivalent to the short range inequality that is the same as $\eqref{bilinear Strichartz on R^2}$ but with $L^2_{t,x}([0,N_1^{-1}]\times \R^2)$ replacing $L_{t,x}^2(\R\times \R^2)$. In \cite{H}, it is proved that this short range estimate holds true on any compact Riemannian manifold without boundary. The exact statement is as follows (From here on, we use the notation $\Delta$ to denote the Laplace-Beltrami operator $\Delta_g$):

\begin{theorem}\label{semiclassical bilinear Strichartz estimate}\cite{H}
Let $M$ be a compact Riemannian 2-manifold without boundary. Suppose that $u_0, v_0 \in L^2(M)$ are spectrally localized on dyadic scales $N_1$ and $N_2$ respectively, i.e. $u_0=\mathbf{1}_{[N_1,2N_1)}(\sqrt{-\Delta})u_0$ and $v_0=\mathbf{1}_{[N_2,2N_2)}(\sqrt{-\Delta})v_0$ with $N_2 \leq N_1$. Then

\begin{equation}\label{semiclassical bilinear estimate}
\|e^{it\Delta} u_0 e^{it\Delta}v_0\|_{L^2([0,\frac{1}{N_1}]\times M)} \lesssim \left(\frac{N_2}{N_1}\right)^{1/2}\|u_0\|_{L^2(M)}\|v_0\|_{L^2(M)}.
\end{equation}

The same estimate holds if $v_0$ is supported at frequencies $\lesssim N_2$.
\end{theorem}

Since we won't prove this estimate in this paper, we remark that the numerology in $\eqref{semiclassical bilinear estimate}$ can be understood (heuristically at least) by a simple back-of-the-envelope calculation. Thinking of $e^{it\Delta}u_0$ as a ``bump function" localized in frequency at scale $N_1$ and \emph{initially} (at $t=0$) localized in space at scale $\frac{1}{N_1}$. The Schr\"odinger evolution moves this bump function at a speed $N_1$ thus expanding its support at this rate while keeping the $L^2$ norm conserved. Similarly, $e^{it\Delta}v_0$ could be thought of as a ``bump function" that is initially concentrated in space at scale $\sim \frac{1}{N_2}$ and moving (expanding) at speed $N_2$. A simple schematic diagram allows to estimate the space-time overlap of the two expanding ``bump functions" thus giving the estimate $\frac{N_2^{(d-1)/2}}{N_1^{1/2}}$ for the $L^2_{t,x}([0,N_1^{-1}]\times \R^d)$ of the product.

By splitting the time interval into pieces of length $N_1^{-1}$, one can use $\eqref{semiclassical bilinear estimate}$ to be get a bilinear Strichartz estimate over any time interval $[0,T]$. We will not be interested in the case of $T=1$ (which is relevant for the local theory and for which the gain from the factor $\left(\frac{N_2}{N_1}\right)^{1/2}$ is lost in the summation\footnote{Unless one is able to find a smarter way to sum the contribution of each individual subinterval.}(see \eqref{Lambda T})), but rather we will be concerned with estimates on small intervals $[0,T]$ with $T\ll 1$ (actually $N_1^{-1}\leq T \ll N_2^{-1}$ in most important cases). More precisely, our main use of Theorem \ref{semiclassical bilinear Strichartz estimate} will be to overcome the difficulty of derivative loss in the bilinear Strichartz estimates on general closed manifolds. The idea to recapture the loss of derivatives is to ``inflate" the manifold $M$ and move the analysis to a dilated manifold $M_\lambda:=\lambda M$ with its rescaled metric $g_\lambda$. If $\lambda$ is large, one expects frequency-localized Strichartz estimates to become better behaved and exhibit a gain manifested by a negative power of $\lambda$. This can be heuristically explained by the fact that, in the limit $\lambda \to \infty$, $M_\lambda$ becomes flat and hence the Strichartz estimates should resemble those on $\R^2$ (where no derivative loss occurs). Theorem \ref{semiclassical bilinear Strichartz estimate} yields bilinear Stichartz estimates on $[0,T]\times M$ with $(TN_2)^{1/2}$ on the R.H.S. of $\eqref{semiclassical bilinear estimate}$. By moving the analysis to the rescaled manifold $M_\lambda=\lambda M$ (with $\lambda\gg 1$ chosen appropriately as dictated by a scaling argument in the I-method), this translates to a gain of $\lambda^{-1/2}$ in the bilinear Strichartz estimates on $[0,1]\times M_\lambda$ (cf. $\eqref{the N_2/lambda decay}$), a fact which allows to compensate completely for the extra loss of derivative on closed manifolds as compared to that on the torus $\T^2$. In short, at the scale relevant to the I-method, the Strichartz estimates deduced from Theorem \ref{semiclassical bilinear Strichartz estimate} above are essentially equally ``powerful" to those on the torus.

\subsection{Spectral localization of products of eigenfunctions:} The I-method is based on a multilinear analysis of a multilinear spectral multiplier. Typically, in such analysis one splits all functions into Littlewood-Paley pieces and analyses each generic ``product" by itself. One of the main reasons why frequency localization is such a useful tool in analyzing multilinear operators (ranging from products of functions $fg$ to more general multilinear multipliers) on $\R^d$(or more generally on any abelian group) is the fact that if $f$ is frequency localized in $\{|\xi|\in[N,2N)\}$ and $g$ is frequency localized in $\{|\xi|\in [M,2M)\}$ with (say) $N\geq M$, then the product $fg$ is frequency localized in the region $\{|\xi|\in [2N+2M,N-2M]$. Ultimately, this is due to the fact that $e^{i\xi_1.x}e^{i\xi_2.x}=e^{i(\xi_1+\xi_2).x}$ and the latter is also a character at frequency $\xi_1+\xi_2$. The same thing is true on the torus $\T^d$ with $\xi_1,\xi_2 \in \Z^d$ in which case the product of the two eigenfunctions $e^{i\xi_1.x}$ and $e^{i\xi_2x}$ is \emph{also} an eigenfunction.

Once we move to the case of general compact manifolds\footnote{The results here apply to general compact manifolds $M$ possibly with boundary and the eigenfunctions considered satisfy either Dirichlet or Neumann boundary conditions.}, this sharp localization is not as straightforward. Consequently, the following question becomes of importance: given two eigenfunctions $f(x)$ and $g(x)$ of the Laplace-Beltrami operator $-\Delta$ on a compact $d-$manifold $M^d$ with eigenvalues $\mu^2$ and $\nu^2$ respectively ($\mu,\nu \in \R^+$), what can be said about the spectral localization of the product function $f(x)g(x)$, in other words what can be said about $\pi_{\eta^2}(fg)$ where $\pi_{\eta^2} $ is the projection on the $\eta^2-$eigenspace of $-\Delta$? For definiteness, we call this problem that of spectral localization of products of eigenfunctions.

In the range when $\eta \gg \max(\mu,\nu)$, an answer can be given using the parametrix representation of eigenfunctions and crude integration by parts. Essentially one obtains that $||\pi_{\eta^2}(fg)||_{L^2(M)}\lesssim_p (\eta)^{-p}$ for any $p$ (see \cite{BEE} and section \ref{part 1}). However, this gives a very large range of localization of eigenvalues for the product $fg$, namely all eigenvalues $\eta^2$ satisfying $\eta \leq C\max(\mu,\nu)$ for some large constant $C$. 

Studying the indicative case of the sphere suggests that the product should be localized at a much smaller range similar to that on the torus. On the sphere $S^d$ the eigenfunctions have a special form: they are given by spherical harmonics. Multiplying a spherical harmonic of degree $k$ with another of degree of degree $l$, one can expand the product in terms of spherical harmonics of degree $\leq k+l$. The eigenvalue associated with a spherical harmonic of degree $k$ is $k(k+d-1)$ where $d$ is the dimension of the sphere. This suggests that the spectrum of the product of two eigenfunctions $f$ and $g$ with eigenvalues $\mu^2$ and $\nu^2$ respectively, and satisfying $\mu \geq \nu$ should be sharply concentrated around the range $\sqrt{-\Delta}\in[\mu+\nu,\mu-\nu]$, in the sense that $\pi_{\eta^2}fg$ should decay when $\eta \geq \mu+\nu$. 

Noticing that $\pi_{\eta^2}(fg)$ is captured by the correlation integral $\int_M h(x)f(x) g(x) dx$ where $h$ is an eigenfunction of $-\Delta$ with eigenvalue $\eta^2$, the problem is easily seen to be equivalent to obtaining decay estimates for such an integral. The most general result will be given in Theorem \ref{A_0 and A_n theorem} and the remarks following it (where the more general problem of the spectral localization of a product of multiple eigenfunctions is studied). Theorem \ref{A_0 and A_n theorem} provides an \emph{identity} satisfied by the correlation integral (see equation $\eqref{A_0 and A_n}$). As a consequence of this identity one gets that if $\eta=\mu+K\nu$ with $K>1$ (in other words if $K:=\frac{\eta-\mu}{ \nu}>1$) then the integral $\int_M h(x)f(x) g(x)dx$ decays like $K^{-J}$ for any integer $J$ (see Theorem \ref{A_0 and A_n theorem} and remarks thereafter). The identity is proved using the Ricci commutation identities for covariant derivatives acting on tensors\footnote{All needed notions from Riemannian geometry are revised in section \ref{part 2} for completeness.} and an iteration argument. The proof is presented in section 7.

We leave the statement of Theorem \ref{A_0 and A_n theorem} to section \ref{part 1}. Here we cite a particular example of the estimates that can be deduced from identity $\eqref{A_0 and A_n}$. We note that the exact form of such estimates depends on the norms in which one would like to estimate the eigenfunctions or Littlewood-Paley pieces. For example, if one would like to estimate all eigenfunctions in $L^2$ (thus overruling the exponent distribution dictated by H\"older's inequality), we have the following corollary that can be understood as a refinement of the bilinear eigenfunction (bilinear Sogge) estimates proved by Burq, Gerard, and Tzvetkov in \cite{BEE, MEE}.

\begin{corollary}\label{spectral cluster lemma 1}
Let $M^d$ be a compact $d-$manifold possibly with boundary. Suppose that $\nu^2, \lambda^2, \mu^2$ are eigenvalues of the operator $-\Delta$ with Dirichlet or Neumann boundary conditions satisfying $\nu \geq \lambda \geq \mu$. If $\nu=\lambda+K\mu+2$ \footnote{the additive factor of 2 is used purely for technical reasons.}for some $K>1$ (i.e. $K=\frac{\nu-\lambda-2}{\mu}>1$), then for any $f,g \in L^2(M)$ and any $J\in \N$:

\begin{equation}\label{spectral cluster lemma 1 equation}
\pi_{\nu}\left(\mathbf{1}_{[\lambda,\lambda+1]}(\sqrt{-\Delta})f\mathbf{1}_{[\mu,\mu+1]}(\sqrt{-\Delta})g\right) \lesssim_J \frac{\Lambda(d,\mu)}{K^J}\|f\|_{L^2(M)}\|g\|_{L^2(M)}
\end{equation}

where
\begin{equation}\label{def of Lambda(d,mu)}
\Lambda(d,\mu):= 
\begin{cases}
\mu^{1/2} &\text{if } d=2\\
\mu^{1/2}(\log\mu)^{1/2} &\text{if } d=3\\
\mu^{\frac{d-2}{2}} &\text{if } d\geq 4
\end{cases}.
\end{equation}

More generally, one has:
\begin{equation}
\int_M \mathbf{1}_{[\nu,\nu+1]}(\sqrt{-\Delta})h\,\mathbf{1}_{[\lambda,\lambda+1]}(\sqrt{-\Delta})f\,\mathbf{1}_{[\mu,\mu+1]}(\sqrt{-\Delta})g \,dx\lesssim_J \frac{\Lambda(d,\mu)}{K^J}\|h\|_{L^2(M)}\|f\|_{L^2(M)}\|g\|_{L^2(M)}.
\end{equation}

Interchanging the roles of $\lambda$ and $\nu$, the same estimates hold if $\nu <\lambda$ and $\nu=\lambda-K\mu -2$ with $K>1$ (i.e. if $K:=\frac{\lambda-\nu-2}{\mu}>1$).
\end{corollary}

It is worth mentioning that the factor of $\Lambda(d,\mu)$ appears because we choose to estimate all eigenfunctions in $L^2$ and use a variant of the bilinear eigenfunction estimates of Burq, Gerard, and Tzvetkov\cite{MEE}. In particular, the identity we establish in Theorem \ref{A_0 and A_n theorem} does not exhibit any loss of derivatives, a fact which is often essential for performing an efficient analysis of the spectral localization of products of Littlewood-Paley pieces as well as products of eigenfunction clusters (see subsection \ref{decay of modified energy section} where we apply Theorem \ref{A_0 and A_n theorem} to bound products of Littlewood-Paley pieces in a way that crucially does not lose any derivatives).

The decay of $K^{-J}$ mentioned above and guaranteed by Theorem \ref{A_0 and A_n theorem} is enough to conduct the usual fine scale-course scale Littlewood-Paley interactions on a general manifold (even with boundary) similar to how it is done on $\R^n$ and $\T^n$ once one is able to overcome an additional difficulty mentioned previously and treated in section \ref{multilinear spectral multiplier lemma}. In particular, it allows us to bound the notorious low-high frequency interaction in the I-method (see subsection \ref{decay of modified energy section})

The paper is organized as follows: in section \ref{preliminaries} we fix the notation and explain our rescaling of the manifold $M$ and its effects on eigenvalues, eigenfunctions, and linear Strichartz  and $X^{s,b}$ estimates which we recall as well. In section \ref{bilinear Strichartz estimates}, we recall the bilinear Strichartz estimates proved in \cite{H} and show that the estimates on the time interval $[0,T]$ imply refined estimates on the rescaled manifold $\lambda M$. We also prove a version of these estimates involving $X^{s,b}$ spaces and with differential operators applied to the propagator $e^{it\Delta}u_0$. Such estimates will be needed after we apply the spectral localization machinery mentioned above and explained further in section \ref{part 1}. The proofs of the theorems in this latter section require recalling some ideas from Riemannian geometry and are delayed to the last section \ref{part 2} in order not to distract the reader. In section \ref{multilinear spectral multiplier lemma}, we deal with the problem of bounding multilinear spectral multipliers and formulate and prove a lemma pertaining to multipliers obeying Coifman-Meyer type estimates. In section \ref{I-method section}, we run the I-method strategy and prove Theorem \ref{main theorem}. Finally, in section \ref{part 2}, we review the Ricci commutation identities and present the proofs of section \ref{part 1}.

\textbf{Acknowledgements:} The author is deeply grateful to his advisor, Prof. Terence Tao, for his invaluable help, encouragement, and guidance. He is also thankful to the anonymous referees for their careful reading of the manuscript and their constructive comments and suggestions that considerably improved the presentation.

\section{Preliminaries}\label{preliminaries}

\subsection{Notation:} 

$M$ denotes a closed (compact without boundary) $C^\infty$ Riemannian 2-manifold (surface). In some sections (namely sections \ref{part 1}, \ref{multilinear spectral multiplier lemma}, and \ref{part 2}) the analysis applies to general $d-$ Riemannian manifolds possibly with boundary and $M$ will also be used to denote such a manifold. We let $g_{\alpha \beta}$ denote the Riemannian metric, $g^{\alpha \beta}$ its inverse, $\nabla$ the induced (metric compatible torsion-free) connection, and $\Delta=\nabla_\alpha \nabla^\alpha$ the connection Laplacian. 

Throughout the paper, we use the notation $X \lesssim Y$ to denote $X\leq CY$ for some constant $C>0$ and $A\sim B$ to denote $A\lesssim B \lesssim A$. All implicit constants are allowed to depend on $M$ and on its dimension (in sections where we consider $d-$manifolds). Often times we will attach a subscript to $\lesssim$ as in $\lesssim_\lambda$ to denote the possible dependence of the implicit constant $C$ on $\lambda$. We also use the notation $J+$ when $J\in \R$ to denote $J+\epsilon$ for a fixed arbitrarily small positive number $\epsilon$. Similarly, $J-$ refers to $J-\epsilon$. For functions $f\in C_0^\infty(\R\times M)$ we denote by $\widehat f (.,x)=\mathcal{F}_t f(.,x)$ the Fourier transform in time of the function $t \mapsto f(t,x)$ given by $\widehat f(\tau,x)=\int_\R e^{-it\tau}f(t,x)dt$. We will often take the liberty of omitting the $2\pi$ factors in the definition of the inverse transform as these constants are inconsequential in the analysis. 

The Laplace-Beltrami operator on a compact manifold is a non-positive self adjoint operator with compact resolvent\footnote{In section \ref{part 1} and \ref{part 2}, we consider the Laplace-Beltrami operators on a compact $d-$manifold $X$ with Dirichlet or Neumann boundary conditions. This also is an operator with compact resolvent on $L^2$.}. This gives an orthonormal basis of eigenfunctions $\{e_k\}_{k=1}^\infty$ corresponding to a non-decreasing sequence of non-negative eigenvalues $\mu_k$ of $-\Delta$. We will denote by $\nu_k$ the strictly increasing sequence of eigenvalues and by $\pi_{\nu_k}$ the orthogonal projection on the $\nu_k$ eigenspace. Often times (especially when we study the spectral localization of the operator $\sqrt{-\Delta}$), it will be convenient to denote $\nu_k=n_k^2$ where $n_k \in \R$ and $\pi_{n_k}=\pi_{\nu_k}$. For any interval in $I\subset [0,\infty)$ we denote by $P_I$ the orthogonal projection operator onto eigenvalues $\nu_k$ with $\sqrt{\nu_k}\in I$, or equivalently $P_{I}=\mathbf{1}_{I}(\sqrt{-\Delta})$. In most cases, we will be interested in the case when $I$ is a dyadic interval of the form $[N,2N)$ where $N\in 2^{\N}$ in which case we denote $P_N=P_{[N,2N)}$ for $N=2^{j}$ with $j=1,2,...$ and $P_1=P_{[0,2)}$. $H^s(M)$ is the natural Sobolev space associated with $(\operatorname{Id}-\Delta)^{1/2}$ with the norm:

$$
||u||_{H^s(M)}=\left(\sum_k \langle \nu_k\rangle ^{s}||\pi_{\nu_k} u||_{L^2(M)}^2\right)^{1/2}\sim \left( \sum_{N\in 2^{\N}} N^{2s}||P_N u||_{L^2(M)}\right)^{1/2}.
$$

We also define the space $X^{s,b}(\R \times M)$ as the completion of $C_0^\infty(\R \times M)$ under the norm:

\begin{align*}
||u||_{X^{s,b}(\R \times M)}=&\left(\sum_k \int_{\R_\tau} \langle \tau -\nu_k\rangle^{2b} \langle \nu_k\rangle^s||\widehat{\pi_{\nu_k} u}(\tau)||_{L^2(M)}\right)^{1/2}\\
=&||e^{-it\Delta}u||_{H_t^bH_x^s(\R\times M)}
\end{align*}
where $H_t^bH_x^s(\R \times M)=H_t^b(\R_t;H_x^s(M))$ (similarly we use mixed Lebesgue spaces $L_t^qL_x^r(\R \times M)= L_t^q\left(\R;L_x^r(M)\right)$.

For every compact interval  $I \subset \R$, we define the restriction space $X^{s,b}(I\times M)$ as the space of functions $u$ on $I\times M$ that admit extensions to $\R \times M$ in $X^{s,b}(\R \times M)$. $X^{s,b}(I\times M)$ is equipped with the restriction norm:

$$
||u||_{X^{s,b}(I\times M)}=\inf_{U \in X^{s,b}(\R \times M)}\{||U||_{X^{s,b}(\R \times M)} \text{ with $U=u$ on $I$}\}.
$$

\subsection{Rescaling $M$:} Let $M$ be a $C^\infty$ closed Riemannian surface. Any such surface can be thought of as being embedded in some ambient space $\R^L$. Using this embedding one can define a $\lambda$-rescaled version of $M$ for any $\lambda>0$, which we denote by $\lambda M=:M_\lambda$ and is given by $M_\lambda=D_\lambda M$ where $D_\lambda$ is the dilation by $\lambda$ in $R^L$: $x \mapsto \lambda x$. $M_\lambda$ inherits from $\R^L$ (with its Euclidean metric) a metric $g_{\alpha,\beta}^\lambda$.\footnote{From the differential topology point of view, the two manifolds $M$ and $M_\lambda$ are the same. However, they are different from the point of view of Riemannian geometry. In fact, the rescaled manifold $(M_\lambda,g_{\alpha \beta}^\lambda)$ is isometric to the Riemannian manifold $(M,\lambda^2g_{\alpha \beta})$. The dilation map $D_\lambda: M \to M_\lambda$ is a diffeomorphism from the manifold $M$ to the manifold $M_\lambda$. It is crucial to note that the metric on $M_\lambda$ is not equal to $(D_\lambda)_*g_{\alpha,\beta}$ but rather satisfies the relationship ${D_\lambda}_* g_{\alpha,\beta}=\frac{1}{\lambda^2}g^\lambda_{\alpha,\beta}$. The sample situation to keep in mind here is that of rescaling the interval $[0,1]$ to $[0,\lambda]$.}. When confusion might arise, we will distinguish tensors and operators on $M_\lambda$ with a $\lambda$ subscript or superscript, e.g. $g^\lambda_{\alpha, \beta}, \Delta_\lambda$, etc... For a function $f:M_\lambda \to \C$, we will often denote $\tilde f:= D_\lambda^* f$ the pull-back function given for $y \in M$ by $\tilde f(y) =f(D_\lambda y)=:f(\lambda y)$. It is easy to see that $||f||_{L^p(M_\lambda)}=\lambda^{\frac{2}{p}}||\tilde f||_{L^p(M)}$ and $\operatorname{Volume}(M_\lambda)=\lambda^2\operatorname{Volume}(M)$. 

Since $\Delta_\lambda f (x)=\frac{1}{\lambda^2}\Delta \tilde f(D_\lambda^{-1}x)=\frac{1}{\lambda^2}(\Delta \tilde f) (\frac{x}{\lambda})$\footnote{This can either be seen from formula for the Laplacian in local coordinate given by $\frac{1}{\sqrt{|g^\lambda|}}\partial_i \sqrt{|g^\lambda|}{g^{\lambda}}^{ij}\partial_j=\frac{1}{\lambda^2}\frac{1}{\sqrt{|g|}}\partial_i \sqrt{|g|}g^{ij}\partial_j$ since ${g^\lambda}^{ij}=\frac{1}{\lambda^2}g^{ij}$ or by noting that $g$ and $\lambda^2 g$ induce the same connection $\nabla$ and hence $\Delta_\lambda=\nabla^\alpha \nabla_\alpha ={g^\lambda}^{\alpha \beta}\nabla_{\beta}\nabla_{\alpha}=\frac{1}{\lambda^2}g^{\alpha \beta}\nabla_\beta \nabla_\alpha$.} for any $C^2$ function $f$ on $M_\lambda$, we get that the functions $\frac{1}{\lambda}e_k(\frac{x}{\lambda})$ form an orthonormal basis of $L^2(M_\lambda)$ with corresponding eigenvalues $\frac{\mu_k}{\lambda^2}$. 
As a result, the orthogonal spectral projection operator $\pi_{(\nu / \lambda^2)}$ on $M_\lambda$ is related to $\pi_{\nu}$ on $M$ by the relation

\begin{equation}\label{rescaling projections}
\pi_{(\nu/\lambda^2)} f(x) =\pi_{\nu} \tilde f \circ D_\lambda^{-1}=\left(\pi_{\nu} \tilde f\right)(\frac{x}{\lambda}).
\end{equation}

\subsection{Linear Strichartz estimates}\label{linear Strichartz estimates}

We state some of the linear Strichartz estimates needed. These were obtained in \cite{BGT} (see also \cite{ST}). We say that a pair of exponents $(q,r)$ is $2d-$Schr\"odinger admissible if $2<q\leq \infty$, $2\leq r <\infty$ satisfy $\frac{1}{q}+\frac{1}{r}=\frac{1}{2}$. For any two admissible pairs $(q,r)$, the following estimates hold for the linear propagator $e^{it\Delta}u_0$

\begin{equation}\label{linear Strichartz on M}
||e^{it\Delta}u_0||_{L_t^qL_x^r([0,1]\times M)}\lesssim_{q,r,M} ||u_0||_{H^{1\over q}(M)}
\end{equation}
(cf. \cite{BGT}). This estimate is derived from the following semi-classical Strichartz estimate: Suppose that $u_0=P_{[N,2N)}u_0$, then 

\begin{equation}\label{semiclassical linear Strichartz estimate}
||e^{it\Delta}u_0||_{L_t^qL_x^r([0,N^{-1}]\times M)}\lesssim ||u_0||_{L^2_x(M)}.
\end{equation}

Notice that this estimate has the same form as the corresponding estimate that holds on $\R^2$ (where it holds on $\R\times \R^2$). This theme will also show up when we treat bilinear Strichartz estimates in the next section. $\eqref{linear Strichartz on M}$ follows by splitting the interval $[0,1]$ into $N$ pieces of length $1\over N$ each and applying $\eqref{semiclassical linear Strichartz estimate}$ on each of them using the conservation of mass (a type of square function estimate is also needed, see \cite{BGT}). 

Most of the analysis we will do is on the rescaled manifold $M_\lambda$. This requires us to know the dependence on $\lambda$ of the implicit constant in the linear (and bilinear) Strichartz estimates. This is possible by scaling and using the semi-classical estimate $\eqref{semiclassical linear Strichartz estimate}$ rather than $\eqref{linear Strichartz on M}$. In fact, we will see here and in the next section that time $1$ Strichartz estimates on $M_\lambda$ follow from time $\frac{1}{\lambda^2}$ Strichartz estimates on $M$, which are obtained using $\eqref{semiclassical linear Strichartz estimate}$ by splitting the time interval into pieces of relevant lengths.

If $u_0 \in L^2(M_\lambda)$ is spectrally localized on the interval $[N,2N)$, by which we mean that $u_0=P_{N}u_0$, then the pull-back function $\tilde u_0 \in L^2(M)$ is spectrally localized on the interval $[\lambda N,2\lambda N)$. As a result, we calculate:

\begin{align*}
||e^{it\Delta_\lambda}u_0||_{L_t^qL_x^r([0,1]\times M_\lambda)}=\lambda^{2\over q} ||e^{it\lambda^2\Delta_\lambda}u_0||_{L_t^qL_x^r([0,\frac{1}{\lambda^2}]\times M_\lambda)}=\lambda^{\frac{2}{q}+\frac{2}{r}}||e^{it\Delta}\tilde u_0||_{L_t^qL_x^r([0,\frac{1}{\lambda^2}]\times M)}
\end{align*}

since $\lambda^2 \Delta_\lambda u_0(x) =(\Delta \tilde u_0)(\frac{x}{\lambda})$. If the length of the time interval $\frac{1}{\lambda^2} \leq \frac{1}{\lambda N}$, then $\eqref{semiclassical linear Strichartz estimate}$ gives the bound (recall that $\frac{2}{q}+\frac{2}{r}=1$):

$$
||e^{it\Delta_\lambda}u_0||_{L_t^qL_x^r([0,1]\times M_\lambda)}\lesssim \lambda ||e^{it\Delta}\tilde u_0||_{L_t^qL_x^r([0,\frac{1}{\lambda N}]\times M)}\lesssim \lambda ||\tilde u_0||_{L^2(M)}=||u_0||_{L^2(M_\lambda)}.
$$

If, on the other hand, $\frac{1}{\lambda^2}\geq \frac{1}{\lambda N}$, we split the interval $[0,\frac{1}{\lambda^2}]$ into $\frac{N}{\lambda}$ pieces of length $\frac{1}{\lambda N}$ and apply $\eqref{semiclassical linear Strichartz estimate}$ on each of them to get:

$$
||e^{it\Delta_\lambda}u_0||_{L_t^qL_x^r([0,1]\times M_\lambda)}\lesssim \left(\frac{N}{\lambda}\right)^{1/q}\lambda ||\tilde u_0||_{L^2(M)}=\left(\frac{N}{\lambda}\right)^{1/q}||u_0||_{L^2(M_\lambda)}.
$$

As a result, we get that:

\begin{equation}\label{linear Strichartz estimate on M_lambda}
||e^{it\Delta}u_0||_{L_t^qL_x^r([0,1]\times M_\lambda)} \lesssim\left\{
    \begin{array}{ll}
        ||u_0||_{L^2(M)} & \mbox{if } \lambda \geq N \\
       \left(\frac{N}{\lambda}\right)^{1/q}||u_0||_{L^2(M)} & \mbox{if } \lambda \leq N
    \end{array}
\right.
\end{equation}

Notice that in the limit $\lambda \to \infty$, $M_\lambda$ becomes flat and the Strichartz estimates at a fixed frequency scale become the same as those satisfied on $\R^2$. 

We finally cite an additional Strichartz estimate that will be of use to us. Combining Bernstein's inequality for spectrally localized functions (which is true on compact manifolds, cf. Corollary 2.2 of \cite{BGT}) with the $L_t^8L_x^{8\over 3}$ semiclassical estimate in $\eqref{semiclassical linear Strichartz estimate}$ one gets for $\tilde u \in L^2(M)$ spectrally localized in the interval $[N,2N)$:

$$
||e^{it\Delta}\tilde u_0||_{L_{t,x}^8([0,\frac{1}{N}]\times M)}\lesssim N^{2/4}||e^{it\Delta}\tilde u_0||_{L_t^8L_x^{8/3}([0,\frac{1}{N}]\times M)}\lesssim N^{1/2}||\tilde u_0||_{L_x^2(M)}.
$$

Applying the same rescaling argument as the one used to get $\eqref{linear Strichartz estimate on M_lambda}$ we get:

\begin{equation}\label{L^8 estimate}
||e^{it\Delta}u_0||_{L_{t,x}^8([0,1]\times M_\lambda)} \lesssim\left\{
    \begin{array}{ll}
        N^{1/2}||u_0||_{L^2(M)} & \mbox{if } \lambda \geq N \\
       \left(\frac{N}{\lambda}\right)^{1/8}N^{1/2}||u_0||_{L^2(M)} & \mbox{if } \lambda \leq N
    \end{array}
\right.
\end{equation}

It is well known that any estimate of the form $||e^{it\Delta}u_0||_Y\lesssim ||u_0||_{L_x^2}$ where $Y$ is a Banach space of space-time functions satisfying $||e^{it\theta}f||_Y \lesssim ||f||_Y$, translates directly into an embedding of $X^{0,1/2+}$ into $Y$ (cf \cite{G}, lemma 2.9 of \cite{T}, or the next section for a bilinear version of this). As a result of this, we get the following consequences of $\eqref{linear Strichartz estimate on M_lambda}$ and $\eqref{L^8 estimate}$:

\begin{equation}\label{L^4 estimate Xsb}
||u||_{L_{t,x}^4([0,1]\times M_\lambda)}\lesssim \left\{
    \begin{array}{ll}        
    ||u||_{X^{0,1/2+}([0,1]\times M)} & \mbox{if } \lambda \geq N \\
    \left(\frac{N}{\lambda}\right)^{1/4}||u||_{X^{0,1/2+}([0,1]\times M)} & \mbox{if } \lambda \leq N
    \end{array} 
    \right.
\end{equation}

\begin{equation} \label{L^8 estimate Xsb}
||u||_{L_{t,x}^8([0,1]\times M_\lambda)}\lesssim \left\{
    \begin{array}{ll}

        N^{1/2}||u||_{X^{0,1/2+}([0,1]\times M)} & \mbox{if } \lambda \geq N \\
       \left(\frac{N}{\lambda}\right)^{1/8}N^{1/2}||u||_{X^{0,1/2+}([0,1]\times M)} & \mbox{if } \lambda \leq N
    \end{array}
\right.
\end{equation}

whenever $u(t)=P_N u(t)$. 

\section{Bilinear Strichartz estimates}

In this section, we use the short-range/semi-classical bilinear Strichartz estimate $\eqref{semiclassical bilinear estimate}$ proved in \cite{H} to derive bilinear estimates on $[0,T]\times M$ for any $T>0$. These estimates over the interval $[0,T]$ translate by rescaling to bilinear estimates on $[0,1]\times M_\lambda$. Since the spectral localization machinery of section \ref{part 1} will eventually require us to bound products of the form $P(D)e^{it\Delta}u_0Q(D)e^{it\Delta}v_0$ for differential operators $P
(D)$ and $Q(D)$, we will have to obtain estimates for such products as well. Luckily, this follows directly from the parametrix representation of $e^{it\Delta}$ and the results in \cite{H}. We also translate those bilinear Strichartz estimates into bilinear $X^{s,b}$ estimates as was done in section \ref{linear Strichartz estimates}.

The estimate $\eqref{semiclassical bilinear estimate}$ will serve as a building block for the estimate over the interval [0,1] on the rescaled manifold $M_\lambda$. In fact, by splitting the time interval $[0,T]$ into $\sim N_1T$ pieces each of length $\sim \frac{1}{N_1}$ and using the conservation of $L^2$ norm one easily gets the following corollary:

\begin{corollary}\label{time T bilinear estimate on M}
Suppose the $u_0, v_0 \in L^2(M)$ and $N_2 \leq N_1$ are dyadic scales. Then
\begin{equation}\label{bilinear estimate}
||e^{it\Delta} P_{N_1} u_0 e^{it\Delta}P_{N_2}v_0||_{L^2([0,T]\times M)} \leq \Lambda(T,N_1,N_2)||u_0||_{L^2(M)}||v_0||_{L^2(M)}
\end{equation}

where

\begin{equation}\label{Lambda T}
\Lambda(T,N_1,N_2) \lesssim\left\{
    \begin{array}{ll}
        \left(\frac{N_2}{N_1}\right)^{1/2}  & \mbox{if } T \leq N_1^{-1} \\
       \left(TN_2\right)^{1/2} & \mbox{if } T \geq N_1^{-1}
    \end{array}
\right.
\end{equation}

The same estimate holds if $v_0$ is supported at frequencies $\lesssim N_2$.

\end{corollary}

Next we translate this estimate via scaling into a bilinear Strichartz esimate on $[0,1]\times M_\lambda$. Recall that if $e_k(x)$ is an $L^2$ normalized eigen-basis of the Laplacian of $M$ with eigenvalues $\mu_k$, then $\frac{1}{\lambda}e_k(\frac{x}{\lambda})$ is an $L^2(M_\lambda)$-normalized eigen-basis of the Laplacian on $M_\lambda$ corresponding to eigenvalues $\frac{\mu_k}{\lambda^2}$. As a
corollary, we obtain the following:

\begin{corollary}[Time $T$ estimate on $M$ implies time $1$ estimate on $\lambda M$]\label{Strichartz on rescaled manifold}
Let $N_1,N_2 \in 2^{\Z}$ and suppose $u_0,v_0 \in L^2(M_\lambda)$. If $N_2 \leq N_1$ then:

\begin{eqnarray}\label{bilinear estimate}
||e^{it\Delta_{\lambda}} P_{N_1}u_0 e^{it\Delta_{\lambda}}P_{N_2}v_0||_{L^2([0,1]\times  M_\lambda)}
\lesssim& \Lambda(\lambda^{-2},\lambda N_1,\lambda N_2)||u_0||_{L^2(M_\lambda)}||v_0||_{L^2(M_\lambda)}\\
\lesssim& \left\{
    \begin{array}{ll}
        \left(\frac{N_2}{N_1}\right)^{1/2} ||u_0||_{L^2(M_\lambda)}||v_0||_{L^2(M_\lambda)} & \mbox{if } \lambda \geq N_1 \\
       \left(\frac{N_2}{\lambda}\right)^{1/2}||u_0||_{L^2(M_\lambda)}||v_0||_{L^2(M_\lambda)} & \mbox{if } \lambda \leq N_1
    \end{array}
\right.\label{the N_2/lambda decay}
\end{eqnarray}

The same estimate holds if $v_0$ is supported on frequencies $\lesssim N_2$.

\end{corollary}

Notice that in the ``flat space" limit $\lambda \to \infty$, the bilinear estimate becomes the same as that on $\R^2$.

\proof The frequency localizations of $u_0$ and $v_0$ mean that $\pi^\lambda_{\frac{\nu}{\lambda^2}}(u_0)=0$ unless
${\frac{\sqrt \nu}{\lambda}}\in [N_1,2N_1)$ and
$\pi^{\lambda}_{\frac{\nu}{\lambda^2}}(v_0)=0$ unless
${\frac{\sqrt \nu}{\lambda}}\in [N_2,2N_2)$.

\begin{align*}
||e^{it\Delta_\lambda} u_0
e^{it\Delta_\lambda}v_0||_{L^2([0,1]\times M_\lambda)}
=&\lambda||e^{i\lambda^2t\Delta_\lambda} u_0 e^{i\lambda^2t\Delta_\lambda}v_0||_{L^2([0,\lambda^{-2}]\times M_\lambda)}\\
=&\lambda ||\sum_{\sqrt \nu_k \sim \lambda N_1, \sqrt \nu_l \sim \lambda N_2}e^{-it(\nu_k+\nu_l)} \pi^\lambda_{\nu_k/\lambda^2}u_0(x) \pi^\lambda_{\nu_l/\lambda^2}v_0(x)||_{L^2([0,\lambda^{-2}]\times M_\lambda)}\\
\end{align*}
where we denoted by $\pi^{\lambda}_{\nu_l/\lambda^2}$ the orthogonal
projection operator in $L^2(M_\lambda)$ onto the eigenspace
corresponding to the eigenvalue $\nu_l/\lambda^2$. Define $\tilde u:
M\to \C$ as $\tilde u(x)=u(\lambda x)$. Using the fact that

$$
\pi_{\nu}\tilde u(x) = \pi^\lambda_{\nu/\lambda^2} u(\lambda x)
$$

we get:
\begin{align*}
||e^{it\Delta_\lambda} u_0e^{it\Delta_\lambda}v_0||_{L^2([0,1]\times M_\lambda)}\\
=&\lambda ||\sum_{\sqrt \nu_k \sim \lambda N_1, \sqrt \nu_l \sim \lambda N_2}e^{-it(\nu_k+\nu_l)} \pi_{\nu_k}\tilde u_0(\frac{x}{\lambda}) \pi_{\nu_l}\tilde v_0(\frac{x}{\lambda})||_{L^2([0,\lambda^{-2}]\times M_\lambda)}\\
=&\lambda^2 ||\sum_{\sqrt \nu_k \sim \lambda N_1, \sqrt \nu_l \sim \lambda N_2}e^{-it(\nu_k+\nu_l)} \pi_{\nu_k}\tilde u_0(x) \pi_{\nu_l} \tilde v_0(x)||_{L^2([0,\lambda^{-2}]\times M)}\\
\end{align*}

and hence

\begin{equation}\label{rescaling bilinear estimate}
||e^{it\Delta_\lambda} u_0e^{it\Delta_\lambda}v_0||_{L^2([0,1]\times M_\lambda)}= \lambda^2 ||e^{it\Delta} \tilde u_0e^{it\Delta}\tilde v_0||_{L^2([0,\lambda^{-2}]\times M)}.
\end{equation}

Applying \eqref{bilinear estimate} we get:
\begin{align*}
||e^{it\Delta_\lambda} u_0e^{it\Delta_\lambda}v_0||_{L^2([0,1]\times M_\lambda)} \lesssim& \lambda^2 \Lambda(\lambda^{-2}, \lambda N_1, \lambda N_2)||\tilde u_0||_{L^2(M)}||\tilde v_0||_{L^2(M)}\\
=&\Lambda(\lambda^{-2}, \lambda N_1, \lambda N_2)||u_0||_{L^2(M_\lambda)}|| v_0||_{L^2(M_\lambda)}.
\end{align*}

as claimed.

\endproof

\begin{remark} 
The  main gain provided by this lemma (at least for our purposes of applying the I-method) will be in the regime $N_2 \lesssim \lambda \ll N_1$ in which case the second inequality in $\eqref{the N_2/lambda decay}$ will be crucial in order to get the full $s>\frac{2}{3}$ well-posedness range. We should note that using the linear estimates alone is not sufficient even if one uses endpoint type estimates. In fact, in \cite{J}, it is proved that if $v_0$ is spectrally localized at $\sqrt{-\Delta}\in [N_2,2N_2)$ then:

$$
||e^{it\Delta}v_0||_{L_t^2L_x^\infty([0,N_2^{-1}]\times M)}\lesssim (\log N_2)^{1/2}||v_0||_{L^2(M)}.
$$

Combining this with the trivial $L_t^\infty L_x^2$ and H\"older one gets the following short-range estimate:

\begin{equation}\label{Jiang's bilinear estimate at scale N_2}
||e^{it\Delta} u_0 e^{it\Delta}v_0||_{L^2([0,N_2^{-1}]\times M)} \leq (\log N_2)^{1/2} ||u_0||_{L^2(M)}||v_0||_{L^2(M)}
\end{equation}

which would give at the time scale $T$:

\begin{equation}\label{Jiang's bilinear estimate at scale T}
||e^{it\Delta} u_0 e^{it\Delta}v_0||_{L^2([0,T]\times M)} \leq \left\{
    \begin{array}{ll}
        (\log N_2)^{1/2} ||u_0||_{L^2(M)}||v_0||_{L^2( M)} & \mbox{if }T \ll N_2^{-1}\\
       T^{1/2}N_2^{1/2+}||u_0||_{L^2(M)}||v_0||_{L^2( M)} & \mbox{if } T \gtrsim N_2
    \end{array}
\right.
\end{equation}

which, in turn, would translate to the following estimate on $M_\lambda$:

\begin{equation}\label{Jiang's estimate on lambdaM}
||e^{it\Delta_{\lambda}} u_0 e^{it\Delta_{\lambda}}v_0||_{L^2([0,1]\times M_\lambda)}
\lesssim \left\{
    \begin{array}{ll}
        [\log(\lambda N_2)]^{1/2} ||u_0||_{L^2(M_\lambda)}||v_0||_{L^2(M_\lambda)} & \mbox{if } \lambda \gg N_2 \\
       \frac{N_2^{1/2+}}{\lambda^{1/2-}}||u_0||_{L^2(M_\lambda)}||v_0||_{L^2(M_\lambda)} & \mbox{if } \lambda \lesssim N_2
    \end{array}
\right.
\end{equation}

When comparing $\eqref{the N_2/lambda decay}$ to $\eqref{Jiang's estimate on lambdaM}$, one first notices a substantial improvement in the range $\lambda \gtrsim N_1$ ($\left(\frac{N_2}{N_1}\right)^{1/2}$ versus $\left(\log{\lambda N_2}\right)^{1/2}$). This is similar to the improvement provided by the bilinear refinement to Strichartz estimate in \cite{B2} to linear Strichartz estimate on $\R^d$. Another crucial improvement (especially for our purposes of running the I-method)  happens in the range $N_2\ll\lambda\ll N_1$. In this range, $\eqref{the N_2/lambda decay}$ gives a bound of $\left(\frac{N_2}{\lambda}\right)^{1/2}\ll 1$ whereas $\eqref{Jiang's estimate on lambdaM}$ gives the large constant $(\log \lambda N_2)^{1/2}$ which is not enough to get the global well-posedness result.
\end{remark}

A couple of words about the proof of $\eqref{semiclassical bilinear estimate}$ seem to be in order. This will also allow us to justify a version of these bilinear estimates involving differential operators applied to the propagator $e^{it\Delta}$\footnote{These are stated in Corollary 3.6 of \cite{H}.}. The proof starts with the Burq, Gerard, Tzvetkov parametrix \cite{BGT} of $e^{it\Delta}$. Using this parametrix, one translates estimates like that in $\eqref{semiclassical bilinear estimate}$ into bilinear oscillatory integral operator estimate of the form $||T_\nu f S_\mu g||_{L^2_{t,x}(\R\times \R^d)}$ for operators of the form:

\begin{equation}\label{Fourier integral operator I}
T_\nu f(t,x)=\int_{\R^d}e^{i\nu \phi(t,x,\xi)}a(t,x,\xi)f(\xi)d\xi
\end{equation}

and 

\begin{equation}\label{Fourier integral operator II}
S_{\mu}g(t,x)=\int_{\R^d}e^{i\mu \psi(t,x,\xi)}b(t,x,\xi) g(\xi)d\xi
\end{equation}

where $\nu, \mu >0$, $a,b \in C_0^\infty(\R \times \R^d \times \R^d)$ and $\phi, \psi \in C^{\infty}$ are real-valued phase functions satisfying a non-degeneracy condition and another crucial transversality condition (see \cite{H} for details). This transversality condition is satisfied by the main terms of the parametrices considered at least when $N_2 \ll N_1$. One then applies a bilinear oscillatory integral estimate (Theorem 1 of \cite{H}) in order to obtain $\eqref{semiclassical bilinear estimate}$ in the range $N_1 \gg N_2$. The case when $N_2 \sim N_1$ follows from the linear Strichartz estimates.

The parametrices of $e^{i{t\over N_1}\Delta}u_0$ and $e^{i{t\over N_2}\Delta} v_0$ in \cite{BGT} allow us to write\footnote{Strictly speaking this representation only holds in an open neighborhood of $x_0\in M$. Since $M$ is compact, we can cover it by finitely many of such neighborhood, and hence we only need to prove the estimate on each one of them.}:

$$
e^{i\frac{t}{N_1}\Delta}u_0(x)=\tilde T_{N_1}u_0(t,x)+R_{N_1}u_0(t,x)
$$

and 
$$
e^{i{t\over N_2}\Delta}v_0(x)=\tilde S_{N_2} v_0(t,x)+R_{N_2} v_0(t,x)
$$

with

\begin{equation}\label{def of tilde T}
\tilde T_{N_1} u_0(t,x)=\frac{N_1^d}{(2\pi)^d}\int_{\R^d}e^{iN_1\tilde\phi(t,x,\xi)}a_1(t,x,\xi,N_1)\widehat{\tilde u_0}(N_1\xi)d\xi
\end{equation}

and 

\begin{equation}\label{def of tilde S}
\tilde S_{N_2} v_0(t,x)=\frac{N_2^d}{(2\pi)^d}\int_{\R^d}e^{iN_2 \tilde\phi(t,x,\xi_2)}a_2(t,x,\xi_2,N_2)\widehat {\tilde v_0}(N_2\xi_2)d\xi_2.
\end{equation}

Here $\tilde u_0$ and $\tilde v_0$ are the respective microlocalizations of $u_0$ and $v_0$ in the considered coordinate patch (in particular $||\tilde u_0||_{L^2(\R^d)}\lesssim ||u_0||_{L^2(M)}$ and similarly for $\tilde v_0$)(cf. \cite{BGT},\cite{H}) and $a_1, a_2 \in C_0^\infty(\R \times \R^d \times \R^d)$ are polynomials in $1\over N_1$ and $1 \over N_2$ respectively. The remainder operators $R_{N_1}$ and $R_{N_2}$ are smoothing operators that satisfy:

\begin{equation}\label{Remainder terms}
||R_{N_1} u_0||_{L^\infty_t H^\sigma([-\alpha,\alpha]\times M)}\lesssim_N N_1^{-N}||u_0||_{L^2(M)} \textrm{ and }
||R_{N_2} v_0||_{L^\infty_t H^\sigma([-\alpha,\alpha]\times M)}\lesssim_N N_2^{-N}||v_0||_{L^2(M)}
\end{equation}

for any $N$.

If $P(D)$ is a differential operator on $M$ of degree $n$, then $P(D)e^{i{t\over N_1}\Delta} u_0$ has the following expression:

$$
P(D)e^{i\frac{t}{N_1}\Delta}u_0(x)=N_1^n \tilde T'_{N_1}u_0(t,x)+R'_{N_1}u_0(t,x)
$$

where $\tilde T '_{N_1}$ and $R'_{N_1}$ are operators of the same form as $T_{N_1}$ and $R_{N_1}$. In particular, $T'_{N_1}$ has an expression as in $\eqref{def of tilde T}$ (just with different $a$) and $R'_{N_1}$ obeys the same estimates as in $\eqref{Remainder terms}$ (by choosing $N$ large enough). Similar expressions for $e^{i\frac{t}{N_2}}v_0$ allow us to deduce that following from the exact same analysis used to prove $\eqref{semiclassical bilinear estimate}$ (see Corollary 3.6 of \cite{H}):

\begin{corollary}\label{Strichartz on rescaled manifold with differential operators}
Suppose the $u_0, v_0 \in L^2(M)$ are spectrally localized around $N_1,N_2 \in 2^{\Z}$ respectively as in Corollary \ref{Strichartz on rescaled manifold}. Let $P(D)$ and $Q(D)$ be differential operators on $M$ of orders $n$ and $m$ respectively:

\begin{equation}\label{bilinear estimate with differential operators}
||P(D)e^{it\Delta} u_0 Q(D)e^{it\Delta}v_0||_{L^2([0,T]\times M)} \leq N_1^n N_2^m \Lambda(T,N_1,N_2)||u_0||_{L^2(M)}||v_0||_{L^2(M)}
\end{equation}

where $\Lambda(T,N_1,N_2)$ is given in $\eqref{Lambda T}$.
\end{corollary}

As before, we need to translate this ``time $T$" estimate on $M$ into a ``time 1" estimate on $M_\lambda$. In order to simplify the scaling, we will only define define differential operators on $M_\lambda$ as rescalings of differential operators on $M$. Suppose $\tilde P$ and $\tilde Q$ are differential operators on $M$: for any $f\in C^{\infty}(M_\lambda)$ we define the operator $P$ acting on $f$ as $P(D)f(x)=\left(D_\lambda^{-1}\right)^*\circ \left(\tilde P(D)\tilde f\right)=\left(\tilde P(D) \tilde f \right)(\frac{x}{\lambda})$ where $\tilde f \in C^{\infty}(M)$ is defined as before by $\tilde f(y)=f(\lambda y)$ for every $y\in M$ (In words, $P(D)f$ is obtained by pulling $f$ back to $M$ to get $\tilde f$, applying $\tilde P(D)$ to $\tilde f$, and finally pushing forward the resulting function to $M_\lambda$). With such conventions we have:

\begin{corollary}\label{time 1 bilinear estimate on lambda M with differential operators}
Suppose the $u_0, v_0 \in L^2(M_\lambda)$ are spectrally localized around $N_1,N_2 \in 2^{\Z}$ respectively as in Corollary \ref{time T bilinear estimate on M}. Let $\tilde P(D)$ and $\tilde Q(D)$ be differential operators on $M$ of orders $n$ and $m$ respectively and define $P(D)$ and $Q(D)$ on $C^\infty(M_\lambda)$ as indicated above. Then:

\begin{eqnarray}\label{bilinear estimate with differential operators}
\|P(D)e^{it\Delta} u_0 Q(D)e^{it\Delta}v_0\|_{L^2([0,1] \times M_\lambda)} \lesssim (\lambda N_1)^n (\lambda N_2)^m \Lambda(\lambda^{-2},\lambda N_1,\lambda N_2)\|u_0\|_{L^2(M_\lambda)}\|v_0\|_{L^2(M_\lambda)}\\
\lesssim  (\lambda N_1)^n (\lambda N_2)^m \left(\frac{N_2}{\lambda}\right)^{1/2}\|u_0\|_{L^2(M_\lambda)}\|v_0\|_{L^2(M_\lambda)} \;\; \textrm{if $\lambda \leq N_1$}
\end{eqnarray}

where $\Lambda(T,N_1,N_2)$ is given in $\eqref{Lambda T}$.
\end{corollary}

The proof is merely a rescaling of Corollary \ref{bilinear estimate with differential operators} performed as in the proof of Corollary \ref{Strichartz on rescaled manifold} which we will not repeat.

As for linear estimates, bilinear Strichartz estimates can be reformulated in terms of bilinear $X^{s,b}$ estimates. The reformulation of $\eqref{bilinear estimate}$ and $\eqref{bilinear estimate with differential operators}$ is the following:

\begin{corollary}\label{bilinear Xsb estimate lemma}
For any $b>1/2$ and any $f,g \in X^{0,b}([0,1]\times \lambda M)$ spectrally localized in dyadic regions around $N_1$ and $N_2$ respectively (i.e. $\mathbf{1}_{[N_1,2N_1]}(\sqrt{-\Delta}) f=f$ and $\mathbf{1}_{[N_2,2N_2]}(\sqrt{-\Delta}) g=g$), we have:

\begin{equation}\label{bilinear Xsb estimate}
\|f g\|_{L^2_{t,x}([0,1]\times M_\lambda)}\lesssim \Lambda(\lambda^{-2},\lambda N_1,\lambda N_2)\|f\|_{X^{0,b}([0,1]\times M_\lambda)}\|g\|_{X^{0,b}([0,1]\times M_\lambda)}.
\end{equation}

If $P(D)$ and $Q(D)$ are differential operators of orders $n$ and $m$ respectively defined as in Corollary \ref{time 1 bilinear estimate on lambda M with differential operators} then:

\begin{equation}\label{bilinear Xsb estimate with differential operators}
\|P(D)f\,Q(D)g\|_{L^2_{t,x}([0,1]\times M_\lambda)}\lesssim (\lambda N_1)^n(\lambda N_2)^m \Lambda(\lambda^{-2},\lambda N_1,\lambda N_2)\|f\|_{X^{0,b}([0,1]\times M_\lambda)}\|g\|_{X^{0,b}([0,1]\times M_\lambda)}.
\end{equation}

\end{corollary}

\proof We will only prove $\eqref{bilinear Xsb estimate with differential operators}$ as $\eqref{bilinear Xsb estimate}$ is merely a special case.  Without loss of generality, it is enough to assume that $f,g \in C_0^\infty([-2,2]\times M_\lambda)$. Let 

$$
F(t)=e^{-it\Delta}f(t) \textrm{ and } G(t)=e^{-it\Delta}g(t)
$$

Then 
\begin{align*}
P(D)f(t)=&P(D)e^{it\Delta}F(t)=\int_\R e^{it\tau_1}P(D)e^{it\Delta}\widehat{F}(\tau_1)d\tau_1 \textrm{ and }\\
Q(D)g(t)=&Q(D)e^{it\Delta}G(t)=\int_\R e^{it\tau_2}Q(D)e^{it\Delta}\widehat{G}(\tau_2)d\tau_2 .
\end{align*}

As a result,
\begin{align*}
|| P(D)f Q(D)g||_{L^2(\R \times M_\lambda)}=&\left|\left| \int_{\R_{\tau_1}}\int_{\R_{\tau_2}} e^{it(\tau_1+\tau_2)}P(D)e^{it\Delta}\widehat F(\tau_1)Q(D)e^{it\Delta}\widehat G(\tau_2)d\tau_1d\tau_2\right|\right|_{L^2([-2,2]\times M_\lambda)}\\
\leq& \int_{\R_{\tau_1}}\int_{\R_{\tau_2}}\left|\left| P(D) e^{it\Delta}\widehat F(\tau_1)Q(D) e^{it\Delta}\widehat G(\tau_2)\right|\right|_{L^2([-2,2]\times M_\lambda)}d\tau_1d\tau_2\\
\lesssim& (\lambda N_1)^n(\lambda N_2)^m \Lambda(\lambda^{-2},\lambda N_1,\lambda N_2)\int_{\R_{\tau_1}}\int_{\R_{\tau_2}}||\widehat F(\tau_1)||_{L^2(M_\lambda)}||\widehat G(\tau_2)||_{L^2(M_\lambda)}d\tau_1d\tau_2\\
\lesssim& (\lambda N_1)^n(\lambda N_2)^m   \Lambda(\lambda^{-2},\lambda N_1,\lambda N_2) ||\langle \tau_1\rangle\widehat F(\tau_1)||_{L_{\tau_1,x}^2(\R \times M_\lambda)}||\langle \tau_2\rangle\widehat G(\tau_2)||_{L_{\tau_2,x}^2(\R\times M_\lambda)}\\
=&(\lambda N_1)^n(\lambda N_2)^m   \Lambda(\lambda^{-2},\lambda N_1,\lambda N_2) ||f||_{X^{0,b}}||g||_{X^{0,b}}.
\end{align*}

\endproof

We conclude this section with a statement of a standard trilinear $X^{s,b}$ estimate that follows from $\eqref{the N_2/lambda decay}$. In fact, it is well known (see \cite{BEE} for example) that any bilinear Strichartz estimate of the form $||e^{it\Delta}u_0 \,e^{it\Delta}v_0||_{L^2_{t,x}}\lesssim \langle N_2\rangle^{s_0}||u_0||_{L^2}||v_0||_{L^2}$ (where $u_0$ and $v_0$ are spectrally localized dyadicly around frequencies $N_1\geq N_2$) would imply a trilinear $X^{s,b}$ estimate of the form $\eqref{trilinear estimate on M}$. This is made precise in the following lemma borrowed from \cite{BEE} (Proposition 2.5):

\begin{lemma}(Proposition 2.5 of \cite{BEE}) Suppose that a bilinear estimate of the form 
$$
||e^{it\Delta}u_0 \,e^{it\Delta}v_0||_{L^2_{t,x}([0,1]\times M)}\lesssim \langle N_2\rangle ^{s_0}||u_0||_{L^2(M)}||v_0||_{L^2(M)}
$$

holds whenever $u_0=\mathbf{1}_{[N_1,2N_1]}(\sqrt{-\Delta}) u_0$ and $v_0=\mathbf{1}_{[N_2,2N_2]}(\sqrt{-\Delta})v_0$ with $N_1\geq N_2$. Then for any $s>s_0$ there exists $(b,b')\in \R^2$ satisfying

\begin{equation}\label{condition on b and b'}
0<b'<\frac{1}{2} <b,\,\textrm{     }b+b'<1
\end{equation}

such that for any three functions $u_1,u_2,u_3 \in X^{s,b}(\R \times M)$:

\begin{equation}\label{trilinear estimate on M}
||u_1 \overline{u_2} u_3||_{X^{s,b'}}\lesssim ||u_1||_{X^{s,b}}||u_2||_{X^{s,b}}||u_3||_{X^{s,b}}
\end{equation}.

\end{lemma}

The fact that $\eqref{the N_2/lambda decay}$ holds on $M_\lambda$ (since $\Lambda(\lambda^{-2},\lambda N_1,\lambda N_2) \leq \langle N_2 \rangle^{1/2}$ for $\lambda >1$) implies, with the same proof as in \cite{BEE}, that there exists $(b,b')\in \R^2$ as in $\eqref{condition on b and b'}$ such that the following estimate holds:

\begin{equation}\label{trilinear estimate}
||u_1 \overline{u_2} u_3||_{X^{s,b'}([0,1]\times M_\lambda)}\lesssim ||u_1||_{X^{s,b}([0,1]\times M_\lambda)}||u_2||_{X^{s,b}([0,1]\times M_\lambda)}||u_3||_{X^{s,b}([0,1]\times M_\lambda)}
\end{equation}

for any $s>\frac{1}{2}$ where the implicit constant is independent of $\lambda$ for $\lambda>1$.

\section{Spectral localization: Part I}\label{part 1}

In this section, we deal with the problem of the spectral localization of products of eigenfunctions: given two eigenfunctions $f$ and $g$ of the Laplace-Beltrami operator on a compact $d-$manifold $M$, where is the product $fg$ spectrally localized? In other words, what can be said about $\pi_\mu (fg)$ where $\pi_\mu$ is the projection on the $\mu-$eigenspace of $(-\Delta)$. Of course, one can pose the same question for the product of any number of eigenfunctions. For our purposes, we will be most interested in the spectral concentration of the product of three eigenfunctions, but the same analysis carries on for any number of eigenfunctions. The results for two eigenfunctions mentioned in the introduction can be easily obtained by setting the third to be the constant eigenfunction 1. We should also remark as well that the results in this section apply for any smooth compact $d-$ dimensional Riemannian manifold $M^d$ including those with smooth\footnote{The exact regularity requirements will increase as a factor of the sharpness of the spectral localization one would like to prove. This is quantified by the number of iterations $n$ in Theorem \ref{A_0 and A_n theorem} which requires the eigenfunctions to be $C^{2n}(M)$ and hence it would be sufficient for the boundary (if it exists) to be $C^{2n}$.} boundary as long as one imposes either Dirichlet or Neumann boundary conditions on the eigenfunctions.

As mentioned in the introduction, this problem is trivial in the case of the torus and the sphere because of special eigen-function bases in those two cases (characters $e^{in.x}$ for $\T^d$ and spherical harmonics for  $S^d$). On a general compact manifold, the spectral localization of the product of two eigenfunctions  $e_2$ and $ e_3$ with eigenvalues $\mu_2$ and $\mu_3$ respectively on the $\mu_1$-eigenspace is detected via the inner products $\langle e_1, e_2e_3\rangle_{L^2(M)}$ where $e_1$ is an eigenfunction with eigenvalue $\mu_1$. As a result, the above problem is reduced to that of estimating integrals of the form $\int_M e_1(x)e_2(x) e_3(x) dx$ with $e_1,e_2,e_3$ being eigenfunctions of $-\Delta$ with eigenvalues $\mu_1, \mu_2, \mu_3$. Of course, if one is interested in the spectral localization of a product of more that two eigenfunctions $e_2$ and $e_3$, one needs to estimate integrals coming from $\langle e_1, \left(e_2\ldots e_k\right)\rangle$, namely $\int_M e_1(x)e_2(x)\ldots e_k(x)dx$. For the purpose of applying the I-method to the cubic nonlinear Schr\"odinger equation, we will be mainly interested in estimates for $k=4$.

The problem of identifying the spectral localization of the product of three eigenfunctions was encountered in \cite{BEE} where the following crude estimate was obtained:

\begin{lemma}\label{BGT localization}
Let $e_1,\ldots e_4$ be $L^2$ normalized eigenfunctions of the Laplacian with eigenvalues $\mu_1,
\ldots \mu_4$ respectively. There exists $C>0$ such that, if for $j=2,3,4$ we have $C\mu_j \leq \mu_1$, then for every $p>0$ there exists $C_p>0$ such that:

\begin{equation}\label{BGT localization equation}
\left|\int_M e_1(x) e_2(x)e_3(x) e_4(x) dx\right| \leq C_p \mu_1^{-p}.
\end{equation}
\end{lemma}

A proof of this lemma can be found in \cite{BEE} (lemma 2.6) and is based on the parametrix expression of the eigenfunctions and a simple non-stationary phase (integration by parts) argument. 

Unfortunately, the above estimate is way too crude for our purposes. To explain this point, we introduce the following notation. Since we will be interested in localizing $\sqrt{-\Delta_g}$ rather than $-\Delta_g$, it will be notationally convenient to denote, for $i=1,2,3,4$, the eigenvalues $\mu_i=n_i^2$ where $n_i\in [0,\infty)$. We also assume without any loss of generality that $n_2\geq n_3 \geq n_4$ and that $n_3>0$ (otherwise the answer is trivial). The suggestive cases of the torus $\T^d$ and the sphere $S^d$ suggest that the integral $\int_M e_1(x)e_2(x)e_3(x)e_4(x)dx$ should vanish (or at least present some  sort of decay) if $n_1 > n_2 +n_3 +n_4$ rather than $n_1>C(n_2+n_3+n_4)$ suggested by lemma \ref{BGT localization}. The importance of such an improvement, in comparison to lemma \ref{BGT localization}, is most crucial when $n_2 \gg n_3$, in which case we are multiplying a high frequency eigenfunction $e_2$ by two low frequency eigenfunctions $e_3$ and $e_4$. We expect that the resulting function, $e_2 e_3 e_4$, to be spectrally localized in the region $n_2 \pm O(n_3)$ rather than the much larger region $O(n_2)$ suggested by lemma \ref{BGT localization}. 

As a result, we are reduced to proving decay of the integral $\int_M e_1(e_2 e_3 e_4)dx$ when $n_1\geq n_2 +Kn_3$. We will be able to prove very fast (polynomial of any order) decay in terms of $K$ of this integral. This will follow from the following theorem whose proof is postponed to the section \ref{part 2}, as it requires the revision of some ideas from Riemannian geometry.

\begin{theorem}\label{A_0 and A_n theorem}
Let $M^d$ be a compact $d-$dimensional $C^\infty$ Riemannian manifold (possibly with boundary) and let $e_1,\ldots, e_4$ be eigenfunctions of the Laplace-Beltrami operator (with Dirichlet or Neumann boundary conditions) corresponding to eigenvalues $n_1^2 \ldots n_4^2$ respectively. Denote by $A_0$:

\begin{equation}\label{def of A_0}
A_0=\int_M e_1(x)\ldots e_4(x)dx.
\end{equation}

Then for any $n\in \N$:

\begin{equation}\label{A_0 and A_n}
A_0=\frac{(-2)^nA_n}{(n_1^2-n_2^2-n_3^2-n_4^2)^n}
\end{equation}

where $A_n$ is given by:

\begin{equation}\label{def of A_n}
A_n=\int_M e_1\left(B_n(e_2,e_3,e_4)+C_n(e_2,e_3,e_4)\right)dx
\end{equation}

and $B_n(f,g,h)$ and $C(f,g,h)$ are \emph{trilinear} operators given by:

\begin{equation}\label{def of B_n}
B_n(f,g,h)=\mathcal{O}_{\substack{i+j+k=2n \\ 0\leq i,j,k \leq n}}\left(\nabla^i f * \nabla^j g * \nabla^k h \right)
\end{equation}

\begin{equation}\label{def of C_n}
C_n(f,g,h)=\mathcal{O}_{\substack{i+j+k\leq 2(n-1) \\ 0\leq i,j,k \leq n-1}}\left(\tilde R_n*\nabla^i f * \nabla^j g * \nabla^k h \right).
\end{equation}

Here $\nabla^i f$ denotes the $i-th$ covariant derivative of $f$ (which is a (i,0) tensor) and for any two tensors $A$ and $B$, $A*B$ denotes some contraction of $A\otimes B$ and $\mathcal{O}(A*B)$ denotes a linear combination of contractions of $A\otimes B$. Here $\tilde R_n$ denotes a tensor obtained from the Riemann curvature tensor by contracting and differentiating it a bounded (in terms of $n$) number of times.

\end{theorem}

To see why this theorem provides the decay advertised above, we set $n_1=n_2+Kn_3$. This gives $n_1^2-(n_2^2+n_3^2+n_4^2) \geq 2Kn_2n_3 +(K^2-2)n_3^2\geq 2Kn_2n_3$ if $K \geq \sqrt 2$ and hence $\eqref{A_0 and A_n}$ gives that $|A_0|\leq \frac{1}{K^n}\left(\frac{|A_n|}{n_2^n n_3^n}\right)$. Given the structure of $A_n$ in $\eqref{def of A_n}$,$\eqref{def of B_n}$ and $\eqref{def of C_n}$ each derivative falling on $e_2$ is accompanied by a factor of $n_2$ in the denominator and each derivative falling on $e_3$ or $e_4$ is accompanied by a factor of $n_3$ in the denominator, which makes the term $\frac{|A_n|}{n_2^n n_3^n}$ essentially bounded at least from a heuristic point of view\footnote{The exact type of bound on $\frac{|A_n|}{n_2^n n_3^n}$ will depend on the context and the spaces involved. See for example Corollary \ref{spectral cluster lemma 1}. An estimate that does not involve any loss of derivatives can be obtained if one uses H\"older's inequality and estimates the eigenfunctions (or eigenfunction clusters) in $L^p$ spaces rather than $L^2$.}. The proof of Theorem \ref{A_0 and A_n theorem} is based on manipulations with covariant derivatives of the eigenfunction and is left to the last section \ref{part 2} in order not to distract the reader.

We should also mention that the fact that the operators $B_n$ and $C_n$ are multilinear allows one to use Theorem \ref{A_0 and A_n theorem}, which is a statement about single eigenfunctions $e_2, e_3, e_4$, to derive statements about eigenfunction clusters (see corollary below) or even Littlewood-Paley pieces of functions (see subsection \ref{decay of modified energy section}). A first instance of this is illustrated in the following proposition (cited in the introduction) which can be understood as a refinement of the bilinear Sogge estimates proved in \cite{BEE} by Burq, Gerard, and Tzvetkov:

\begin{corollary}\label{spectral cluster lemma 1}
Suppose that $\nu^2, \lambda^2, \mu^2$ are eigenvalues of the operator $-\Delta_g$ satisfying $\nu \geq \lambda \geq \mu$. If $\nu=\lambda+K\mu+2$ for some $K>1$ (i.e. $K=\frac{\nu-\lambda-2}{\mu}>1$), then for any $f,g \in L^2(M)$ and any $J\in \N$:

\begin{equation}\label{spectral cluster lemma 1 equation}
\pi_{\nu}\left(\mathbf{1}_{[\lambda,\lambda+1]}(\sqrt{-\Delta})f\mathbf{1}_{[\mu,\mu+1]}(\sqrt{-\Delta})g\right) \lesssim_J \frac{\Lambda(d,\mu)}{K^J}\|f\|_{L^2(M)}\|g\|_{L^2(M)}
\end{equation}

where
\begin{equation}\label{def of Lambda(d,mu)}
\Lambda(d,\mu):= 
\begin{cases}
\mu^{1/2} &\text{if } d=2\\
\mu^{1/2}(\log \mu)^{1/2} &\text{if } d=3\\
\mu^{\frac{d-2}{2}} &\text{if } d\geq 4
\end{cases}
\end{equation}

More generally, one has:
\begin{equation}
\int_M \mathbf{1}_{[\nu,\nu+1]}(\sqrt{-\Delta})h\,\mathbf{1}_{[\lambda,\lambda+1]}(\sqrt{-\Delta})f\,\mathbf{1}_{[\mu,\mu+1]}(\sqrt{-\Delta})g \,dx\lesssim_J \frac{\Lambda(d,\mu)}{K^J}\|h\|_{L^2(M)}\|f\|_{L^2(M)}\|g\|_{L^2(M)}.
\end{equation}

Interchanging the roles of $\lambda$ and $\nu$, the same estimates hold if $\nu <\lambda$ and $\nu=\lambda-K\mu -2$ with $K>1$ (i.e. if $K:=\frac{\lambda-\nu-2}{\mu}>1$).
\end{corollary}

The proof of this corollary is also left to the appendix as it requires some ideas from the next section. It is essentially a consequence of Theorem \ref{A_0 and A_n theorem} and the bilinear eigenfunction cluster estimates in \cite{BEE} and \cite{MEE}. It is worth mentioning that the spectral localization operator used here $\mathbf{1}_{[\lambda,\lambda+1]}(\sqrt{-\Delta})$ can be replaced by smoother versions like $\chi(\sqrt{-\Delta}-\lambda)$ with $\chi \in \mathcal{S}(\R)$. As before, this proposition is particularly useful when $\mu \ll \lambda$. It says that that the product $\mathbf{1}_{[\lambda,\lambda+1]}(\sqrt{-\Delta})f\,\mathbf{1}_{[\mu,\mu+1]}(\sqrt{-\Delta})g$ is spectrally localized (as measured by the $L^2$ norm of its projection onto various eigenspaces) in the region $\lambda +O(\mu)$ and starts to decay rapidly (faster than any polynomial power of $K$) as we move away from this region.

\section{A multilinear spectral multiplier lemma}\label{multilinear spectral multiplier lemma}

As is customary in previous applications of the I-method either on $\R^d$ or $\T^d$(cf. \cite{CKSTT1},\cite{dSPST}), one is faced with estimating $k-$linear multiplier forms: 

\begin{equation}\label{k linear convolution}
\begin{split}
\Lambda(f_1,\ldots,f_k)=&\int_{\xi_1+\ldots+\xi_k=0}\bar m(\xi_1,\xi_2,\ldots,\xi_k)\widehat{f_1}(\xi_1)\ldots\widehat{f_k}(\xi_k)d\xi_1\ldots\xi_k\\
=& \int \tilde m \left(\xi_1, \ldots, \xi_{k-1}, -\xi_1-\ldots-\xi_{k-1}\right) \widehat{f_1}(\xi_1)\ldots\widehat{f_k}(-\xi_1-\ldots -\xi_{k-1})d\xi_1\ldots\xi_{k-1}
\end{split}
\end{equation}

One is then interested in proving $L^2$ or $X^{s,b}$ type estimates for such forms, which, thanks to tools like Plancherel's theorem and Fourier inversion, can be done by proving weighted multilinear convolution estimates  (in $L^2$) for $\eqref{k linear convolution}$ either on $\R^d$ or $\Z^d$. This is due to the fact that $\eqref{k linear convolution}$ is actually a weighted convolution in Fourier space. Such convolution estimates make sense if we replace $\R^d$ or $\Z^d$ by any other additive abelian group $G$ and a systematic study of such estimates was done in \cite{Tmc}.

Once we move to the realm of general compact Riemannian manifolds and away from the category of abelian groups, multilinear spectral multipliers cease to be expressible as multilinear convolution operators. In fact, the operators with which we will be concerned have the form:

\begin{equation}\label{multilinear multiplier}
\Lambda(f_1,\ldots,f_k)=\sum_{n_1,\ldots, n_k}\bar m(n_1,\ldots,n_k)\int_{M}\pi_{n_1}f_1(x)\ldots \pi_{n_k}f_k(x)\, dx
\end{equation}

where $f_i(x)=\sum_{n_i}\pi_{n_i}f_i(x)$ is the spectral expansion of $f$ and $\pi_{n_i}f$ is the projection on the $n_i^2-$eigenspace. Note that in the case of the torus, one can use the Fourier expansion of $f$ and Fourier inversion to write $\eqref{multilinear multiplier}$ in the form $\eqref{k linear convolution}$ with $\xi_i\in \Z^d$. The estimates we will be interested in establishing take the form:

$$
|\Lambda(f_1,\ldots,f_k)|\lesssim_{\bar m} ||f_1||_Y\ldots ||f_k||_Y
$$

where $Y$ is some Banach space of functions like an $L^2$ Sobolev space or an $X^{s,b}$ space (in the latter case, one would be considering functions that depend on time and integrating in $t$ as well as in $x$ in the definition of $\Lambda$ in $\eqref{multilinear multiplier}$). 

In the case of convolution expressions like $\eqref{k linear convolution}$, one can use Fourier inversion to get such bounds. For example, in the case of the torus $\T^d$, one can bound $\eqref{multilinear multiplier}$ as follows:

$$
|\Lambda(f_1,\ldots, f_k)|\leq ||m||_{L^{\infty}}\int_{\xi_1+\ldots+\xi_k=0}|\widehat{f_1}(\xi_1)|\ldots |\widehat{f_k}(\xi_k)|d\xi_1\ldots d\xi_k=||m||_{L^\infty}\int_{\T^d} \tilde{f_1}(x) \ldots \tilde{f_k}(x) dx
$$

where $\tilde f_i(x)=\sum_{\xi}|f(\xi)|e^{2\pi i \xi.x}$. One is then reduced to estimating the above space integral. This can be done by H\"older's inequality for example (or any bilinear or quadrilinear estimate available). This would eventually give the desired estimate if one can show that $||\tilde f||_Y\lesssim ||f||_Y$, which is trivial if $Y$ is some $L^2$ based space like $H^s$ or $X^{s,b}$.

One could try to follow the same strategy above in estimating $\eqref{multilinear multiplier}$. Unfortunately, this leads directly to failure (except when $k=2$) because taking absolute values $|\pi_{n_i}f_i(x)|$ forbids us to recover $f_i$ again partly because $|\pi_{n_i}f_i(x)|$ is not an eigenfunction of the Laplacian anymore.

The key to estimating expressions like $\eqref{multilinear multiplier}$ for some bounded multipliers $\bar m$ is the following observation. If $\bar m$ is given as a tensor product of functions depending only on one variable at a time, i.e. $\bar m(\xi_1,\ldots,\xi_k)=m_1(\xi_1)\ldots m_k(\xi_k)$, one can postpone taking absolute values till after ``inverting" the spectral decomposition. More precisely,

$$
\Lambda(f_1,\ldots, f_k)=\sum_{n_i}\int_M \prod_{i=1}^k\left(m(n_i)\pi_{n_i}f_i(x)\right)dx=\int_M \prod_{i=1}^k \left(\sum_{n_i}m_i(n_i)\pi_{n_i}f_i(x)\right)dx=\int_M \tilde f_1(x)\ldots \tilde f_k(x)dx
$$

where $\tilde f_i(x)=\sum_{i}m_i(n_i)\pi_{n_i}f_i(x)$ is just a linear spectral multiplier applied to $f$. This would reduce the problem of estimating the multilinear multiplier into that of estimating the integral $\int \prod \tilde f_i dx$ in terms of $||\tilde f||_Y$. 

If $\bar m$ is not a multiplier in tensor form, one can still estimate $\eqref{multilinear multiplier}$ if $m$ satisfies some smoothness and decay properties. The idea is to first split the $n_i$ sums into dyadic pieces $n_i\sim N_i$ and then use Fourier series to write $\bar m(n_1,\ldots, n_k)=\sum_{\theta \in\Z^k}A(\theta)e^{i(\theta_1n_1+\ldots+ \theta_k n_k)}$. This reduces to the case considered above since each $e^{i(\theta_1n_1+\ldots+ \theta_k n_k)}$ is obviously in tensor form (cf. \cite{Chr} \cite{Hw}).

We won't formulate here the most general multilinear multiplier estimate in order not to distract the reader, but focus on a conditional lemma which is the case that will show up in the treatment of the I-method in the next section. For this, we assume that $Y$ is Banach space of functions that satisfies $||\tilde f||_Y\lesssim ||f||_Y$ whenever $\tilde f=\sum_{n_i} e^{i\theta_in_i}\pi_{n_i}f$ is a frequency modulation of $f$. In this case, we will say that $Y$ satisfies the ``\emph{modulation stability property}". In particular, this is the case for $X^{s,b}$ spaces and more generally $L^2$ Sobolev spaces. Such spaces are the analogues of translation-invariant Banach spaces on $\R^d$ and $\T^d$.

\begin{lemma}\label{multilinear multiplier lemma} Let $\Lambda$ is a k-linear multiplier as in $\eqref{multilinear multiplier}$ associated with the multiplier $\bar m$ and let $Y$ be a Banach space that satisfies the ``modulation stability property" as above. Assume that $\bar m$ satisfies the following symbol-type estimates\footnote{As is the case with usual multiplier theorems, one does not need to include derivatives of all orders for the theorem to be true. For this lemma, symbol estimates up to second derivatives are enough.}:
\begin{equation}\label{symbol estimates on bar m}
|\partial_{\xi_1}^{\alpha_1}\ldots\partial_{\xi_k}^{\alpha_k}\bar m(\xi_1,\ldots, \xi_k)|\lesssim \langle \xi_1\rangle^{-\alpha_1}\ldots \langle\xi_k\rangle^{-\alpha_k}.
\end{equation}

Suppose that the following estimate holds:
\begin{equation}\label{conditional estimate}
\left|\int_0^t\int_Mf_1(t),\ldots, f_k(t) dxdt\right| \leq  B ||f_1||_{Y}\ldots ||f_k||_{Y}
\end{equation}
whenever $f_i(t,.)=P_{N_i}f_i(t,.)$ are spectrally localized to frequency scales $N_i$. Then there exists a constant $C$ (depending only on implicit constants in $\eqref{conditional estimate}$) such that,

\begin{equation}\label{multilinear multiplier estimate}
|\int_0^t\Lambda(f_1(t),\ldots, f_k(t))dt|\leq CB ||f_1||_{Y}\ldots||f_k||_{Y}.
\end{equation}

\end{lemma}

\proof Since $n_i \sim N_i$, one can write $n_i=N_i \tilde n_i$ where $n_i \in [0,2)$. As a result, one can define a smooth function $\Psi \in C_c^\infty([-4,4])$ such that:

\begin{equation}\label{definition of multiplier Psi}
\Psi(\tilde n_1, \ldots, \tilde n_k)=\bar m(N_1 \tilde n_1,\ldots, N_k \tilde n_k)\;\;\;\textrm{ on $[0,2)^k$}.
\end{equation}

Thanks to $\eqref{symbol estimates on bar m}$ it is easy to see that $\Psi$ has bounded derivatives of all orders (only finitely many orders are needed actually). Extending $\Psi$ as a periodic function to $\R^k$ allows us to express it as a Fourier series:

$$
\Psi(\tilde n_1,\ldots, \tilde n_k)=\sum_{\theta_i \in \Z/4}A(\theta_1,\ldots, \theta_k)e^{i(\theta_1\tilde n_1+\ldots+\theta_k \tilde n_k)}.
$$

As a result of this, one can express $\Lambda$ as follows:

\begin{align*}
\Lambda(\phi_1,\ldots,\phi_k)=&\sum_{n_i\sim N_i}\sum_{\theta_i\in\Z/4}A(\theta_1,\ldots,\theta_k)e^{i(\theta_1 \frac{n_1}{N_1}+\ldots+\theta_k \frac{n_k}{N_k})}\int_M\pi_{n_i}f_i(t)\ldots\pi_{n_k}f_k(t)dx\\
=& \sum_{\theta_i \in \Z/4}A(\theta_1,\ldots,\theta_k)\int_M\left(\sum_{n_1\sim N_1}e^{i\theta_1 \frac{n_1}{N_1}}\pi_{n_1}f_1(t)\right)\ldots \left(\sum_{n_k\sim N_k}e^{i\theta_k \frac{n_k}{N_k}}\pi_{n_k}f_k(t)\right)dx\\
=& \sum_{\theta_i \in \Z/4}A(\theta_1,\ldots,\theta_k)\int_M\tilde f_1^{\theta_1}(t,x)\ldots \tilde f_k^{\theta_k}(t,x)dx\\
\end{align*}

where $\tilde f_i^{\theta_i}= \sum_{n_i\sim N_i}e^{i\theta_i \frac{n_i}{N_i}}\pi_{n_i}f_i(t)$ for $1\leq i \leq k$. As a result,

\begin{align*}
\left|\int_0^t \Lambda(f_1,\ldots,f_k)dt\right| \leq& \sum_{\theta_i\in \Z/4}|A(\theta_1,\ldots,\theta_k)|\left|\int_0^t\int_M \tilde f_1^{\theta_1}(t,x)\ldots \tilde f_k^{\theta_k}(t,x)dxdt\right|\\
\leq& B\sum_{\theta_i\in \Z/4}|A(\theta_1,\ldots,\theta_k)|\,||\tilde f_1^{\theta_1}||_Y\ldots ||\tilde f_k^{\theta_k}||_Y\\
\leq& C_1B\sum_{\theta_i\in \Z/4}|A(\theta_1,\ldots,\theta_k)|\,||f_1||_Y\ldots || f_k||_Y
\leq CB|| f_1||_Y\ldots || f_k||_Y
\end{align*}

where we have used $\eqref{conditional estimate}$ in the second inequality, the stability of the the $Y$ norm under spectral modulation in the third inequality, and finally the $L^1$ summability of $A$
which comes from the fact that $\Psi$ is $C^2([-4,4]^k)$ (see for example \cite{Gr}).\endproof

\section{The I-Method}\label{I-method section}

In this section we prove Theorem \ref{main theorem} by applying the I-method machinery. We first notice that if $U(t,x)$ solves $\eqref{cubic NLS}$ on $M$ over the interval $[0,T]$, then the function:

\begin{equation}\label{scaling on M}
u(t,x)=\frac{1}{\lambda}\tilde U(\frac{t}{\lambda^2},\frac{x}{\lambda})
\end{equation}

will solve the cubic nonlinear Schr\"odinger equation posed on the rescaled manifold $M_\lambda$, that is:

\begin{eqnarray}\label{NLS on M_lambda}
i\partial_t u +\Delta_\lambda u =&|u|^2u\\
u(0,x)=&u_0(x) \in H^s(M_\lambda)
\end{eqnarray}

over the interval $[0,\lambda^2 T]$\footnote{We may assume without loss of generality that the initial data $U_0$ is in $C^\infty(M)$. One can remove this assumption after proving the polynomial bound on $E$ using a standard limiting argument.}. We will be most interested in the range where $0<s<1$. Let $N\gg 1$ be a fixed large number to be specified later (depending only on $T$) and denote the eigenvalues of $-\Delta_{\lambda}$ by $n_1^2 < n_2^2 <\ldots$ where $n_i \in \R^+$ (Recall that $n_i^2 =\frac{\nu_i}{\lambda^2}$ where $\nu_i$ are the eigenvalues of $-\Delta_g$). We define the spectral multiplier:

\begin{equation}\label{def of Iu}
Iu=\sum_{n_i} m(n_i)\pi_{n_i}u
\end{equation}
where $m(n_i)$ is given by:

\begin{equation}\label{m(k)}
m(k)=\left\{
    \begin{array}{ll}
        1  & \mbox{if } k \leq N \\
        \left(\frac{N}{k}\right)^{1-s} & \mbox{if } k \geq 2N
    \end{array}
\right.
\end{equation}

and a smooth interpolent in between. For technical reasons, it will
be preferable to specify $m(k)=m_0(\frac{k}{N})$ where
$m_0: \R \to [0,1] \in C^{\infty}$ is non-increasing and satisfies:

\begin{equation}\label{def of m_0}
m_0(t)=\left\{
    \begin{array}{ll}
        1  & \mbox{if } t \leq 1 \\
        t^{-(1-s)} & \mbox{if } t \geq 2.
    \end{array}
\right.
\end{equation}

With this convention,
$Iu=\sum_{n_i}m_0(\frac{n_i}{N})\pi_{n_i}(u)$. Notice that $I$ is the identity operator on frequencies $\leq N$ and is an integration operator for frequencies $ \geq 2N$. As a result, $I$ is smoothing of order $1-s$, in fact for any $s_0 \in \R$:
\begin{equation}\label{u vs Iu in H^s}
||u||_{H^{s_0}(M)} \leq ||Iu||_{H^{s_0+(1-s)}} \leq N^{1-s}||u||_{H^{s_0}}.
\end{equation}

In particular, for $u \in H^s(M_\lambda)$, $Iu \in H^1 (M_\lambda)$ and we are justified to define the modified energy:

\begin{equation}\label{modified energy}
\tilde{E}[u]=E[Iu]=\int_{M_\lambda} \frac{1}{2}|\nabla_g{Iu}(t,x)|^2+\frac{1}{4}|Iu(t,x)|^4 dx.
\end{equation}

Notice that $\tilde{E}[u(t)]$ controls $||u(t)||_{\dot H^s}$ and
hence the boundedness of $\tilde{E}[u]$ over the interval $[0,\lambda ^2T]$ implies that $||u||_{\dot H^s}$ remains bounded (in particular it does not blowup). Our goal is then to show that for any $T>0$ the modified energy $\eqref{modified energy}$ remains bounded. 

Despite the fact that $\tilde{E}[u]$ is not a conserved quantity in general (since $Iu$ does not solve $\eqref{NLS on M_lambda}$ in general), it is almost conserved: in the sense that its rate of change will be a negative power of $N$. This will allow us to prove polynomial (in $T$) bounds for $\tilde{E}[u]$.

The proof is in steps: First we prove that the equation satisfied by $Iu$ is locally well-posed in $X^{1,b}([0,\delta]\times M_\lambda)$ with $b=\frac{1}{2}+$ and $\delta \gtrsim 1$. This will allow us to obtain the bound $||Iu||_{X^{1,\frac{1}{2}+}}\lesssim_{||Iu_0||_{H^1}} 1$ (compare to the bound from $\eqref{u vs Iu in H^s}$ obtained from the local well-posedness of $\eqref{NLS on M_lambda}$) that will be used in bounding the increment of the modified energy over the interval $[0,\delta]$. The latter will be proved in Proposition \ref{almost conservation} using a multilinear analysis of $\partial_t \tilde E[u](t)$. Finally, the process is repeated over sub-intervals of length $\delta$ partitioning $[0,\lambda^2 T]$ as long as the $\tilde E[u](t)\lesssim 1$. This latter condition will specify the number of iterations allowed and hence the dependence of $N$ on $T$ as well as the polynomial bounds on $\tilde E[u(t)]$.

\remark In contrast to previous applications of the I-method (cf \cite{CKSTT1}, \cite{B4},\cite{dSPST}), the increment on the modified energy will be bounded not only by a negative power of $N$ but also by a negative power of the scaling parameter $\lambda$ (which will eventually be chosen to be a positive power of $N$). This is the main reason why one is able to prove that global well-posedness holds on the full range $s>2/3$ similar to that on the torus where better Strichartz estimates hold (but with worse dependence on $\lambda$). The reason for this gain is the bilinear Strichartz estimates $\eqref{the N_2/lambda decay}$ and the fact that $\lambda$ itself will be chosen in the end to be a positive power of $N$ (namely $\lambda \sim N^{\frac{1-s}{s}}$).

\subsection{Local well-posedness of the I-system:}\label{well-posedness of the I-system} If $u$ satisfies \eqref{NLS on M_lambda}, then $Iu$ satisfies:

\begin{eqnarray}\label{Iu system}
i\partial_t Iu +\Delta_g Iu =&I(|u|^2u)\\
Iu(0,x)=&Iu_0(x) \in H^1(M_\lambda).
\end{eqnarray}

We shall need a local well-posedness result for the above I-system:

\begin{proposition}\label{LWP of I system}
 Suppose $u_0\in H^s(M_\lambda)$ is such that $||Iu_0||_{H^{1}}\lesssim 1$, then there exists $0<\delta \sim 1$ such that \eqref{Iu system} is locally well-posed on $[0,\delta]$ and $||Iu||_{X^{1,\frac{1}{2}+}([0,\delta]\times M_\lambda)}\lesssim 1$.
\end{proposition}

\proof $Iu$ satisfies the following integral equation on the interval $[0,\delta]$:

$$
Iu(t)=e^{it\Delta}Iu_0 - \int_0^t e^{i(t-s)\Delta}I\left(|u(s)|^2u(s)\right) ds$$

and hence we have for some $0<b'<\frac{1}{2}<b$ with $b+b'<1$:

\begin{equation}\label{energy estimate}
||Iu\|_{X^{1,b}([0,\delta]\times M_\lambda)} \lesssim ||Iu_0||_{H^1(M_\lambda)}+ \delta^{1-(b+b')}||I\left(|u|^2u\right) ||_{X^{1,-b'}([0,\delta]\times M_\lambda)}.
\end{equation}

This follows from the retarded estimate on the Duhamel term in $X^{s,b}$ spaces (See for instance Lemma 3.2 of \cite{G}, section 2.6 of \cite{T}, or Proposition 2.11 of \cite{BEE} in the context of compact manifolds). As a result, we will be done once we show that $\|I\left(|u|^2u\right)||_{X^{1,-b'}} \lesssim ||Iu||_{X^{1,b}}^3$. This will in turn follow from:

\begin{equation}\label{I(cubic) estimate}
||I\left(u_1 \overline{u_2} u_3\right)||_{X^{1,-b'}} \lesssim \prod_{i=1}^3||Iu_i||_{X^{1,b}}.
\end{equation}

We will deduce \eqref{I(cubic) estimate} from the cubic estimate $\eqref{trilinear estimate}$ which we recall for convenience:

\begin{equation}\label{cubic estimate}
||u_1 \overline{u_2} u_3||_{X^{s',-b'}} \lesssim \prod_{i=1}^3||u_i||_{X^{s',b}}
\end{equation}

for all $s'>1/2$.

To prove \eqref{I(cubic) estimate} from $\eqref{cubic estimate}$, we split into two cases:

\emph{Case 1:} First, suppose that $P_{\leq 3N} u_i=u_i$ for all $i=1,2,3$, where $P_{\leq 3N}u_i =P_{\sqrt{-\Delta} \leq 3N}u_i$. As a result:

$$
||I(u_1\overline{u_2}u_3)||_{X^{1,-b'}} \leq ||u_1\overline{u_2}u_3||_{X^{1,-b'}}\lesssim \prod_{i=1}^{3}||u_i||_{X^{1,b}}\lesssim \prod_{i=1}^{3}||Iu_i||_{X^{1,b}}
$$
where in the first inequality we used that $m$ is bounded by 1, in the second we used the trilinear estimate $\eqref{cubic estimate}$, and in the third inequality we used that $m(k)\sim 1$ when $\sqrt{-\Delta}\lesssim N$.

\emph{Case 2:} Now suppose that one of the $u_i$, say $u_1$ for definiteness, satisfies $P_{\leq 2N}u_1=0$, then using equation \eqref{u vs Iu in H^s} we have:
$$
||I(u_1\overline{u_2}u_3)||_{X^{1,-b'}} \lesssim N^{1-s}||u_1\overline{u_2}u_3||_{X^{s,-b'}}\lesssim N^{1-s}||u_1||_{X^{s,b}}||u_2||_{X^{s,b}}||u_3||_{X^{s,b}}\lesssim \prod_{i=1}^{3}||Iu_i||_{X^{1,b}}
$$
where we used $\eqref{u vs Iu in H^s}$ for the first inequality, the trilinear estimate \eqref{cubic estimate} for the second inequality, and for the third we used \eqref{u vs Iu in H^s} and the fact that $||u_1||_{X^{s,b}} \leq N^{-(1-s)}||Iu_1||_{X^{1,b}}$.

Now \eqref{I(cubic) estimate} follows by decomposing each $u_i=P_{\leq 3N}u_i+(u_i-P_{\leq 3N}u_i)$.
\endproof

\begin{remark} The proof of the above proposition explains why we choose to work on the manifold $M_\lambda$ rather than $M$. Rescaling $M$ to $M_\lambda$ allows for the normalization $\|Iu_0\|_{H^1}\sim 1$ (cf. section \ref{polybounds}) which yields for a time of existence $\delta \sim 1$. Without this rescaling, we would need to know the sharp dependence of $\delta$ on $\|Iu\|_{H^1}$. Unfortunately, this sharp dependence is not provided directly\footnote{though one can recover it by the scaling argument we do.} by \eqref{energy estimate} and $\eqref{I(cubic) estimate}$ which lead to sub-optimal results.

\end{remark}

\subsection{Decay of the modified energy}\label{decay of modified energy section}

We are now ready to prove the almost conservation of the modified energy:

\begin{proposition}\label{almost conservation}
Let $s>\frac{2}{3}$, $u_0\in H^s(M_\lambda)$ with $0<\lambda<N$. If $u(t)$ solves \eqref{NLS on M_lambda} and $||Iu||_{X^{1,\frac{1}{2}+}([0,\delta]\times M_\lambda)}\lesssim 1$, then

\begin{equation}\label{decay of modified energy}
\left|\tilde{E}[u(t)]-\tilde{E}[u(0)]\right| \lesssim \frac{1}{\lambda N^{\frac{1}{2}-}}
\end{equation}
for all $0\leq t \leq \delta$
\end{proposition}

\proof Computing the time derivative of $\tilde E[u(t)]$ we get:

\begin{align*}
\partial_t \tilde E[u(t)]=&\partial_t E(Iu(t))=\operatorname{Re} \int_{M_\lambda}
\overline{Iu_t}\left(-\Delta Iu+|Iu|^2Iu\right)dx\\
=&\operatorname{Re} \int_{M_\lambda}
\overline{Iu_t}\left(-iIu_t-\Delta Iu+|Iu|^2Iu\right) dx
=\operatorname{Re} \int_{M_\lambda}
\overline{Iu_t}\left(-I(|u|^2u)+|Iu|^2Iu\right)dx
\end{align*}

where we have used \eqref{Iu system}. Writing $u(t,x)=\sum_{n}\pi_{n}u(t,x)$ and using the fundamental theorem of calculus we get:

\begin{align*}
&\tilde E[u(t)]-\tilde E[u(0)]\\
&=\operatorname{Re}\sum_{n_i}\int_0^t\int_{M_\lambda} m(n_1)\overline{\pi_{n_1}u_t}
\left[m(n_2)m(n_3)m(n_4)-I\right]\pi_{n_2}u\,\overline{\pi_{n_3}u}\,\pi_{n_4}u \,dx\,dt.
\end{align*}

Since $I$ is self-adjoint we can write this as:

\begin{align*}
&\tilde E[u(t)]-\tilde E[u(0)]\\
&=\operatorname{Re}\sum_{n_i}\int_0^t \int_{M_\lambda}
\left[1-\frac{m(n_1)}{m(n_2)m(n_3)m(n_4)}\right]\overline{\pi_{n_1}(Iu_t)}\,\pi_{n_2}(Iu)\,\overline{\pi_{n_3}(Iu)}\,\pi_{n_4}(Iu)\,dx\,dt.
\end{align*}

Using \eqref{Iu system} once more we get that:

\begin{equation}\label{sum of terms}
\tilde E[u(t)]-\tilde E[u(0)]=\operatorname{Re}i\left(\operatorname{Term}_1 -\operatorname{Term}_2\right)
\end{equation}

where

\begin{equation}\label{definition of term1}
\operatorname{Term}_1=\sum_{n_i}\int_0^t\int_{M_\lambda}
\left[1-\frac{m(n_1)}{m(n_2)m(n_3)m(n_4)}\right]\overline{\pi_{n_1}(\Delta Iu)}\pi_{n_2}(Iu)\overline{\pi_{n_3}(Iu)}\pi_{n_4}(Iu) dx\,dt
\end{equation}

\begin{equation}\label{definition of term2}
\operatorname{Term}_2=\sum_{n_i}\int_0^t\int_{M_\lambda}
\left[1-\frac{m(n_1)}{m(n_2)m(n_3)m(n_4)}\right]\overline{\pi_{n_1}\left(I\left(|u|^2u\right)\right)}\pi_{n_2}(Iu)\overline{\pi_{n_3}(Iu)}\pi_{n_4}(Iu)dx\,dt
\end{equation}

\subsubsection{Bound on $\operatorname{Term}_1$}We start by estimating $\operatorname{Term}_1$. Our goal is to prove that:

\begin{equation}\label{decay of term1}
\operatorname{Term}_1 \lesssim_{||Iu||_{X^{1,1/2+}}} \frac{1}{\lambda N^{1/2-}}.
\end{equation}

For this we break $u$ into a dyadic sum $u=\sum_{N} u_N$ where $N=2^j$, for $j=0,1,2, \ldots$ with $u_N=P_{[N,2N)}u$ for $N>1$ and $u_1=P_{[0,2)}u$. This reduces estimating $\operatorname{Term}_1$ into that of an integral involving dyadic frequency pieces at scale $N_j$. We will be able to sum over all those dyadic pieces by making sure that our estimates include a geometric decay factor in the highest frequency involved, that would allow using Cauchy-Schwarz to recover $||u||_{X^{s,b}}$. We present the details:

First notice that:

$$
||\Delta Iu||_{X^{-1,1/2+}} \leq ||Iu||_{X^{1,1/2+}}.
$$

Therefore, in order to estimate $\operatorname{Term}_1$ and prove \eqref{decay of term1}, we only need to show that:

\begin{equation}\label{decay of term1 in dyadic version}
\begin{split}
\left|\sum_{n_i\sim N_i}\int_0^t\int_{M_\lambda}
\left[1-\frac{m(n_1)}{m(n_2)m(n_3)m(n_4)}\right]\overline{\pi_{n_1}(\phi_1)}\pi_{n_2}(\phi_2)\overline{\pi_{n_3}(\phi_3)}\pi_{n_4}(\phi_4) dxdt\right| \\
\lesssim \frac{1}{\lambda N^{1/2-}}(N_1N_2N_3N_4)^{0-}||\phi_1||_{X^{-1,1/2+}}\prod_{i=2}^4||\phi_i||_{X^{1,1/2+}}
\end{split}
\end{equation}

whenever the $\phi_j$ are spectrally supported on $n_i \sim N_i$\footnote{In one particular case (Case 4 below), we won't prove the exact form in $\eqref{decay of term1 in dyadic version}$ but rather use Cauchy-Schwarz in $N_2$ (since we will have that $N_1\sim N_2$) and geometric series summation to sum in $N_3$ and $N_4$.}
Noticing that the above is symmetric\footnote{Strictly speaking the symmetry is broken due to the existence of complex conjugates. However, these will not affect the analysis in any way and hence the treatment of the other cases is similar.} in $n_2,n_3,n_4$, we will assume without loss of generality that

\begin{equation}\label{N2 and N3 N4}
N_2\geq N_3 \geq N_4.
\end{equation}

We now conduct a case by case analysis depending on how $N_2$ compares to $N_1$ and $N$. In the particular case when $M=\T^2$ or $S^2$, one can directly assume that $N_1 \lesssim N_2$ because the right hand side of \eqref{decay of term1 in dyadic version} vanishes if $n_1 \geq n_2+n_3+n_4$ thanks to the sharp spectral localization mentioned in section \ref{part 1} on those domains. Fortunately, this is a minor issue and will be dealt with using the following lemma:

\begin{lemma}\label{N1 >> N2 contribution}
There exists $C>0$ such that, if $N_1 \geq CN_2$, then for every $q>0$,
\begin{equation}\label{dealing with the N1>N2 case}
\text{L.H.S. of } \eqref{decay of term1 in dyadic version} \lesssim_q \frac{1}{N_1^q}\prod_{i=1}^4||\phi_i||_{X^{0,1/2+}}.
\end{equation}
\end{lemma}

\proof The proof of this lemma is straightforward using lemma \ref{BGT localization}. We include it here for completeness. We will use the convention that for any function $f: M_\lambda \to \C$, we shall denote by $\tilde f:M\to \C$ the pull-back function defined for $y\in M$ as $\tilde f(y)=f(\lambda y)$. Rescaling the integral in the L.H.S. of $\eqref{decay of term1 in dyadic version}$ back into an integral over $M$ and using the fact that $\pi_{n_i}f(x)=\pi_{\lambda n_i} \tilde f(\frac{x}{\lambda})$ we get:

\begin{align*}
\text{L.H.S. of } \eqref{decay of term1 in dyadic version}=&\lambda^2 \sum_{n_i\sim N_i}\left|\left(1-\frac{m(n_1)}{m(n_2)m(n_3)m(n_4)}\right) \int_0^t\int_M \overline{\pi_{\lambda n_1}(\tilde \phi_1)}\pi_{\lambda n_2}(\tilde \phi_2)\overline{\pi_{\lambda n_3}(\tilde \phi_3)}\pi_{\lambda n_4}(\tilde \phi_4) dx\,dt\right|\\
\lesssim& \lambda^2 \sum_{n_i\sim N_i}\frac{1}{(\lambda n_1)^p}\left|1-\frac{m(n_1)}{m(n_2)m(n_3)m(n_4)}\right|\int_0^t\prod_{i=1}^4||\pi_{\lambda n_i}\tilde \phi_i(t)||_{L^2(M)}dt\\
\end{align*}

by applying lemma \ref{BGT localization}. Using the fact that $\frac{m(n_1)}{m(n_2)}\leq 1$ since $m$ is non-increasing and the fact that $\frac{m(n_j)}{\langle n_j\rangle }\leq 1$ for $j=3,4$, we can crudely bound the whole multiplier $\left|1-\frac{m(n_1)}{m(n_2)m(n_3)m(n_4)} \right|\lesssim n_1^2$. Also, by applying the Weyl asymptotics:

$$
\# \{ \nu: N_j \leq \sqrt{\nu} \leq 2N_j \text{ with $\nu$ an eigenvalue of $-\Delta_g$}\} \lesssim N_j^2
$$

and hence by Cauchy-Schwarz we get that $\sum_{n_i\sim N_i}||\pi_{\lambda n_i}\tilde \phi_i(t)||_{L^2(M)}\lesssim (\lambda N_i)||\tilde \phi_i||_{L^2(M)}$. As a result, we get:

\begin{align*}
&\text{L.H.S. of }\eqref{decay of term1 in dyadic version} 
\lesssim \frac{\lambda^2}{\lambda^p N_1^{p-2}}\int_0^t \prod_{i=1}^4 (\lambda N_i)||\tilde \phi_i||_{L^2(M)}dt
\lesssim \frac{1}{\lambda^{p-2}N_1^{p-6}}\int_0^t\prod_{i=1}^4||\phi_i(t)||_{L^2( M_\lambda)}dt\\
&\leq  \frac{1}{\lambda^{p-2}N_1^{p-6}}\|\phi_1(t)\|_{L_{t,x}^2([0,t]\times M_\lambda)}\|\phi_2(t)\|_{L_{t,x}^2([0,t]\times M_\lambda)}\|\phi_3(t)\|_{L_{t}^\infty L_x^2([0,t]\times M_\lambda)}\|\phi_4(t)\|_{L_{t}^\infty L_x^2([0,t]\times M_\lambda)}\\
&\lesssim \frac{1}{N^{p-6}}\prod_{i=1}^4||\phi_i||_{X^{0,1/2+}}.
\end{align*}

\endproof

Since the multiplier on the LHS of $\eqref{decay of term1 in dyadic version}$ vanishes when $N_2 \ll N_1 \leq N$, the above lemma is more than enough to give the needed decay on the RHS of $\eqref{decay of term1 in dyadic version}$.

 \remark The same argument (with the roles of $N_1$ and $N_2$ interchanged) gives that:

$$
\text{L.H.S. of } \eqref{decay of term1 in dyadic version} \lesssim_q \frac{1}{N_2^q}\prod_{i=1}^4||\phi_i||_{X^{0,1/2+}}
$$

unless $N_2 \sim \max\{N_1,N_3\}$.

As a result, we will assume from now on that

\begin{equation}\label{N2 and N1}
N_1 \lesssim N_2.
\end{equation}
The analysis will be divided into several cases by comparing $N_2$ to $N$ and the other frequencies:

\textbf{Case 1. $N_2 \ll N$:} The bound is trivially true since by \eqref{N2 and N1} and \eqref{N2 and N3 N4} the multiplier on the left hand side of \eqref{decay of term1 in dyadic version} is zero.

\textbf{Case 2. $N_2 \gtrsim N \gg N_3\geq N_4$:} This is the most delicate case. Applying lemma \ref{N1 >> N2 contribution} (and interchanging the roles of $N_1$ and $N_2$) we may assume without any loss of generality that $N_2 \sim N_1$.

Split the dyadic interval $[N_1,2N_1)$ into $J$ intervals $I_\alpha$ of length $N_3$ each, and the interval $[N_2,2N_2)$ into $K$ intervals $I_\beta$ of length $N_3$ each. Then $J\sim K \sim \frac{N_2}{N_3}$. As a result of this, we have:

\begin{equation}\label{case 2 split into subintervals}
\begin{split}
&\textrm{L.H.S. of \eqref{decay of term1 in dyadic version} }\\
&=\left|\sum_{I_\alpha,I_\beta} \sum_{\substack{n_1 \in I_\alpha, n_2 \in I_\beta\\ n_3 \sim N_3,n_4 \sim N_4}} \int_0^t\int_{M_\lambda}
\left[1-\frac{m(n_1)}{m(n_2)m(n_3)m(n_4)}\right]\overline{\pi_{n_1}(\phi_1)}\pi_{n_2}(\phi_2)\overline{\pi_{n_3}(\phi_3)}\pi_{n_4}(\phi_4) dxdt\right|.
\end{split}
\end{equation}

We analyze \eqref{case 2 split into subintervals} differently according to the respective locations of the two intervals $I_\alpha$ and $I_\beta$. To be more precise, we will call $S_1$ the sum in \eqref{case 2 split into subintervals} when the interval $I_\alpha$ is at distance at least $8N_3$ \emph{to the right of} $I_\beta$ (i.e. elements in $I_\alpha$ are at least $8N_3$ larger than those in $I_\beta$), $S_2$ that when $I_\alpha$ is at distance at least $8N_3$ \emph{to the left of} $I_\beta$, and $S_3$ when the two intervals are distances $\leq 8N_3$ apart.

We start by estimating $S_1$: First, we fix some notation. In this case $N_1\geq N_2$, hence we can write $N_1=N_2+RN_3$ with $R \in \N \cup 0$ (since $N_1,N_2$ are dyadic integers bigger than $N_3$). As a result, $[N_1,2N_1)=\cup_{\alpha=R}^{R+J-1}I_\alpha$ with $I_\alpha=[N_2+\alpha N_3, N_2+(\alpha+1)N_3)$. Similarly, $[N_2,2N_2)=\cup_{\beta=0}^{K-1} I_\beta=[N_2+\beta N_3,N_2+(\beta+1)N_3)$. Since we are summing with $(\alpha,\beta) \in S_1$, we have $\alpha -\beta > 8$.

\begin{align*}
&\sum_{S_1} \int_0^t\int_{M_\lambda}
\left[1-\frac{m(n_1)}{m(n_2)m(n_3)m(n_4)}\right]\overline{\pi_{n_1}(\phi_1)}\pi_{n_2}(\phi_2)\overline{\pi_{n_3}(\phi_3)}\pi_{n_4}(\phi_4) dxdt\\
&=\sum_{S_1}\lambda^2\int_0^t\int_M \left[1-\frac{m(n_1)}{m(n_2)m(n_3)m(n_4)}\right]\overline{\widetilde{\pi_{n_1}\phi_1}}\widetilde{\pi_{n_2}\phi_2}\overline{\widetilde{\pi_{n_3}\phi_3}}\widetilde{\pi_{n_4}\phi_4} dydt\\
\end{align*}

where, as before, for any function $f:M_\lambda \to \C$, $\tilde f :M \to \C$ is given by $\tilde f (y)=f(\lambda y)$ for every $y\in M$. Using the fact that $\widetilde{ \pi_{n} f}(x)=\pi_{\lambda n}\tilde f(x)$ from $\eqref{rescaling projections}$, we get by applying Theorem \ref{A_0 and A_n theorem} that the latter expression is equal to:

\begin{align*}
&\sum_{S_1} \int_0^t\int_{M_\lambda}
\left[1-\frac{m(n_1)}{m(n_2)m(n_3)m(n_4)}\right]\overline{\pi_{n_1}(\phi_1)}\pi_{n_2}(\phi_2)\overline{\pi_{n_3}(\phi_3)}\pi_{n_4}(\phi_4) dxdt\\
=&\sum_{S_1}\left[1-\frac{m(n_1)}{m(n_2)m(n_3)m(n_4)}\right]\frac{\lambda^2 (-2)^l}{\left((\lambda n_1)^2-(\lambda n_2)^2-(\lambda n_3)^2-(\lambda n_4)^2\right)^l}\int_0^t\int_M \overline{\widetilde{\pi_{n_1}\phi_1}}\\
&\times \left(\tilde B(\widetilde{\pi_{n_2}\phi_2},\overline{\widetilde{\pi_{n_3}\phi_3}},\widetilde{\pi_{n_4}\phi_4})+\tilde C(\widetilde{\pi_{n_2}\phi_2},\overline{\widetilde{\pi_{n_3}\phi_3}},\widetilde{\pi_{n_4}\phi_4})\right) dy dt\\
\end{align*}
 
where
\begin{equation}\label{definition of tilde Bn1n2n3l}
\tilde B_l(f,g,h)=\mathcal{O}_{\substack{i+j+k=2l \\ 0\leq i,j,k \leq l}}\left(\nabla^i f * \nabla^j g * \nabla^k h \right)
\end{equation}

\begin{equation}\label{definition of tilde Cn1n2n3l}
\tilde C_l(f,g,h)=\mathcal{O}_{\substack{i+j+k\leq 2(l-1) \\ 0\leq i,j,k \leq l-1}}\left(\nabla^aR*\nabla^i f * \nabla^j g * \nabla^k h \right).
\end{equation}

In particular, both $\tilde B$ and $\tilde C$ are trilinear operators that are linear combinations of products of differential operators on $M$ of the form $\tilde Q _1(f)\,\tilde Q_2 (g)\,\tilde Q_3 (h)$ whose respective orders $i,j,k$ satisfy $i+j+k \leq 2l$ and $0\leq i,j,k \leq l$ (where $l$ will be chosen large enough so that the decay factor in $|\alpha-\beta|$ coming from the denominator $n_1^2-n_2^2-n_3^2-n_4^2$ would cancel a growth factor in $|\alpha-\beta|$ that comes from the multiplier thus giving a summable contribution (cf. $\eqref{estimate on the sum}$)). Redoing the scaling to go back to $M_\lambda$ we get:

\begin{equation}\label{case 2 s1 sum}
\begin{split}
\sum_{S_1} \int_0^t\int_{M_\lambda}
\left[1-\frac{m(n_1)}{m(n_2)m(n_3)m(n_4)}\right]\overline{\pi_{n_1}(\phi_1)}\pi_{n_2}(\phi_2)\overline{\pi_{n_3}(\phi_3)}\pi_{n_4}(\phi_4) dxdt\\
=\operatorname{l.c. }\frac{1}{\lambda^{2l}}\sum_{S_1}
\left[1-\frac{m(n_1)}{m(n_2)m(n_3)m(n_4)}\right]\frac{(-2)^l}{(n_1^2-n_2^2-n_3^2-n_4^2)^l}
\int_0^t\int_{M_\lambda}\overline{\pi_{n_1}(\phi_1)}\big[Q_2(\pi_{n_2}\phi_2)\overline{Q_3(\pi_{n_3}\phi_3)}\\
Q_4(\pi_{n_4}\phi_4)\big]dx\,dt.
\end{split}
\end{equation}

where $\operatorname{l.c. }$ is a shorthand for the statement ``a linear combination of" and $Q_i(f)(x):=[\tilde Q (\tilde f)](\frac{x}{\lambda})$ are differential operators of the type discussed in Corollary \ref{bilinear Xsb estimate lemma}.

Let us denote:

\begin{equation}\label{def of bar m}
\bar m (n_1,n_2,n_3,n_4)=\left[1-\frac{m(n_1)}{m(n_2)m(n_3)m(n_4)}\right]\frac{(-2)^l}{(n_1^2-n_2^2-n_3^2-n_4^2)^l}.
\end{equation}

Then the sum of $S_1$ can be written as a finite $O(1)$ linear combination of:

\begin{equation}\label{the S_1 sum}
\frac{1}{\lambda^{2l}}\sum_{S_1} \bar m(n_1,n_2,n_3,n_4) \int_0^t\int_{M_\lambda}\overline{\pi_{n_1}(\phi_1)}
Q_2(\pi_{n_2}\phi_2)\overline{Q_3(\pi_{n_3}\phi_3)}Q_4(\pi_{n_4}\phi_4)dx\,dt.
\end{equation}

The estimate for this sum is obtained in two steps:

\textbf{Step 1:}

The first step is to apply a similar analysis to that in section \ref{multilinear spectral multiplier lemma} to bound the above integral using estimates on the multiplier $\bar m$. In fact, we will be able to write $\eqref{the S_1 sum}$ as follows:

\begin{equation}\label{Step 1 equation}
\frac{1}{\lambda^{2l}}\sum_{(\alpha,\beta)\in S_1} \sum_{\theta_i \in \Z/4} A(\theta_1,\ldots,\theta_4) \int_0^t\int_{M_\lambda} P_{I_\alpha} \phi^{\theta_1}_1
Q_2(P_{I_{\beta}} \phi^{\theta_2}_2)Q_3( \phi^{\theta_3}_3)Q_4( \phi^{\theta_4}_4)dx\,dt
\end{equation}

where $A(\theta_1,\ldots,\theta_4)$ is a summable sequence with

\begin{equation}\label{estimate on the sum}
\sum_{\theta_i \in \Z/4}|A(\theta_1,\ldots,\theta_4)| \lesssim (\alpha-\beta)\frac{N_3}{N_2}\frac{1}{\left(N_2N_3(\alpha-\beta)\right)^l}
\end{equation}

 and $\phi_j^{\theta_j}$ is a frequency modulation of $\phi_j$ (see $\eqref{the tilda phis}$). In particular, $\phi_j^{\theta_j}$ has the same $X^{s,b}$ norms and frequency support as $\phi_j$.

\textbf{Step 2}
In the second step, we prove that for each fixed $(\theta_1,\ldots,\theta_4)\in (\Z/4)^4$, we have:

\begin{equation}\label{Step 2 equation}
\begin{split}
\frac{1}{\lambda^{2l}}\left|\int_0^t\int_{M_\lambda} \tilde P_{I_\alpha} \phi^{\theta_1}_1
Q_2( \phi^{\theta_2}_2)Q_3( \phi^{\theta_3}_3)Q_4(\phi^{\theta_4}_4)dx\,dt\right| \lesssim N_2^l N_3^l \frac{(N_3N_4)^{1/2}}{\lambda}\\
\times||P_{I_\alpha}\phi_1||_{X^{0,1/2+}} ||P_{I_\beta}\phi_2||_{X^{0,1/2+}}\prod_{i=3}^4||\phi_3||_{X^{0,1/2+}}.
\end{split}
\end{equation}

Combining $\eqref{Step 1 equation}$, $\eqref{estimate on the sum}$, and $\eqref{Step 2 equation}$, one gets:

\begin{align*}
\eqref{case 2 s1 sum} \lesssim &\sum_{\alpha,\beta} \left(\sum_{\theta_i \in \Z/4} |A(\theta_1,\ldots,\theta_4)|\right)(N_2N_3)^l\frac{(N_3N_4)^{1/2}}{\lambda}
||P_{I_\alpha}\phi_1||_{X^{0,1/2+}} ||P_{I_\beta}\phi_2||_{X^{0,1/2+}}\prod_{i=3}^4||\phi_i||_{X^{0,1/2+}}\\
\lesssim&\sum_{\alpha,\beta} |\alpha-\beta|\frac{N_3}{N_2}\frac{(N_2N_3)^l}{\left(N_2N_3(\alpha-\beta)\right)^l}\frac{(N_3N_4)^{1/2}}{\lambda}
||P_{I_\alpha}\phi_1||_{X^{0,1/2+}} ||P_{I_\beta}\phi_2||_{X^{0,1/2+}}\prod_{i=3}^4||\phi_i||_{X^{0,1/2+}}\\
\lesssim& \frac{N_3}{N_2}\frac{(N_3N_4)^{1/2}}{\lambda}\sum_{\alpha,\beta} \frac{1}{\left(\alpha-\beta\right)^{l-1}}
||P_{I_\alpha}\phi_1||_{X^{0,1/2+}} ||P_{I_\beta}\phi_2||_{X^{0,1/2+}}\prod_{i=3}^4||\phi_4||_{X^{0,1/2+}}.
\end{align*}

By taking $l\geq 3$, we notice that by Schur's test (for example) and the fact that $\alpha \geq \beta+8$ that

\begin{equation}\label{bound on S_1}
\eqref{case 2 s1 sum} \lesssim \frac{N_3}{N_2}\frac{(N_3N_4)^{1/2}}{\lambda}\prod_{i=1}^4||\phi_i||_{X^{0,1/2+}}
\end{equation}

which will be enough to prove $\eqref{decay of term1 in dyadic version}$. We now turn to proving $\textbf{Step 1}$ and $\textbf{Step 2}$:

\textbf{Proof of Step 1:} The analysis here is almost the same as that in the proof of lemma \ref{multilinear multiplier lemma} except for the localization to the intervals $I_\alpha$ and $I_\beta$ and the presence of the differential operators $Q_i$. We include the details for the convenience of the reader.

As in the proof of lemma \ref{multilinear multiplier lemma}, since $n_1 \in I_{\alpha}=[N_2+\alpha N_3,N_2+(\alpha+1)N_3)$, $n_2 \in I_{\beta}=[N_2+\beta N_3,N_2+(\beta+1)N_3)$, $n_3 \sim N_3$ and $n_4 \sim N_4$ we can write\footnote{with the obvious modifications if $N_3$ or $N_4$ is equal to 1, in which case we write $n_4=\tilde n_4 \in [0,2)$.}:

\begin{equation}\label{case 2 the tildas}
\begin{split}
n_1&=N_2+ \alpha N_3 +N_3 \tilde{n_1} \text{ with $\tilde{n_1} \in [0,1)$}\\
n_2&=N_2+ \beta N_3 +N_3 \tilde{n_2} \text{ with $\tilde{n_2} \in [0,1]$}\\
n_3&=N_3(1+ \tilde{n_3}) \text{ with $\tilde{n_3} \in [0,1]$}\\
n_4&=N_4(1+ \tilde{n_3}) \text{ with $\tilde{n_4} \in [0,1]$}.\\
\end{split}
\end{equation}

We now define the function $\Psi: [0,1]^4 \to \R$ given by: 

\begin{equation}\label{def of psi step 1}
\Psi(\tilde n_1,\tilde n_2,\tilde n_3,\tilde n_4)=\bar m(n_1,n_2,n_3,n_4)
\end{equation}

with $n_i$ given in $\eqref{case 2 the tildas}$. Extend $\Psi$ to a $C^{\infty}$ compactly supported function on $[-2,2]^4$ and then as a $4-$periodic function on $\R^4$. This allows us to express it in Fourier series:

$$
\Psi(\tilde n_1,\ldots,\tilde n_4)=\sum_{\theta_i \in \Z/4} e^{i (\theta_1 \tilde n_1+\ldots+\theta_4 \tilde n_4)}A(\theta_1,\ldots,\theta_4)
$$

with 

$$
\sum |A(\theta_1, \ldots,\theta_4)|\lesssim ||\Psi||_{C^2([0,1]^4)}.\footnote{See for instance \cite{Gr} Theorem 3.2.16.}
$$

Notice that this already gives $\eqref{Step 1 equation}$. In fact, 

\begin{align*}
\eqref{the S_1 sum}=&\frac{1}{\lambda^{2l}}\sum_{(\alpha,\beta)\in S_1}\sum_{\substack{n_1 \in I_\alpha, n_2 \in I_\beta\\ n_3 \sim N_3,n_4 \sim N_4}}\bar m(n_1,n_2,n_3,n_4) \int_0^t \int_{M_{\lambda}}\overline{\pi_{n_1} \phi_1}Q_2( \pi_{n_2}\phi_2)\overline{Q_3( \pi_{n_3}\phi_3)}Q_4( \pi_{n_4} \phi_4)dxdt\\ 
=&\frac{1}{\lambda^{2l}} \sum_{(\alpha,\beta)\in S_1}\sum_{\substack{n_1 \in I_\alpha, n_2 \in I_\beta\\ n_3 \sim N_3,n_4 \sim N_4}}\sum_{\theta_i \in \Z/4}A(\theta_1,\ldots,\theta_4) \\
&\times\int_0^t \int_{M_{\lambda}}e^{i\theta_1 \tilde n_1}\overline{\pi_{n_1} \phi_1}Q_2( e^{i\theta_2 \tilde n_2}\pi_{n_2}\phi_2)\overline{Q_3( e^{i\theta_3 \tilde n_3}\pi_{n_3}\phi_3)}Q_4( e^{i\theta_4 \tilde n_4}\pi_{n_4}\phi_4)dxdt\\ 
=&\frac{1}{\lambda^{2l}}\sum_{(\alpha,\beta)\in S_1} \sum_{\theta_i \in \Z/4}A(\theta_1,\ldots,\theta_4)\int_0^t \int_{M_\lambda}
\overline{P_{I_\alpha} \phi^{\theta_1}_1}
Q_2(P_{I_\beta}\phi^{\theta_2}_2)\overline{Q_3( \phi^{\theta_3}_3)}
Q_4(\phi^{\theta_4}_4))dxdt
\end{align*}

where

\begin{equation}\label{the tilda phis}
\begin{split}
\phi_1^{\theta_1}&=\sum_{n_1}e^{-i\theta_1 (n_1-(N_2+\alpha N_3))\over N_3}\pi_{n_1}\phi_1\\
\phi_2^{\theta_2}&=\sum_{n_2}e^{i\theta_2 (n_2-(N_2+\beta N_3))\over N_3}\pi_{n_2}\phi_2\\
\phi_3^{\theta_3}&=\sum_{n_3}e^{-i\theta_3 (n_3-N_3)\over N_3}\pi_{n_3}\phi_3\\
\phi_4^{\theta_4}&=\sum_{n_4}e^{i\theta_4 (n_4-N_4)\over N_4}\pi_{n_4}\phi_4.\\
\end{split}
\end{equation}

In order to prove $\eqref{estimate on the sum}$ and finish the proof of \textbf{Step 1} all we need to do is prove that:

\begin{equation}\label{C^2 norm of Psi}
||\Psi||_{C^2([0,1]^4)} \lesssim (\alpha-\beta)\frac{N_3}{N_2}\frac{1}{\left(N_2N_3(\alpha-\beta)\right)^l}.
\end{equation}

This estimate follows by direct verification using the following facts:

\begin{enumerate}
\item By the mean value theorem we have that 

\begin{equation}\label{mean value theorem}
\left|1-\frac{m(n_1)}{m(n_2)}\right| \sim \frac{|n_1-n_2|}{N_2} \sim (\alpha-\beta)\frac{N_3}{N_2}.
\end{equation}

\item Similarly
\begin{equation}\label{estimate on derivatives}
\left|\frac{\partial^\gamma}{\partial(\tilde n_1,\tilde n_2)^\gamma}\left(1-\frac{m(N_2+\alpha N_3 +N_3\tilde n_1)}{m(N_2+\beta N_3 +N_3 \tilde n_2)}\right)\right| \lesssim (\frac{N_3}{N_2})^{|\gamma|}\ll (\alpha-\beta)\frac{N_3}{N_2}
\end{equation}

for any multi-index $\gamma$ with $|\gamma|\leq 2$. This follows from the easily verified fact that:

$$
\frac{\partial^l m(\xi)}{m(\xi)}\lesssim \frac{1}{|\xi|^l}
$$
for any $\xi \in \R$ and $l\in \N$ (see $\eqref{def of m_0}$)

\item 

A calculation shows that:

\begin{align*}
n_1^2-(n_2^2+n_3^2+n_4^2)=&2N_2N_3\left(\alpha-\beta+(\tilde n_1-\tilde n_2) \right)+\left(\alpha^2-\beta^2\right)N_3^2+N_3^2(\tilde n_1^2-\tilde n_2^2)\\
&+2N_3^2\left(\alpha \tilde n_1-\beta \tilde n_2\right)-N_3^2(1+\tilde n_3)^2-N_4^2(1+\tilde n_4)^2\\
=&\left(2N_2N_3\left(\alpha-\beta\right)+\left(\alpha^2-\beta^2\right)N_3^2 \right)\psi_2(\tilde n_1,\tilde n_2,\tilde n_3,\tilde n_4)
\end{align*}
with $\psi_2: [0,1]^4 \to \R$ being $C^\infty$, bounded below by $1/2$ (since in the $S_1$ sum we are estimating $\alpha-\beta \geq 8$), and bounded above along with its derivatives by $O(1)$.

This gives that 
\begin{align*}
\left|\frac{\partial^\gamma}{\partial(\tilde n_1,\tilde n_2,\tilde n_3,\tilde n_4)^\gamma}\left(\frac{2^l}{(n_1^2-n_2^2-n_3^2-n_4^2)^l}\right)\right| &\lesssim_l \frac{1}{\left(2N_2N_3 (\alpha-\beta)+\left(\alpha^2-\beta^2\right)N_3^2\right)^l}\\
&\lesssim \frac{1}{\left(N_2N_3 (\alpha-\beta)\right)^l}.
\end{align*}

\end{enumerate}

Combining the above three facts one gets $\eqref{C^2 norm of Psi}$.

\textbf{Proof of Step 2} The proof of $\eqref{Step 2 equation}$ is an easy consequence of Corollary \ref{bilinear Xsb estimate lemma}. In fact,

\begin{align*}
\operatorname{L.H.S.} \eqref{Step 2 equation}\lesssim&
\frac{1}{\lambda^{2l}}||P_{I_\alpha}\phi_1^{\theta_1} Q_4(\phi_4^{\theta_4})||_{L^2_{t,x}([0,t]\times M_\lambda)}||Q_2(P_{I_\beta}\phi_2^{\theta_2}) Q_3(\phi_3^{\theta_3})||_{L^2_{t,x}([0,t]\times M_\lambda)}\\
\lesssim& N_2^l N_3^l \frac{(N_3N_4)^{1/2}}{\lambda}
\times||P_{I_\alpha}\phi_1||_{X^{0,1/2+}} ||P_{I_\beta}\phi_2||_{X^{0,1/2+}}\prod_{i=3}^4||\phi_i||_{X^{0,1/2+}}
\end{align*}

where we used $\eqref{bilinear Xsb estimate with differential operators}$ in the second inequality along with the fact that $\deg(Q_2)+\deg(Q_3)+\deg(Q_4)\leq 2l$ and $0\leq \deg(Q_i)\leq l$ for $i=2,3,4$.

Similarly, one gets the same bound for $S_2$ just by interchanging the roles of $N_1$ and $N_2$ above. The bound for $S_3=\{I_\alpha, I_\beta \textrm{ are at a distance $\leq 8N_3$ apart}\}$ (near diagonal terms) is simpler. For each $I_\beta$, let $S_3^\beta$ be the set of $I_\alpha$ intervals that are at a distance less than $\leq 8N_3$ from $I_\beta$. Clearly there are $\leq 19$ elements in $S_3^\beta$:

\begin{align*}
&\sum_{S_3} \int_0^t\int_{M_\lambda}
\left[1-\frac{m(n_1)}{m(n_2)m(n_3)m(n_4)}\right]\overline{\pi_{n_1}(\phi_1)}\pi_{n_2}(\phi_2)\overline{\pi_{n_3}(\phi_3)}\pi_{n_4}(\phi_4) dx\,dt\\
&=\sum_{I_\beta}\sum_{I_\alpha\in S_3^\beta} \sum_{\substack{n_1 \in I_\alpha, n_2 \in I_\beta\\ n_3 \sim N_3,n_4 \sim N_4}} \int_0^t\int_{M_\lambda}
\left[1-\frac{m(n_1)}{m(n_2)m(n_3)m(n_4)}\right]\overline{\pi_{n_1}(\phi_1)}\pi_{n_2}(\phi_2)\overline{\pi_{n_3}(\phi_3)}\pi_{n_4}(\phi_4) dx\,dt.
\end{align*}

The sum

$$
\sum_{\substack{n_1 \in I_\alpha, n_2 \in I_\beta\\ n_3 \sim N_3,n_4 \sim N_4}} \int_0^t\int_{M_\lambda}
\left[1-\frac{m(n_1)}{m(n_2)m(n_3)m(n_4)}\right]\overline{\pi_{n_1}(\phi_1)}\pi_{n_2}(\phi_2)\overline{\pi_{n_3}(\phi_3)}\pi_{n_4}(\phi_4) dx\,dt
$$

is bounded using lemma \ref{multilinear multiplier lemma}. In fact, with 

$$
\bar m(\xi_1,\xi_2,\xi_3,\xi_4)=\frac{N_2}{N_3}\left[1-\frac{m(n_1)}{m(n_2)m(n_3)m(n_4)}\right]
$$

condition $\eqref{symbol estimates on bar m}$ was already verified in $\eqref{mean value theorem}$ and $\eqref{estimate on derivatives}$. The second condition $\eqref{conditional estimate}$ in lemma \ref{multilinear multiplier lemma} follows from $\eqref{bilinear estimate}$ since:

\begin{align*}
\int_0^t\int_{M_\lambda} f_1(t,x) f_2(t,x) f_3(t,x) f_4(t,x) dxdt\lesssim ||f_1 f_4||_{L^2_{t,x}([0,t]\times M_\lambda)}||f_2f_3||_{L^2_{t,x}([0,t]\times M_\lambda)}\\
\lesssim \frac{(N_3N_4)^{1/2}}{\lambda}\prod_{i=1}^4||f_i||_{X^{0,1/2+}}
\end{align*}

whenever $f_1,\ldots, f_4$ have the same frequency localizations as those under consideration.
Therefore,

\begin{equation}\label{bound on S_3}
\begin{split}
&\left|\sum_{S_3} \int_0^t\int_{M_\lambda}
\left[1-\frac{m(n_1)}{m(n_2)m(n_3)m(n_4)}\right]\overline{\pi_{n_1}(\phi_1)}\pi_{n_2}(\phi_2)\overline{\pi_{n_3}(\phi_3)}\pi_{n_4}(\phi_4) dxdt\right|\\
&\lesssim\sum_{I_\beta}\sum_{I_\alpha\in S_3^\beta} \frac{N_3}{N_2}\frac{(N_3N_4)^{1/2}}{\lambda}||P_{I_\alpha}\phi_1||_{X^{0,1/2+}}||P_{I_\beta}\phi_2||_{X^{0,1/2+}}\prod_{i=3}^4||\phi_i||_{X^{0,1/2+}}\\
&\lesssim \frac{N_3}{N_2}\frac{(N_3N_4)^{1/2}}{\lambda} \prod_{i=1}^4||\phi_i||_{X^{0,1/2+}}
\end{split}
\end{equation}

where we applied lemma \ref{multilinear multiplier lemma} in the first inequality and used Cauchy-Schwarz in the $I_\beta$ sum and the fact that $|S_3^\beta|\lesssim 1$ in the second.

Combining \eqref{bound on S_1} and $\eqref{bound on S_3}$ we get:

\begin{equation}\label{case 2 total contribution}
\begin{split}
\text{L.H.S. of \eqref{decay of term1 in dyadic version}}\lesssim& \frac{N_3}{N_2}\frac{N_1}{N_2N_3N_4}\frac{(N_3N_4)^{1/2}}{\lambda} ||\phi_1||_{X^{-1,1/2+}}\prod_{i=2}^4||\phi_i||_{X^{1,1/2+}}\\
&\lesssim \frac{1}{\lambda N^{1/2-}N_2^{0+}}||\phi_1||_{X^{-1,1/2+}}\prod_{i=2}^4||\phi_i||_{X^{1,1/2+}}
\end{split}
\end{equation}

which is what we want.

\textbf{Case 3. $N_2 \sim N_3 \gtrsim N$} 

The bound is obtained by invoking lemma \ref{multilinear multiplier lemma}. In fact, by $\eqref{the N_2/lambda decay}$ we have:

\begin{equation}\label{conditional estimate case 3}
\begin{split}
|\int_0^t\int_M f_1(t,x)\ldots f_4(t,x)dxdt| \lesssim ||f_2 f_1||_{L^2_{t,x}([0,\delta]\times M)}||f_3 f_4||_{L^2_{t,x}([0,\delta]\times M)}\\\lesssim \left(\frac{N_1}{\lambda}\right)^{1/2}\left(\frac{N_4}{\lambda}\right)^{1/2}\prod_{i=1}^4||f_i||_{X^{0,1/2+}}
\end{split}
\end{equation}

whenever $f_i$ are spectrally localized to frequencies $\sim N_i$, $N_1\lesssim N_2$, and  $N_4\leq N_3$. This verifies the conditional estimate $\eqref{conditional estimate}$ of lemma \ref{multilinear multiplier} with $B=\frac{(N_1N_4)^{1/2}}{\lambda}$ and $Y=X^{0,1/2+}([0,\delta]\times M)$. 

The symbol type conditions $\eqref{symbol estimates on bar m}$ are satisfied by the multiplier:
$$
\bar m := \left(\frac{m(N_1)}{m(N_2)m(N_3)m(N_4)}\right)^{-1}\left[1-\frac{m(n_1)}{m(n_2)m(n_3)m(n_4)}\right]
$$

since 

$$
\left|1-\frac{m(n_1)}{m(n_2)m(n_3)m(n_4)}\right| \lesssim  \frac{m(n_1)}{m(n_2)m(n_3)m(n_4)}
\sim \frac{m(N_1)}{m(N_2)m(N_3)m(N_4)}.
$$

This follows from the fact that $\frac{m(n_1)}{m(n_2)}\gtrsim 1$ since $N_2 \gtrsim N_1$ and $m$ is non-increasing. Estimates of the derivatives follow as well since $m$ itself satisfies the symbol-type estimate:

\begin{equation}\label{derivative of m}
\frac{d^l}{dt^l}m(t)\lesssim \frac{m(t)}{t^l}
\end{equation}

which is due to $\eqref{def of m_0}$\footnote{As mentioned in the footnote to lemma \ref{multilinear multiplier lemma}, one only needs to verify the symbol estimates for two derivatives only.}. As a result,

\begin{equation}\label{Case 3 split}
\begin{split}
\operatorname{L.H.S. of }\eqref{decay of term1 in dyadic version}\lesssim&
\frac{m(N_1)}{m(N_2)m(N_3)m(N_4)} \frac{(N_1N_4)^{1/2}}{\lambda}\prod_{i=1}^4||\phi_i||_{X^{0,1/2+}([0,\delta]\times M)}\\
\lesssim&
\frac{m(N_1)}{m(N_2)m(N_3)m(N_4)} \frac{N_1}{N_2N_3N_4}\frac{(N_1N_4)^{1/2}}{\lambda}||\phi_1||_{X^{-1,1/2+}([0,\delta]\times M)}\prod_{i=2}^4||\phi_i||_{X^{1,1/2+}([0,\delta]\times M)}.\\
\end{split}
\end{equation}

We will estimate $\frac{m(N_1)}{m(N_2)m(N_3)m(N_4)}$ depending on the relative position of $N_1$ relative to $N$:

\textbf{Case 3a $N_2 \sim N_3$, $N_1\gg N$}: In this case, we use the explicit expressions for $m(N_1), m(N_2),m(N_3)$ to get:

\begin{equation}\label{Case 3a split}
\begin{split}
\operatorname{L.H.S. of }\eqref{decay of term1 in dyadic version}\lesssim&
\frac{1}{\lambda} \frac{N^{1-s}N_1^{-(1-s)}}{N^{2(1-s)}N_2^{-(1-s)}N_3^{-(1-s)}}\frac{N_1^{3/2}}{N_2N_3m(N_4)N_4^{1/2}}||\phi_1||_{X^{-1,1/2+}}\prod_{i=2}^4||\phi_i||_{X^{1,1/2+}}\\
\lesssim& \frac{1}{\lambda} \frac{N_1^{1/2+s}}{N^{1-s}N_2^{2s}m(N_4)N_4^{1/2}}||\phi_1||_{X^{-1,1/2+}([0,\delta]\times M)}\prod_{i=2}^4||\phi_i||_{X^{1,1/2+}([0,\delta]\times M)}\\
\lesssim& \frac{1}{\lambda} \frac{1}{N^{1-s}N_2^{s-1/2}}||\phi_1||_{X^{-1,1/2+}([0,\delta]\times M)}\prod_{i=2}^4||\phi_i||_{X^{1,1/2+}([0,\delta]\times M)}\\
\lesssim& \frac{1}{\lambda N^{1/2-}N_2^{0+}}||\phi_1||_{X^{-1,1/2+}([0,\delta]\times M)}\prod_{i=2}^4||\phi_i||_{X^{1,1/2+}([0,\delta]\times M)}
\end{split}
\end{equation}
where we have used the fact that $N_2\sim N_3$ and $N_1 \lesssim N_2$ in the second and third inequalities respectively, and that $m(N_4)N_4^{1/2-}\gtrsim 1$ for $s>\frac{1}{2}$ in the third inequality as well.

\textbf{Case 3b: $N_2\sim N_3 \gtrsim N$, $N_1\lesssim N$} In this case, bound $m(N_1)$ by 1 and use the explicit expressions for $m(N_2)$ and $m(N_3)$ to get:

Here the key estimate goes as:

\begin{equation}\label{Case 3b split}
\begin{split}
\operatorname{L.H.S. of }\eqref{decay of term1 in dyadic version}\lesssim&
\frac{1}{\lambda} \frac{1}{N^{2(1-s)}N_2^{-(1-s)}N_3^{-(1-s)}}\frac{N_1^{3/2}}{N_2N_3m(N_4)N_4^{1/2}}||\phi_1||_{X^{-1,1/2+}}\prod_{i=2}^4||\phi_i||_{X^{1,1/2+}}\\
\lesssim& \frac{1}{\lambda} \frac{N_1^{3/2}}{N^{2(1-s)}N_2^{2s}m(N_4)N_4^{1/2}}||\phi_1||_{X^{-1,1/2+}([0,\delta]\times M)}\prod_{i=2}^4||\phi_i||_{X^{1,1/2+}([0,\delta]\times M)}\\
\lesssim& \frac{1}{\lambda N^{1/2}} \frac{N^{2s}}{N_2^{2s}}||\phi_1||_{X^{-1,1/2+}([0,\delta]\times M)}\prod_{i=2}^4||\phi_i||_{X^{1,1/2+}([0,\delta]\times M)}\\
\lesssim& \frac{1}{\lambda N^{1/2-}N_2^{0+}}
||\phi_1||_{X^{-1,1/2+}([0,\delta]\times M)}\prod_{i=2}^4||\phi_i||_{X^{1,1/2+}([0,\delta]\times M)}
\end{split}
\end{equation}

where we have used the fact that $N_2\sim N_3$ and $N_1 \lesssim N$ in the second and third inequalities respectively, and that $m(N_4)N_4^{1/2-}\gtrsim 1$ for $s>\frac{1}{2}$ in the third inequality as well.

\textbf{Case 4 $N_2\gg N_3\gtrsim N$:} 

Recall that in this case $N_1 \sim N_2$. We will obtain the bound by invoking lemma \ref{multilinear multiplier} and the same bound for the multiplier as in \textbf{Case 3}. However, in this case, since $N_3\ll N_1$, $N_3$ should replace $N_1$ on the left hand side of  $\eqref{conditional estimate case 3}$ :

\begin{equation}\label{Case 4 split}
\begin{split}
\operatorname{L.H.S. of }\eqref{decay of term1 in dyadic version}\lesssim&
\frac{m(N_1)}{m(N_2)m(N_3)m(N_4)} \frac{(N_3N_4)^{1/2}}{\lambda}\prod_{i=1}^4||\phi_i||_{X^{0,1/2+}([0,\delta]\times M)}\\
\lesssim&
\frac{1}{\lambda}\frac{1}{m(N_3)N_3^{1/2}m(N_4)N_4^{1/2}} ||\phi_1||_{X^{-1,1/2+}([0,\delta]\times M)}\prod_{i=2}^4||\phi_i||_{X^{1,1/2+}([0,\delta]\times M)}\\
\lesssim&
\frac{1}{\lambda}\frac{1}{m(N_3)N_3^{1/2}} ||\phi_1||_{X^{-1,1/2+}([0,\delta]\times M)}\prod_{i=2}^4||\phi_i||_{X^{1,1/2+}([0,\delta]\times M)}\\
\lesssim& \frac{1}{\lambda}\frac{1}{N^{1-s}N_3^{\frac{1}{2}-(1-s)}} ||\phi_1||_{X^{-1,1/2+}([0,\delta]\times M)}\prod_{i=2}^4||\phi_i||_{X^{1,1/2+}([0,\delta]\times M)}\\
\lesssim& \frac{1}{\lambda N^{1/2-}N_3^{0+}}||\phi_1||_{X^{-1,1/2+}([0,\delta]\times M)}\prod_{i=2}^4||\phi_i||_{X^{1,1/2+}([0,\delta]\times M)}\\
\end{split}
\end{equation}

where we used the fact that $N_2\sim N_1$ in the second inequality, that $m(N_4)N^{1/2-}\gtrsim 1$ (since $s>1/2$) in the third inequality, the explicit formula for $m(N_3)$ in the forth, and finally the fact that $N_3\gtrsim N$ in the fifth.

In this case, we do not have an exponential decay factor of the form $N_2^{0+}$ in the denominator. However, since $N_1\sim N_2$ on can use Cauchy-Schwarz to sum in $N_1$ and $N_2$ and then use the geometric decay factor $N_3^{0+}$ to sum in $N_3$ and $N_4$.

This finishes the proof of $\eqref{decay of term1 in dyadic version}$ and hence the bound on $\operatorname{Term}_1$

\subsubsection{Bound on $\operatorname{Term}_2$:}

We now turn to proving the decay estimate for $\operatorname{Term}_2$ in $\eqref{definition of term2}$. We will be able to prove better decay for $\operatorname{Term}_2$ than that for $\operatorname{Term}_1$. In fact we will show:
\begin{equation}\label{decay of term 2}
\operatorname{Term_2}\lesssim \frac{1}{\lambda N^{1-}}||Iu||_{X^{1,1/2+}}^6.
\end{equation}

This decay is proved by exploiting the high multiplicity of factors more than the exact shape of the multiplier (though of course the existence of the multiplier is crucial to restrict attention to the case where at least one factor has high frequency). In fact, since the specific form of the multiplier is not as important in this case, we will rewrite $\operatorname{Term}_2$ in a form that is more convenient to do a multilinear analysis in $L^p$ spaces:

\begin{align*}
\operatorname{Term}_2=&\sum_{n_i}\int_0^t\int_{M_\lambda}
\left[1-\frac{m(n_1)}{m(n_2)m(n_3)m(n_4)}\right]\overline{\pi_{n_1}\left(I\left(|u|^2u\right)\right)}\pi_{n_2}(Iu)\overline{\pi_{n_3}(Iu)}\pi_{n_4}(Iu)\,dx\,dt\\
=&\int_0^t\int_{M_\lambda}
I\left(|u|^2u\right)Iu\overline{Iu}Iudxdt-\int_0^t\int_{M_\lambda}
I^2\left(|u|^2u\right)u\,\overline{u}\,u\,dx\,dt.\\
\end{align*}

We now write $u=\sum_{N}u_N$ in Littlewood-Paley pieces as before to get:

\begin{equation}\label{term 2 dyadic}
\begin{split}
\operatorname{Term}_2
=\sum_{N_i \in 2^{\N}}\left(\int_0^t\int_{M_\lambda}
I\left(u_{N_1}\overline{u_{N_2}}u_{N_3}\right)Iu_{N_4}\overline{Iu_{N_5}}Iu_{N_6}\,dx\,dt \right.\\
\left .-\int_0^t\int_{M_\lambda}
I^2\left(u_{N_1}\overline{u_{N_2}}u_{N_3}\right)u_{N_4}\overline{u_{N_5}}u_{N_6}\,dx\,dt\right).
\end{split}
\end{equation}

Let us denote by $N_{max}$ the largest of all $N_1\ldots N_6$ and $N_{med}$ the second largest. As a result of the eigenvalue localization lemma \ref{BGT localization}, and the fact that $I^2 u_{N_i}=I u_{N_i}=u_{N_i}$ if ${N_i}\leq N/8$, the contribution of the terms for which $N_{max}\ll N$ to $\eqref{term 2 dyadic}$ is bounded by $N^{-q}\prod_{i=1}^{6}||Iu||_{X^{1,1/2+}}$ for some large enough $q$. In fact, writing $u_{N_1}\overline{u_{N_3}}u_{N_3}=P_{\leq N}\left(u_{N_1}\overline{u_{N_3}}u_{N_3}\right)+P_{>N}\left(u_{N_1}\overline{u_{N_3}}u_{N_3}\right)$, we get since $I\left(P_{\leq N}u_{N_1}\overline{u_{N_3}}u_{N_3}\right)=I^2\left(P_{\leq N}u_{N_1}\overline{u_{N_3}}u_{N_3}\right)$ that the contribution of the range $N_{max}\ll N$ to $\eqref{term 2 dyadic}$ consists of:

$$
\sum_{N_i \ll N}\left(\int_0^t\int_{M_\lambda}
IP_{>N}\left(u_{N_1}\overline{u_{N_2}}u_{N_3}\right)u_{N_4}\overline{u_{N_5}}u_{N_6}\,dx\,dt-\int_0^t\int_{M_\lambda}
I^2P_{>N}\left(u_{N_1}\overline{u_{N_2}}u_{N_3}\right)u_{N_4}\overline{u_{N_5}}u_{N_6}\,dx\,dt\right).
$$

Estimating $||IP_{>N}\left(u_{N_1}\overline{u_{N_2}}u_{N_3}\right)||_{L_x^2(M)}$ and $\|I^2P_{>N}\left(u_{N_1}\overline{u_{N_2}}u_{N_3}\right)\|_{L_x^2(M)}$ by $N^{-C}\prod_{i=1}^3||u_{N_i}||_{L_x^2(M)}$ for some large $C$ and using H\"older, a crude Sobolev embedding for the other terms, and the embedding $X^{0,1/2+}\subset L_t^\infty L_x^2$ one gets that the contribution of the range $N_{max}\ll N$ is harmless.

A similar argument using lemma \ref{BGT localization} shows that the contribution of $N_{max}\gg N_{med}$ is also harmless so we restrict attention to the case when $N_{max}\gtrsim N$ and $N_{max}\sim N_{med}$.

Using the fact that $I$ is bounded as an $L^p$ multiplier (cf. corollary 4.3.2 of \cite{Sogge}), one can estimate the contribution of $N_{max}\gtrsim N$ using H\"older's inequality as follows:

\begin{align*}
\operatorname{Term}_2^{N_{max}\gtrsim N}\lesssim& \sum_{\substack{N_{max}, N_{med}\gtrsim N,\\ \tilde N_i\leq N_{max}}}||u_{N_{max}}||_{L_{t,x}^4}||u_{N_{med}}||_{L^4_{t,x}}\prod_{i=1}^4||u_{\tilde N_i}||_{L_{t,x}^8}\\
\lesssim& \sum_{\substack{N_{max}, N_{med}\gtrsim N,\\ \tilde N_i\leq N_{max}}} \frac{N_{max}^{1/4}N_{med}^{1/4}}{\lambda^{1/4}\lambda^{1/4}}||u_{N_{max}}||_{X^{0,1/2+}}||u_{N_{med}}||_{X^{0,1/2+}}\prod_{i=1}^4\frac{\tilde N_i^{5/8}}{\lambda^{1/8}}||u_{\tilde N_i}||_{X^{0,1/2+}}
\end{align*}

by the $\eqref{L^4 estimate Xsb}$ and $\eqref{L^8 estimate Xsb}$\footnote{Here we used the second part of $\eqref{L^8 estimate Xsb}$ corresponding to the case when $\tilde N_i \geq \lambda$ ($i=3,\ldots, 6$). The cases when $\tilde N_i \leq \lambda$ are treated similarly and yield the bound $\frac{1}{\lambda^{1/2} N^{3/2-}} \leq \frac{1}{\lambda N^{1-}}$ since $\lambda \leq N$.}. As a result, we have:

\begin{align*}
\operatorname{Term}_2^{N_{max}\gtrsim N}
\lesssim\frac{1}{\lambda} \sum_{N_{max} \in 2^{\N}, N_{med}, \tilde N_i\leq N_{max}} \frac{N_{max}^{-3/4}N_{med}^{-3/4}}{m(N_{max})m(N_{med})}||Iu_{N_{max}}||_{X^{1,1/2+}}||Iu_{N_{med}}||_{X^{1,1/2+}}\\
\times \prod_{i=1}^4\frac{\tilde N_i^{-3/8}}{m(\tilde N_i)}||Iu_{\tilde N_i}||_{X^{1,1/2+}}.
\end{align*}

Since $N_{max}$ and $N_{med}$ are $\gtrsim N$, it is easy to see that for $\alpha> 1-s$, $N_{max}^{\alpha}m(N_{max})\gtrsim N^{\alpha}$ and $N_{med}^{\alpha} m(N_{med})\gtrsim N^{\alpha}$. Similarly, for any $k>0$, $m(k)k^{\beta}\gtrsim 1$ if $\beta > 1-s$. Applying those estimates with $\alpha=3/4->1/3>1-s$ and $\beta=3/8->1/3>1-s$, we get:

\begin{align*}
\operatorname{Term}_2^{N_{max}\gtrsim N}
\lesssim&\frac{1}{\lambda N^{3/2-}} \sum_{N_{max} \in 2^{\N}, N_{med}, \tilde N_i\leq N_{max}} N_{max}^{0-}||Iu_{N_{max}}||_{X^{1,1/2+}}||Iu_{N_{med}}||_{X^{1,1/2+}}\prod_{i=1}^4||Iu_{\tilde N_i}||_{X^{1,1/2+}}\\
\lesssim& \frac{1}{\lambda N^{3/2-}}||Iu||_{X^{1,1/2+}}^6\leq \frac{1}{\lambda N^{1-}}||Iu||_{X^{1,1/2+}}^6
\end{align*}

since $\lambda \ll N$.
\endproof

\subsection{Polynomial bounds on $\tilde{E}[u(t)]$ and global well-posedness:}\label{polybounds}

We are now ready for the final step of the argument. Suppose that $U_0 \in H^s(M)$, then consider the function $u_0: M_\lambda \to \C$ given by $u_0(x)=\frac{1}{\lambda}U_0(\frac{x}{\lambda})$. Then $||u_0||_{\dot H^s(M_\lambda)}=\frac{1}{\lambda^s}||U_0||_{\dot H^s(M)}$ and
\begin{align*}
\tilde E[u_0]=&\frac{1}{2}\int_{M_\lambda}|\nabla Iu_0|^2 dv_\lambda+\frac{1}{4}\int_{M_{\lambda}}|Iu_0|^4 dv_\lambda.\\
\end{align*}
Using the fact, that $u_0(x)=\frac{1}{\lambda}U_0(\frac{x}{\lambda})=\frac{1}{\lambda}\sum_{\nu_j}\pi_{\nu_j}U_0(\frac{x}{\lambda})$, we calculate:

\begin{align*}
\int_{M_\lambda}|\nabla Iu_0|^2 dv_\lambda=& \langle -\Delta_\lambda Iu_0, Iu_0\rangle_{L^2(M_\lambda)}=\langle -\Delta_\lambda I^2u_0, u_0\rangle_{L^2(M_\lambda)}\\
=&\sum_{\nu_j}\frac{\nu_j}{\lambda^2}\left(m_0\left(\frac{1}{N}\sqrt{\frac{\nu_j}{\lambda^2}}\right)\right)^2  ||\pi_{\frac{\nu_j}{\lambda^2}}u_0||_{L^2(M_\lambda)}^2
=\sum_{\nu_j}\frac{\nu_j}{\lambda^2}m_0\left(\frac{\sqrt{\nu_j}}{\lambda N}\right)^2  ||\pi_{\nu_j}U_0||_{L^2(M)}^2
\end{align*}
where we have used in the last inequality that $||\pi_{\frac{\nu_j}{\lambda^2}}u_0||_{L^2(M_\lambda)}^2=||\pi_{\nu_j}U_0||_{L^2(M)}^2$
which comes from \eqref{rescaling projections} and the fact that $\pi^\lambda_{\nu/\lambda^2}u_0(x)=\frac{1}{\lambda}\pi_{\nu}U_0(\frac{x}{\lambda})$. Splitting the sum in two cases $\nu_j \leq 2(\lambda N)^2$ and $\nu_j \geq 2(\lambda N)^2$ get that:
$$
\int_{M_\lambda}|\nabla Iu_0|^2 dv_\lambda \lesssim \frac{N^{2(1-s)}}{\lambda^{2s}} ||U_0||_{H^s(M)}^2.
$$

Using the Gagliardo-Nirenberg inequality (which is scale-invariant) we get:

$$
||Iu||_{L^4(M_\lambda)}\lesssim ||Iu||_{\dot H^1(M_\lambda)}^{1/2}||Iu||_{L^2(M_\lambda)}^{1/2}\lesssim  ||Iu||_{\dot H^1(M_\lambda)}^{1/2}||U_0||_{L^2(M)}^{1/2} 
$$
(the scaling leaves the $L^2$norm dimensionless). As a result, one gets:

\begin{equation}\label{rescaling modified energy}
\tilde E [u_0] \lesssim \frac{N^{2(1-s)}}{\lambda^{2s}}||U_0||_{H^s(M)}^2.
\end{equation}

Therefore, choosing $\lambda \sim N^{\frac{1-s}{s}}$ ensures that $\tilde E[u(0)] \leq \frac{1}{2}$. By Proposition \ref{LWP of I system}
there exists a $\delta >0$ such that:

$$
\tilde E[u(\delta)]\leq \tilde E[u(0)] + O\left(\frac{1}{\lambda N^{1/2-}}\right).
$$

We can keep on using Proposition \ref{LWP of I system} to iterate the above inequality as long as $\tilde E[u(k\delta)] \lesssim 1$, and hence we can repeat the above procedure $\sim \lambda N^{1/2-}$ times.
Given $T \gg 1$, we choose $N=N(T)\gg 1$ so that
\begin{equation}\label{T in terms of N}
T\sim \frac{\delta \lambda N^{1/2-}}{\lambda^2}\sim
\frac{N^{1/2-}}{\lambda}\sim N^{\frac{3s-2}{2s}-}.
\end{equation}

Since the exponent is positive for $s>\frac{2}{3}$, $N \gg 1$ is defined for all times $T \gg 1$. Moreover, for $0\leq t \leq T$, we have:

$$
||U(t)||_{H^s(M)}=\lambda^s||u(\lambda^2t)||_{H^s(M_\lambda)}\lesssim \lambda^s \tilde E[u(\lambda^2t)]^{1/2}\lesssim \lambda^s\sim N^{1-s}\|U_0\|_{H^s(M)}
$$
since for $0\leq t\leq T$ we have $\tilde E[u(\lambda^2t)]\lesssim 1$. Using \eqref{T in terms of N} we get:

\begin{equation}\label{polynomial bound from I method}
||U(t)||_{H^s(M)}\lesssim T^{\frac{2s(1-s)}{3s-2}+}\|U_0\|_{H^s(M)}
\end{equation}

finishing the proof of Theorem \ref{main theorem}.

\endproof

\section{Spectral localization: Part II}\label{part 2}

In this section, we give a proof of the results of Section \ref{part 1}. Here $M^d$ is a compact Riemannian manifold (possibly with boundary) with Riemannian metric $g_{\alpha \beta}$ and $\Delta=\Delta_g$ is the Laplace-Beltrami operator. In what follows, we assume that $e_i$ are eigenfunctions of Laplace-Beltrami operator with Dirichlet or Neumann boundary conditions corresponding the the eigenvalues $n_i^2$. The calculation leading to the proof of Theorem \ref{A_0 and A_n theorem} and Corollary \ref{spectral cluster lemma 1} is based on Bochner-type calculations involving the Ricci commutation identities used to commute covariant derivatives. This is basically the way to generalize some integration by parts manipulations to the case when the involved functions are contracted tensors rather than just functions. We will do calculations without resorting to a prefered coordinate system. The notation we will use is fairly standard: We use abstract index notation\footnote{in particular the indices do not correspond to any preferred coordinate system.}, repeated indices are summed, and $g$ is used to raise and lower indices. We will quickly review some elementary concepts from Riemannian geometry for the sole purpose of fixing notation. For a more comprehensive treatment we refer the reader to any of the standard texts on the subject \cite{Petersen}\cite{GHL} or the first chapter of \cite{CLN} for a review of the commutation identities we will use.

\subsection{Review of some elementary Riemannian geometry concepts} We denote by $\nabla$ the Riemannian connection associated to $g$. This is the unique connection that is torsion free and for which the metric $g_{\alpha \beta}$ is parallel. In other words, 

$$
\nabla g_{\alpha \beta}=0 \textrm{ (Metric Compatibility)}
$$

and for any vector fields $X$ and $Y$:

$$
\nabla_XY-\nabla_Y X=[X,Y] \textrm{ (Torsion Free)}
$$

where $[X,Y]=XY-YX$ is the Lie Bracket.

\textbf{Ricci commutation formulas:} For any $C^{\infty}(M)$ function $f$, it follows from torsion-freeness that:

$$
\nabla_\alpha \nabla_\beta f=\nabla_\beta \nabla_\alpha f. 
$$

Two covariant derivatives only commute when acting on functions. In particular, if one takes two covariant derivatives of higher order tensors (for example taking two covariant derivatives of $\nabla f$), they need not commute. This leads to the definition of the Riemann curvature tensor $R$ which acts on tensors to measure the error incurred from commuting two covariant derivatives. We start with the case of a vector field $Z$, in this case the Riemann curvature tensor is defined as:

$$
R(X,Y)Z:= \nabla_X (\nabla_Y Z)-\nabla_Y (\nabla_X Z)-\nabla_{[X,Y]}Z
$$

or more simply using abstract index notation:

$$
R_{\alpha \beta \gamma}{}^\delta Z^\gamma=(\nabla_\alpha \nabla_\beta-\nabla_\beta \nabla_\alpha)Z^\delta. \footnote{The coordinate-free notation and the index notation are related by:
$$
\left(R(X,Y)Z\right) ^{\delta}=R_{\alpha \beta \gamma}{}^\delta X^\alpha Y^{\beta} Z^\gamma
$$

for any vector fields $X=X^\alpha$, $Y=Y^{\beta}$, and $Z=Z^{\gamma}$. }
$$

What will be important to us is that $R$ is actually a tensor, in the sense that $Z$ is not differentiated when computing $R(X,Y)Z$\footnote{This is in contrast to the fact that the commuting two pseudo-differential operators of order 1 gives a pseudo differential operator of order 1.}. 

The Riemann curvature tensor $R_{\alpha \beta \gamma}^\delta$ does not only give the commutation rules for covariant derivatives acting on vectors, but also for those acting on general $(k,l)$ tensors $T=T^{\gamma_1 \ldots \gamma_l}_{\delta_1\ldots \delta_k}$. We first start with 1-forms: when $T=\omega_\gamma$ is a 1-form, since $\nabla g=0$, we have: 

$$
(\nabla_\alpha \nabla_\beta  -\nabla_\beta \nabla_\alpha)\omega_\gamma= -R_{\alpha \beta \gamma}^\delta \omega_\delta
$$

More generally we have for a $(k,l)$ tensor $T=T_{\delta_1\ldots \delta_k}^{\gamma_1 \ldots \gamma_l}$:

\begin{equation}\label{commuting covariant derivatives on tensors}
(\nabla_\alpha \nabla_\beta-\nabla_\beta \nabla_\alpha)T_{\delta_1\ldots \delta_k}^{\gamma_1 \ldots \gamma_l}=\sum_{i=1}^l R_{\alpha \beta m}^{\gamma_i}T_{\delta_1\ldots \delta_k}^{\gamma_1 \ldots \gamma_{i-1} m \gamma_{i+1} \ldots \gamma_l} 
-\sum_{j=1}^k R_{\alpha \beta \delta_j}^m T_{\delta_1 \ldots \delta_{j-1}m\delta_{j+1}\ldots \delta_k}^{\gamma_1 \ldots \gamma_l}
\end{equation}

(see for example the first chapter of \cite{CLN}). As an application of the above Ricci commutator identities we have:

\textbf{Commutator of $\Delta$ and $\nabla$:}

Recall that the (rough)Laplace operator (also known as the Bochner Laplacian\footnote{This is in contrast with the closely related Hodge Laplacian which will not be of concern for us.}) is defined on general tensors as:

$$
\Delta= \operatorname{div} \nabla=\operatorname{trace}_g \nabla^2=\nabla_\alpha \nabla^{\alpha}=g^{\alpha \beta}\nabla_\alpha\nabla_\beta.
$$

\begin{lemma}\label{commuting Delta and nabla}

\begin{enumerate}

\item For any function $f$
$$
(\Delta \nabla_\alpha-\nabla_\alpha \Delta)f=\operatorname{Ric}_{\alpha \beta}\nabla^{\beta}f
$$
where $\operatorname{Ric}_{\alpha \beta}:= R_{\gamma \alpha \beta}{}^\gamma$ is the Ricci tensor.

\item For any tensor $T$

$$
\Delta \nabla T -\nabla \Delta T= \mathcal{O}\left(R*\nabla T\right) +\mathcal{O}\left((\nabla \operatorname{Ric})*T\right)
$$

where, given any tensors $A$ and $B$, $A*B$ denotes some contraction of $A\otimes B$ and $\mathcal{O}(A*B)$ denotes a linear combination of contractions of $A\otimes B$.

\end{enumerate}
\end{lemma}

\proof We present the simple proof as a warm up for the calculations to come. The proof of (i) follows easily from:

$$
\Delta \nabla_\alpha f=\nabla_{\beta}\nabla^\beta \nabla_\alpha f=\nabla_\beta \nabla_\alpha \nabla^\beta f=\nabla_\alpha \nabla_\beta \nabla^\beta f +R_{\beta \alpha \gamma}{}^\beta \nabla^\gamma f
=\nabla_\alpha \Delta f+\operatorname{Ric}_{\alpha \beta}\nabla^\beta f
$$

where in the second equality we used that covariant derivatives commute on functions, while in the third and forth we used the definition of $R$ and $\operatorname{Ric}$ respectively. 

The proof of (ii) follows similarly using $\eqref{commuting covariant derivatives on tensors}$ (see the first chapter of \cite{CLN} for a more comprehensive review).

\endproof

\subsection{Proof of Theorem \ref{A_0 and A_n theorem}} We are now ready to present the calculation leading to $\eqref{A_0 and A_n}$. For $i=1,\ldots, 4$, let $e_i$ be an eigenfunction of the Laplace-Beltrami operator with Dirichlet or Neumann boundary conditions corresponding to eigenvalues $-n_i^2$. Denote:

$$
A_0=\int_M e_1(x)\ldots e_4(x) dx.
$$

Then by Green's theorem\footnote{Recall that Green's theorem states that $\int_M(u\Delta v-v\Delta u)dx=\int_{\partial M}(u \frac{\partial v}{\partial n}-v\frac{\partial u}{\partial n})dS$ where $\partial/\partial n$ denotes the normal derivative on the boundary and $dS$ is the induced measure on $\partial M$. Since we are either assuming Dirichlet or Neumann boundary conditions, all boundary integrals vanish.}
$$
n_1^2 A_0=\int_M (-\Delta e_1)e_2 e_3 e_4 dx=\int_{M}e_1(-\Delta)[e_2 e_3 e_4]dx.
$$ 

But

\begin{align*}
\int_{M}e_1(-\Delta)[e_2 e_3 e_4]dx=&\int_M e_1(x)\left((-\Delta e_2)e_3 e_4+ e_2(-\Delta e_3)e_4+e_2 e_3 (-\Delta e_4)\right)dx\\
&-2\int_M e_1\left(\nabla_\alpha e_2 \nabla^\alpha e_3 e_4 +\nabla_\alpha e_2 e_3 \nabla^\alpha e_4 +e_2 \nabla_\alpha e_3 \nabla^\alpha e_4\right)dx \\
=& (n_2^2+n_3^2+n_4^2)A_0 -2\int_M e_1\left(\nabla_\alpha e_2 \nabla^\alpha e_3 e_4 +\nabla_\alpha e_2 e_3 \nabla^\alpha e_4 +e_2 \nabla_\alpha e_3 \nabla^\alpha e_4\right)dx .
\end{align*}

As a result, we get:
$$
A_0=\frac{-2A_1}{n_1^2 -n_2^2-n_3^2-n_4^2}
$$

where

$$
A_1=\int_M e_1\left(\nabla_\alpha e_2 \nabla^\alpha e_3 e_4 +\nabla_\alpha e_2 e_3 \nabla^\alpha e_4 +e_2 \nabla_\alpha e_3 \nabla^\alpha e_4\right)dx .
$$

Now we repeat the same procedure for $A_1$:

\begin{align*}
n_1^2 A_1=&\int_M (-\Delta e_1)\left(\nabla_\alpha e_2 \nabla^\alpha e_3 e_4 +\nabla_\alpha e_2 e_3 \nabla^\alpha e_4 +e_2 \nabla_\alpha e_3 \nabla^\alpha e_4\right)dx \\
=&\int_M e_1(-\Delta)\left[\nabla_\alpha e_2 \nabla^\alpha e_3 e_4 +\nabla_\alpha e_2 e_3 \nabla^\alpha e_4 +e_2 \nabla_\alpha e_3 \nabla^\alpha e_4\right]dx 
\end{align*}

by Green's formula. We now compute $(-\Delta)(\nabla_\alpha e_2 \nabla^\alpha e_3 e_4)$:

\begin{equation}\label{first term in second iteration}
\begin{split}
(-\Delta)(\nabla_\alpha e_2 \nabla^\alpha e_3 e_4)=&(-\Delta \nabla_\alpha e_2)\nabla^\alpha e_3 e_4+\nabla_\alpha e_2 (-\Delta \nabla^\alpha e_3)e_4 +\nabla_\alpha e_2 \nabla^\alpha e_3 (-\Delta e_4)\\
&-2\left(\nabla_{\alpha_1}\nabla_{\alpha_0}e_2\nabla^{\alpha_1}\nabla^{\alpha_0}e_3 e_4 +\nabla_{\alpha_1}\nabla_{\alpha_0}e_2 \nabla^{\alpha_0}e_3 \nabla^{\alpha_1}e_4+\nabla_{\alpha_0}e_2 \nabla_{\alpha_1}\nabla^{\alpha_0}e_3 \nabla^{\alpha_1}e_4\right)\\
=&(-\nabla_\alpha \Delta e_2)\nabla^\alpha e_3 e_4+\nabla_\alpha e_2 (-\nabla^\alpha \Delta e_3)e_4 +\nabla_\alpha e_2 \nabla^\alpha e_3 (-\Delta e_4)\\
&-\operatorname{Ric}_{\alpha \beta}\nabla^{\beta}e_2 \nabla^\alpha e_3 e_4- \operatorname{Ric}_{\alpha \beta}\nabla^\alpha e_2 \nabla^\beta e_3 e_4\\
&-2\left(\nabla_{\alpha_1}\nabla_{\alpha_0}e_2\nabla^{\alpha_1}\nabla^{\alpha_0}e_3 e_4 +\nabla_{\alpha_1}\nabla_{\alpha_0}e_2 \nabla^{\alpha_0}e_3 \nabla^{\alpha_1}e_4+\nabla_{\alpha_0}e_2 \nabla_{\alpha_1}\nabla^{\alpha_0}e_3 \nabla^{\alpha_1}e_4\right)\\
=&(n_2^2+n_3^2+n_4^2)\left(\nabla_\alpha e_2\nabla^\alpha e_3 e_4\right)-2\operatorname{Ric}_{\alpha \beta}\nabla^{\beta}e_2 \nabla^\alpha e_3 e_4\\
&-2\left(\nabla_{\alpha_1}\nabla_{\alpha_0}e_2\nabla^{\alpha_1}\nabla^{\alpha_0}e_3 e_4 +\nabla_{\alpha_1}\nabla_{\alpha_0}e_2 \nabla^{\alpha_0}e_3 \nabla^{\alpha_1}e_4+\nabla_{\alpha_0}e_2 \nabla_{\alpha_1}\nabla^{\alpha_0}e_3 \nabla^{\alpha_1}e_4\right)\\
\end{split}
\end{equation}

where we have used lemma \ref{commuting Delta and nabla} for the second inequality. Similarly, one gets:

\begin{equation}\label{second term in second iteration}
\begin{split}
(-\Delta)(\nabla_\alpha e_2 e_3 \nabla^\alpha e_4)=&(n_2^2+n_3^2+n_4^2)\left(\nabla_\alpha e_2 e_3\nabla^\alpha e_4\right)-2\operatorname{Ric}_{\alpha \beta}\nabla^{\beta}e_2 e_3 \nabla^\alpha e_4\\
&-2\left(\nabla_{\alpha_1}\nabla_{\alpha_0}e_2\nabla^{\alpha_1}e_3 \nabla^{\alpha_0} e_4 +\nabla_{\alpha_1}\nabla_{\alpha_0}e_2 e_3 \nabla^{\alpha_1}\nabla^{\alpha_0} e_4+\nabla_{\alpha_0}e_2 \nabla_{\alpha_1}e_3 \nabla^{\alpha_1}\nabla^{\alpha_0} e_4\right)\\
\end{split}
\end{equation}

and

\begin{equation}\label{third term in second iteration}
\begin{split}
(-\Delta)(e_2 \nabla_\alpha e_3 \nabla^\alpha e_4)=&(n_2^2+n_3^2+n_4^2)\left(e_2 \nabla_\alpha e_3\nabla^\alpha e_4\right)-2\operatorname{Ric}_{\alpha \beta}e_2 \nabla^\beta e_3 \nabla^\alpha e_4\\
&-2\left(\nabla^{\alpha_1}e_2\nabla_{\alpha_1}\nabla_{\alpha_0} e_3 \nabla^{\alpha_0} e_4 +\nabla_{\alpha_1}e_2 \nabla_{\alpha_0} e_3 \nabla^{\alpha_1}\nabla^{\alpha_0} e_4+e_2 \nabla_{\alpha_1}\nabla_{\alpha_0}e_3 \nabla^{\alpha_1}\nabla^{\alpha_0} e_4\right).\\
\end{split}
\end{equation}

Adding $\eqref{first term in second iteration},\eqref{second term in second iteration},\eqref{third term in second iteration}$, we get:

$$
n_1^2 A_1=(n_2^2+n_3^2+n_4^2)A_1 -2A_2
$$

where $A_2$ is of the following form:

$$
A_2:=\int_M e_1\left(B_2(e_2,e_3, e_4)+C_2(e_2,e_3,e_4)\right)dx
$$

where $B_2(f,g,h)$ and $C_2(f,g,h)$ are \emph{trilinear} operators that can be expressed schematically as:

\begin{equation}\label{def of B_2}
B_2(f,g,h)=\mathcal{O}_{\substack{i+j+k=4 \\ 0\leq i,j,k \leq 2}}\left(\nabla^i f * \nabla^j g * \nabla^k h \right)
\end{equation}

\begin{equation}\label{def of C_2}
C_2(f,g,h)=\mathcal{O}_{\substack{i+j+k\leq 2 \\ 0\leq i,j,k \leq 1}}\left(R*\nabla^i f * \nabla^j g * \nabla^k h \right).
\end{equation}

Now suppose, as an induction hypothesis, that:

$$
A_{n-1}=\frac{-2}{n_1^2-n_2^2-n_3^2-n_4^2} A_n
$$

where

$$
A_n=\int_{M}e_1\left(B_n(e_2,e_3,e_4)+C_n(e_2,e_3,e_4)\right)dx
$$

and 

\begin{equation}\label{def of B_n 2}
B_n(f,g,h)=\mathcal{O}_{\substack{i+j+k=2n \\ 0\leq i,j,k \leq n}}\left(\nabla^i f * \nabla^j g * \nabla^k h \right)
\end{equation}

\begin{equation}\label{def of C_n 2}
C_n(f,g,h)=\mathcal{O}_{\substack{i+j+k\leq 2(n-1) \\ 0\leq i,j,k \leq n-1}}\left(\nabla^aR*\nabla^i f * \nabla^j g * \nabla^k h \right)
\end{equation}

where $a$ is some exponent(which can be calculated explicitly in terms of $n$) signifying a number of derivatives applied to the curvature tensor.

Then, as before, by Green's theorem:

$$
n_1^2 A_n=\int_{M}(-\Delta)e_1\left(B_n(e_2,e_3,e_4)+C_n(e_2,e_3,e_4)\right)dx=\int_{M}e_1\left((-\Delta)B_n(e_2,e_3,e_4)+(-\Delta)C_n(e_2,e_3,e_4)\right)dx
$$

and our goal is to write this as 

$$
n_1^2 A_n=(n_2^2+n_3^2+n_4^2)A_n -2 \int_M e_1\left(B_{n+1}(e_2,e_3,e_4)+C_{n+1}(e_2,e_3,e_4)\right)dx
$$

with $B_{n+1}, C_{n+1}$ as in $\eqref{def of B_n 2}$ and $\eqref{def of C_n 2}$ with $n$ replaced by $n+1$. For this we use the rules to commute $\Delta$ and $\nabla$ in lemma \ref{commuting Delta and nabla}(ii). In fact,

\begin{align*}
(-\Delta) B_n(e_2,e_3,e_4)= \mathcal{O}_{\substack{i+j+k=2n \\ 0\leq i,j,k \leq n}}(-\Delta)\left(\nabla^i e_2 * \nabla^j e_3 * \nabla^k e_4 \right)\\
=\mathcal{O}_{\substack{i+j+k=2n \\ 0\leq i,j,k \leq n}}\bigg\{\left((-\Delta)\nabla^i e_2 * \nabla^j e_3 * \nabla^k e_4 +\nabla^i e_2 * (-\Delta)\nabla^j e_3 * \nabla^k e_4 +\nabla^i e_2 \nabla^j e_3 (-\Delta)\nabla^k e_4\right)\\
+2\left(\nabla_\alpha \nabla^i e_2* \nabla^\alpha \nabla^j e_3*\nabla^k e_4 +\nabla_\alpha\nabla^i e_2* \nabla^j e_3 *\nabla^\alpha \nabla^k e_4+\nabla^i e_2* \nabla_\alpha \nabla^j e_3* \nabla^\alpha \nabla^k e_4\right)\bigg\}\\
=\mathcal{O}_{\substack{i+j+k=2n \\ 0\leq i,j,k \leq n}}\bigg\{\left(\nabla^i (-\Delta e_2) * \nabla^j e_3 * \nabla^k e_4 +\nabla^i e_2 * \nabla^j (-\Delta e_3) * \nabla^k e_4 +\nabla^i e_2 \nabla^j e_3 \nabla^k (-\Delta e_4)\right)\\
+2\left(\nabla_\alpha \nabla^i e_2 *\nabla^\alpha \nabla^j e_3 *\nabla^k e_4+\nabla_\alpha\nabla^i e_2 *\nabla^j e_3* \nabla^\alpha \nabla^k e_4+\nabla^i e_2 *\nabla_\alpha \nabla^j e_3 *\nabla^\alpha \nabla^k e_4\right)\bigg\}\\
+\mathcal{O}_{\substack{i+j+k
\leq 2n \\ 0\leq i,j,k \leq n}}\left(\nabla^a R*\nabla^i e_2 * \nabla^j e_3 * \nabla^k e_4 \right).
\end{align*}

We should note that for the third equality above, we used lemma \ref{commuting Delta and nabla}(ii) inductively applied to the tensor $\nabla^{i-1}e_2$ (and similarly for $\nabla^{j-1}e_3$ and $\nabla^{k-1}e_4$) to give $\Delta \nabla (\nabla^{i-1}e_2)-\nabla \Delta (\nabla^{i-1} e_2)=\mathcal{O}(R*\nabla^i e_2 +\nabla \operatorname{Ric}*\nabla^{i-1}e_2)$ which is of the form above. Doing the same thing for $\Delta \nabla^{i-1} e_2$ over and over we get the result claimed. As a result, we get:

\begin{equation}\label{Delta B_n}
\begin{split}
(-\Delta) B_n(e_2,e_3,e_4)
=&(n_2^2+n_3^2+n_4^2) B_n(e_2,e_3,e_4)\\
&+\mathcal{O}_{\substack{i+j+k=2n \\ 0\leq i,j,k \leq n}}\bigg(\nabla_\alpha \nabla^i e_2* \nabla^\alpha \nabla^j e_3 *\nabla^k e_4+\nabla_\alpha\nabla^i e_2* \nabla^j e_3* \nabla^\alpha \nabla^k e_4\\
&+\nabla^i e_2* \nabla_\alpha \nabla^j e_3 *\nabla^\alpha \nabla^k e_4\bigg)+\mathcal{O}_{\substack{i+j+k\leq 2n \\ 0\leq i,j,k \leq n}}\left(\nabla^a R*\nabla^i e_2 * \nabla^j e_3 * \nabla^k e_4 \right)\\
=&(n_2^2+n_3^2+n_4^2) B_n(e_2,e_3,e_4)
+\mathcal{O}_{\substack{i+j+k=2(n+1) \\ 0\leq i,j,k \leq n+1}}\left(\nabla^i e_2 *\nabla^j e_3 *\nabla^k e_4\right)\\
&+\mathcal{O}_{\substack{i+j+k
\leq 2n \\ 0\leq i,j,k \leq n}}\left(\nabla^a R*\nabla^i e_2 * \nabla^j e_3 * \nabla^k e_4 \right).\\
\end{split}
\end{equation}

Obviously, the second term above will join $B_{n+1}$ whereas the third will be part of $C_{n+1}$. The computation for $\Delta C_n$ is similar:

\begin{align*}
(-\Delta)C_n(e_2,e_3,e_4)=&\mathcal{O}_{\substack{i+j+k\leq 2(n-1) \\ 0\leq i,j,k \leq n-1}}(-\Delta)\left(\nabla^aR*\nabla^i e_2 * \nabla^j e_3 * \nabla^k e_4 \right)\\
=& \mathcal{O}_{\substack{i+j+k\leq 2(n-1) \\ 0\leq i,j,k \leq n-1}}\bigg(
\nabla^{a}R*(-\Delta)\nabla^i e_2*\nabla^{j}e_3*\nabla^k e_4+ \nabla^{a}R*\nabla^i e_2*(-\Delta)\nabla^{j}e_3*\nabla^k e_4\\
&+ \nabla^{a}R*\nabla^i e_2*\nabla^{j}e_3*(-\Delta)\nabla^k e_4\bigg)
+\mathcal{O}_{\substack{i+j+k\leq 2n \\ 0\leq i,j,k \leq n}}\left(\nabla^{a'} R *\nabla^i e_2 *\nabla^j e_3* \nabla^k e_4\right)\\
=& \mathcal{O}_{\substack{i+j+k\leq 2(n-1) \\ 0\leq i,j,k \leq n-1}}\bigg(
\nabla^{a}R*\nabla^i (-\Delta e_2)*\nabla^{j}e_3*\nabla^k e_4+ \nabla^{a}R*\nabla^i e_2*\nabla^{j}(-\Delta e_3)*\nabla^k e_4\\
&+ \nabla^{a}R*\nabla^i e_2*\nabla^{j}e_3*\nabla^k  (-\Delta e_4)\bigg)
+\mathcal{O}_{\substack{i+j+k\leq 2n \\ 0\leq i,j,k \leq n}}\left(\nabla^a R *\nabla^i e_2 *\nabla^j e_3 *\nabla^k e_4\right)
\end{align*}

by applying lemma \ref{commuting Delta and nabla} inductively as before. As a result, we get:

\begin{equation}\label{Delta C_n}
(-\Delta)C_n(e_2,e_3,e_4)= (n_2^2+n_3^2+n_4^2)C_n(e_2,e_3,e_4)
+\mathcal{O}_{\substack{i+j+k\leq 2n \\ 0\leq i,j,k \leq n}}\left(\nabla^a R *\nabla^i e_2 *\nabla^j e_3* \nabla^k e_4\right).
\end{equation}

The last term will join $C_{n+1}(e_2,e_3,e_4)$ to give that:

$$
n_1^2 A_n =(n_2^2+n_3^2+n_4^2)A_n +\int_M e_1 \left(B_{n+1}(e_2,e_3,e_4)+C_{n+1}(e_2,e_3,e_4)\right)dx
$$

where $B_{n+1}$ and $C_{n+1}$ are trilinear operators as in $\eqref{def of B_n}$ and $\eqref{def of C_n}$ respectively, with $n$ replaced by ${n+1}$. This concludes the induction proof.
\endproof

\subsection{Proof of Corollary \ref{spectral cluster lemma 1}}: The proof of Corollary \ref{spectral cluster lemma 1} will follow by applying Theorem \ref{A_0 and A_n theorem} and a variant of the bilinear eigenfunction cluster estimates of \cite{BEE}\cite{MEE} after dealing with a couple of technical problems. Let us first recall the bilinear eigenfunction cluster estimates (also called bilinear Sogge estimates) from \cite{MEE}. 

\begin{proposition}\label{BGT3 bilinear estimates}(Bilinear eigenfunction cluster estimates \cite{MEE})

Let $\chi \in \mathcal{S}(\R)$. For $\lambda \in \R$, denote by $\chi_\lambda=\chi(\sqrt{-\Delta}-\lambda)$ the spectral projector around $\lambda$. For any $\lambda \geq \mu$,

\begin{equation}\label{BGT3 bilinear estimates}
||\chi_\lambda f \chi_\mu g||_{L^2(M)} \lesssim \Lambda(d,\mu) ||f||_{L^2(M)}||g||_{L^2(M)}
\end{equation}

for all $f,g \in L^2(M)$ where $\Lambda(d,\mu)$ was defined in $\eqref{def of Lambda(d,mu)}$.
\end{proposition}

The proof in \cite{MEE} is based on Sogge's parametrix representation of $\chi_\lambda f$ in local coordinates (see \cite{Sogge}). More precisely, for every $N \geq 1$, one has the splitting:

$$
\chi_\lambda f= \lambda^{\frac{d-1}{2}}T_\lambda f +R_\lambda f
$$

with 
$$
||R_\lambda f||_{H^k(M)}\lesssim_{N,k}\lambda^{k-N}||f||_{L^2(M)}
$$

and in a system of local coordinates around each $x_0\in M$, $T_\lambda$ has the following parametrix representation:

\begin{equation}\label{eigenfunction parametrix}
T_\lambda f(x)=\int_{\R^d} e^{i\lambda \phi(x,y)}a(x,y,\lambda)f(y) dy
\end{equation}

where $a(x,y,\lambda)$ is a polynomial in $\lambda^{-1}$ with smooth coefficients supported in a compact subset $\{(x,y) \in V\times V: V \subset \R^d \textrm{(compact)}\}$ and $-\phi(x,y)=d_g(x,y)$ is the geodesic distance between $x$ and $y$. By taking $N$ large enough, the proof of $\eqref{BGT3 bilinear estimates}$ reduces to proving estimates for oscillatory integrals of the form $\eqref{eigenfunction parametrix}$ (cf. \cite{MEE}]).

As a result of this, one can directly deduce by the exact analysis as that leading to $\eqref{BGT3 bilinear estimates}$ that the following estimate holds:

\begin{equation}\label{BGT3 bilinear estimates in tensor form}
||(\nabla^j \chi_\lambda f)*(\nabla^k \chi_\mu g)||_{L^2(M)} \lesssim \lambda^j \mu^k\Lambda(d,\mu) ||f||_{L^2(M)}||g||_{L^2(M)}
\end{equation}

since any contraction of the form $(\nabla^j \chi_\lambda f)*(\nabla^k \chi_{\mu}g)$ can be written using the above splitting of $\chi_\lambda f = \lambda^{\frac{d-1}{2}}T_\lambda f +R_\lambda f$ and $\chi_\mu g =\mu^{\frac{d-1}{2}}T_\mu g +R_\mu$ as a linear combination of products of operators of the same form for which the analysis in \cite{MEE} applies (what is crucial in this reduction is that $R$ is a smoothing operator for large enough $N$ and so it does not affect the analysis in any significant way. Also, $a$ is just a symbol of order 0 in $x$ and is a polynomial in $\frac{1}{\lambda}$, so any derivative falling on it does not affect the analysis in any way either).

We now turn to the proof of proposition \ref{spectral cluster lemma 1}. Without loss of generality, we assume that $f=\mathbf{1}_{[\lambda,\lambda+1]}(\sqrt{-\Delta})f$, $g=\mathbf{1}_{[\mu,\mu+1]}(\sqrt{-\Delta})g$, and $h=\mathbf{1}_{[\nu,\nu+1]}(\sqrt{-\Delta})h$. Write

$$
\int_M h(x) f(x) g(x)dx= \sum_{\substack{n_1 \in [\nu,\nu+1]\\ n_2 \in [\lambda,\lambda+1]\\n_3 \in [\mu, \mu+1]}} \int_M \pi_{n_1}h\, \pi_{n_2}f\, \pi_{n_3}g\,dx
$$

where $\pi_{n}$ is the projection onto the $n^2$-eigenspace. By lemma \ref{A_0 and A_n} (with $e_4=1$)we get:

\begin{equation}\label{int hfg}
\int_M h(x) f(x) g(x)dx= \sum_{\substack{n_1 \in [\nu,\nu+1]\\ n_2 \in [\lambda,\lambda+1]\\n_3 \in [\mu, \mu+1]}} \frac{(-2)^J}{\left(n_1^2-n_2^2-n_3^2\right)^J} \int_M \pi_{n_1} h \left(B(\pi_{n_2} f, \pi_{n_3} g) +C(\pi_{n_2}f,\pi_{n_3}g\right)dx
\end{equation}

where 

$$
B(\pi_{n_2} f,\pi_{n_3}g)=\left(\nabla^{J}\pi_{n_2} f\right)*\left(\nabla^{J} \pi_{n_3} f\right)
$$

and 

$$
C(\pi_{n_2}f,\pi_{n_3}g)=\mathcal{O}_{\substack{i+j\leq 2(J-1) \\ 0\leq i,j \leq J-1}}\left(\nabla^aR*(\nabla^i \pi_{n_2}f) * (\nabla^j \pi_{n_3}g)\right).
$$

The integral in $\eqref{int hfg}$ is similar to the ones studied in section \ref{multilinear spectral multiplier lemma}. Essentially, one would like to argue as follows: since $\nu=\lambda +K\mu+2$, then $n_1^2-n_2^2-n_3^2 \geq 2K\lambda \mu$ and hence one would like to estimate $\eqref{int hfg}$ as follows:

$$
\text{LHS of $\eqref{int hfg}$} \lesssim_J \frac{1}{(K\lambda \mu)^J}\left|\int_{M}  h(x)(B(f,g)+C(f,g))dx\right| \lesssim_J \frac{\lambda^J\mu^J\Lambda(d,\mu)}{(K\lambda \mu)^J}||h||_{L^2(M)}||f||_{L^2(M)}||g||_{L^2(M)}
$$
where the first inequality will be justified in what follows whereas the second inequality comes from $\eqref{BGT3 bilinear estimates in tensor form}$ and the fact that all derivatives of $R$ are bounded to get the second inequality. To justify the first inequality, we argue similar to what we did in section \ref{multilinear spectral multiplier lemma}. We include the details for completeness. As before, we write $n_1=\nu+r_1$, $n_2=\lambda+r_2$, and $n_3=\mu+r_3$ with $(r_1, r_2, r_3)\in [0,1]^3$. As a result, we get that $
\frac{2^J}{(n_1^2-n_2^2-n_3^2)^J}=\frac{2^J}{\left((\nu+r_1)^2-(\lambda+r_2)^2-(\mu+r_3)^2\right)^J}$. But 
$$
\Xi(x_1,x_2,x_3):=\frac{1}{\left((\nu+x_1)^2-(\lambda+x_2)^2-(\mu+x_3)^2\right)^J}
$$
is a smooth function on $[0,1]^3$ that is bounded along with its derivatives by $\frac{1}{(K\lambda \mu)^J}$. In fact, this follows from the following estimates:

$$
\frac{2^J}{(n_1^2-n_3^2-n_4^2)^J}\leq \frac{2^J}{\left(\nu^2-(\lambda+1)^2-(\mu+1)^2\right)^J}=\frac{2^J}{\left((\lambda+K\mu+2)^2-(\lambda+1)^2+(\mu+1)^2\right)^J}\leq \frac{1}{(K\lambda \mu)^J}
$$

since $K>1$. A similar estimate holds for the first two derivatives. Multiplying $f(x_1,x_2,x_3)$ by a compactly supported function on $[-2,2]^3$ and extending the resulting function periodically to $\R^3$ gives a $4-$periodic function on $\R^3$. Expressing this function in Fourier series we get that $\Xi(x_1,x_2,x_3)=\sum_{\theta_i \in \Z/4}A(\theta_1, \theta_2, \theta_3)e^{i(\theta_1 x_1+\theta_2 x_2 +\theta_3 x_3)}$ with 
\begin{equation}\label{l^1 estimate last}
\sum_{\theta_i \in \Z/4}|A(\theta_1,\theta_2,\theta_3)|\lesssim \frac{1}{(K\lambda \mu)^J}
\end{equation} 
since $A(\theta_1,\theta_2,\theta_3)$ are the Fourier coefficients of a $C^2$ function whose $C^2$ norm is bounded by $\frac{1}{(K\lambda \mu)^J}$.

With this in hand, we write:

\begin{align*}
\int_M h(x)f(x)g(x) dx=& \sum_{\substack{r_1 \in \mathcal{S}-\nu \cap [0,1] \\ r_2 \in S-\lambda \cap [0,1]\\ r_3\in \mathcal{S}-\mu \cap [0,1]}} \sum_{\theta_i \in \Z/4}A(\theta_1,\theta_2,\theta_3)e^{i(\theta_1 r_1+\theta_2 r_2+\theta_3 r_3)}\\ 
&\times \int_M \pi_{\nu+r_1} h \left(B(\pi_{\lambda+r_2} f, \pi_{\mu+r_3} g) +C(\pi_{\lambda+r_2}f,\pi_{\mu+r_3}g\right)dx
\end{align*}

where we used $\mathcal{S}$ to denote the set $\{n \in \R: n^2 \in \text{spectrum of}(-\Delta_g)\}$ and $\mathcal{S}-\alpha=\{n-\alpha:n\in \mathcal{S}\}$. Letting 
$\tilde h_{\theta_1} =\sum_{r_1 \in \mathcal{S}-\nu \cap [0,1]}e^{i\theta_1 r_1}\pi_{\nu+r_1} h$, $\tilde f_{\theta_2} =\sum_{r_2 \in \mathcal{S}-\lambda \cap [0,1]}e^{i\theta_2 r_2}\pi_{\lambda+r_2} f$, and $\tilde g_{\theta_3} =\sum_{r_3 \in \mathcal{S}-\mu \cap [0,1]}e^{i\theta_3 r_3}\pi_{\mu+r_3} g$, we get, using the fact that $B$ and $C$ are multilinear, that:

\begin{align*}
\int_M h(x)f(x)g(x) dx=& \sum_{\theta_i \in \Z/4}A(\theta_1,\theta_2,\theta_3)\int_M \tilde h_{\theta_1} \left(B(\tilde f_{\theta_2}, \tilde g _{\theta_3}) +C(\tilde f_{\theta_2},\tilde g_{\theta_3}\right)dx.
\end{align*}

But for each fixed $(\theta_1,\theta_2,\theta_3)\in \Z/4$, we have the estimate:

\begin{align*}
\left|\int_M \tilde h_{\theta_1} \left(B(\tilde f_{\theta_2}, \tilde g _{\theta_3}) +C(\tilde f_{\theta_2},\tilde g_{\theta_3}\right)dx\right| \leq& ||\tilde h_{\theta_1}||_{L^2(M)}||\left(B(\tilde f_{\theta_2}, \tilde g _{\theta_3}) +C(\tilde f_{\theta_2},\tilde g_{\theta_3}\right)||_{L^2(M)}\\
\lesssim& \lambda^J \mu^J \lambda(d,\mu) ||\tilde h_{\theta_1}||_{L^2(M)}||\tilde f_{\theta_2}||_{L^2(M)}||\tilde g_{\theta_3}||_{L^2(M)}\\
=&\lambda^J \mu^J \lambda(d,\mu) ||h||_{L^2(M)}||f||_{L^2(M)}|| g||_{L^2(M)}
\end{align*}

where we used Cauchy-Schwarz for the first inequality, $\eqref{BGT3 bilinear estimates in tensor form}$ in the second, and the fact that frequency modulation leave that $L^2$ norm invariant in the third.

As a result we get:

\begin{align*}
\left|\int_M h(x)f(x)g(x) dx\right|\lesssim&  \lambda^J \mu^J \lambda(d,\mu)\left(\sum_{\theta_i \in \Z/4}|A(\theta_1,\theta_2,\theta_3)|\right) ||h||_{L^2(M)}||f||_{L^2(M)}|| g||_{L^2(M)}\\
\lesssim& \frac{\lambda^J\mu^J \Lambda(d,\mu)}{(K\lambda \mu)^J}||h||_{L^2(M)}||f||_{L^2(M)}|| g||_{L^2(M)}=\frac{\Lambda(d,\mu)}{K^J}||h||_{L^2(M)}||f||_{L^2(M)}|| g||_{L^2(M)}
\end{align*}

as desired. \endproof

\end{document}